\documentclass{amsart}
\usepackage{amscd,amssymb}
\usepackage[all]{xy}
\xyoption{poly}
\xyoption{arc}

\def\A{\mathbb{A}}
\def\B{\mathbb{B}}
\def\C{\mathbb{C}}

\def\P{\mathbb{P}}
\def\Q{\mathbb{Q}}
\def\N{\mathbb{N}}
\def\R{\mathbb{R}}
\def\Z{\mathbb{Z}}
\def\fatone{{1\hskip-4pt 1}}
\def\I{{\rm I}}

\def\CA{\mathcal{A}}

\def\CD{\mathcal{D}}

\def\CI{\mathcal{I}}

\def\CN{\mathcal{N}}
\def\CP{\mathcal{P}}

\def\CT{\mathcal{T}}

\def\DP{\mathcal{D}P}
\def\Fx{F^\times}
\def\Mot{{\rm Mot}}

\def\id{{\rm id}}
\def\SS{{\rm S}}

\def\Or{{\rm Or}}
\def\H{{\rm{H}}}
\def\B{{\rm{B}}}
\def\rmc{{\rm{c}}}
\def\lceilx{{\,[\,}}
\def\rfloorx{{\,]\,}}
\def\Cycle{{\rm Cycle}}
\def\cycle{{\rm cycle}}

\def\cub{\square}
\def\disjoint{\sqcup}

\def \cal {\mathcal}
\def \smallbulletP{{\hskip -1pt\scriptscriptstyle \bullet\hskip 1pt}}
\def \smallbullet{{\hskip 0pt\scriptscriptstyle \bullet\hskip 0pt}}

\def \pg{{\rm pg}}

\def \barbarpg{{\barbar{\rm pg}}}
\def \tr{{\rm tr}}
\def \th{{\rm th}}

\def \eps{{\varepsilon}}

\def \sgn{{\rm sgn}}

\def \barbar{\overline}
\def \Bar{{\,\big\vert\,}}
\def \Barr{{\,\bigg\vert\,}}
\def \sha{{\,\amalg\hskip -3.6pt\amalg\,}}
\def \Sha#1#2{{\,\coprod\hskip -9.6pt\coprod_{#1\ }^{#2\ }\,}}

\def \Shaonly{{\,\coprod\hskip -6.2pt\coprod\,}}

\def \Li {{\mathbb L}{\rm i}}

\def \geq{{\,\geqslant\,}}
\def \leq{{\,\leqslant\,}}
\def \ms {{\medskip}}
\def \sm {{\smallskip}}

\def \om {\omega}

\def \Kriz {{Kriz}}
\def \dbar{{\overline{d}}}
\def \CTT {{\widetilde{\CT}}}

\def \fourgon#1#2#3#4{{
\xy
\POS(10,4) \ar@{=} +(-10,0)_#4
\ar@{-} +(0,-10)^#3
\POS(10,-6) \ar@{-} +(-10,0)^#2
\POS(0,4) \ar@{-} +(0,-10)_#1
\POS(0,4) *+{\bullet}
\endxy
}}

\def \fourgondotted#1#2#3#4{{
\xy
\POS(10,4) \ar@{=} +(-10,0)_#4
\ar@{-} +(0,-10)^#3
\POS(10,-6) \ar@{-} +(-10,0)^#2
\POS(0,4) \ar@{.} +(0,-10)_#1
\endxy
}}

\def \fourgona{{
\POS(10,5) \ar@{=} +(-10,0)_4
\ar@{-} +(0,-10)^3
\POS(10,-5) \ar@{-} +(-10,0)^2
\POS(0,5) \ar@{-} +(0,-10)_1
\POS(0,5) *+{\bullet}
}}

\def \threegon#1#2#3{{
\xy
\POS(10,4) \ar@{=} +(-10,0)_#3 
\ar@{-} +(-5,-8)^#2  
\POS(0,4) \ar@{-} +(5,-8)_#1
\POS(0,4) *+{\bullet}
\endxy
}}

\def \threegonarrowc#1#2#3{{
\xy  
\POS(10,4) \ar@{=} +(-10,0)_#3
\ar@{-} +(-5,-8)^#2
\POS(0,4) \ar@{-} +(5,-8)_#1
\POS(5,4) \ar@{<<-} +(0,-8)
\POS(0,4) *+{\bullet}
\endxy
}}

\def \threegonarrowb#1#2#3{{
\xy  
\POS(10,4) \ar@{=} +(-10,0)_#3
\ar@{-} +(-5,-8)^#2
\POS(0,4) \ar@{-} +(5,-8)_#1
\POS(7.5,0) \ar@{<<-} +(-7.5,4)
\POS(0,4) *+{\bullet}
\endxy
}}

\def \threegonarrowa#1#2#3{{
\xy  
\POS(10,4) \ar@{=} +(-10,0)_#3
\ar@{-} +(-5,-8)^#2
\POS(0,4) \ar@{-} +(5,-8)_#1
\POS(2.5,0) \ar@{<<-} +(7.5,4)
\POS(0,4) *+{\bullet}
\endxy
}}

\def \twogon#1#2{\hskip 5pt\xy
\POS(10,0) \ar@{=} +(-10,0)_#2
\POS(0,0)*{};
\POS(0,0) *+{\bullet}
\POS(10,0)*{};
**\crv{(5,-5)};
\POS(5,-5) *+{\scriptstyle #1}
\endxy \hskip 5pt}

\def \twogonright#1#2{\hskip 5pt\xy
\POS(10,0) \ar@{=} +(-10,0)_#2
\POS(0,0)*{};
\POS(10,0) *+{\bullet}
\POS(10,0)*{};
**\crv{(5,-5)};
\POS(5,-5) *+{\scriptstyle #1}
\endxy \hskip 5pt}

\newtheorem{lemma}{\bfseries Lemma}[section]
\newtheorem{Theorem}    {Theorem}[section]
\newtheorem{Corollary}  [Theorem]{Corollary}
\newtheorem{Proposition}       [Theorem]{Proposition}

\newtheorem{Lemma}      [Theorem]{Lemma}
\newtheorem{Definition} [Theorem]{Definition}

\newtheoremstyle{myremark}
 {3pt}
 {3pt}
 {\rm}
 {}
 {\bf}
 {:}
 {.5em}
 {}
\theoremstyle{myremark}
\newtheorem{Remark}[Theorem]{Remark}

\newtheorem{Example}    [Theorem]{Example}

\title{Multiple polylogarithms, polygons, trees and algebraic cycles}

\author{H.~Gangl, A.B.~Goncharov and A.~Levin}

\begin{document}

\begin{abstract}
We construct, for a field $F$ and a natural number $n$, algebraic cycles in 
Bloch's cubical cycle group of codimension $n$ cycles in 
$\big(\P^1_F\setminus\{1\}\big)^{2n-1}$, which correspond to weight $n$ 
multiple polylogarithms with generic arguments if $F\subset \C$.
Moreover, we construct out of them a Hopf subalgebra in the 
Bloch-Kriz cycle Hopf algebra $\chi_\cycle$. 
In the process, we are led to other Hopf algebras built from trees and 
 polygons, which are mapped to $\chi_\cycle$. 
We relate the coproducts to the one
for Goncharov's motivic multiple polylogarithms and to the Connes-Kreimer coproduct on plane trees
and produce the associated Hodge realization for polygons.
\end{abstract}
\maketitle
\tableofcontents

\section{Introduction}

\subsection{Multiple polylogarithms}\label{mulpol}
We start from the analytic background. 
The multiple polylogarithm \index{polylogarithm} \index{multiple polylogarithm} functions were defined in \cite{GonICM} by the power series
$$
Li_{n_1,\dots,n_m}(z_1,\dots,z_m)= \sum_{0<k_1<\cdots<k_m}
\frac{z_1^{k_1}}{k_1^{n_1}} \frac{z_2^{k_2}}{k_2^{n_2}}\dots
\frac{z_m^{k_m}}{k_{m\phantom{|}}^{n_m}}\qquad (z_i \in \C, |z_i|<1)\,.
$$
They admit an analytic continuation to a Zariski open subset of $\C^m$. 
When $m=1$, we recover
 the classical polylogarithm function. 
Putting $n_1 = ... = n_m =1$, we get 
the {\it multiple logarithm} \index{multiple logarithm} function. 

Let $x_i$ be complex numbers. Recall that an 
iterated integral \index{iterated integral} is defined as 
\begin{equation}\label{int_def}
I(x_0; x_1,\dots,x_m; x_{m+1}) = \int\limits_{\Delta_\gamma} \frac{dt_1}{t_1-x_1}\wedge 
\cdots \wedge\frac{dt_m}{t_m-x_m}\,,
\end{equation}
where $\gamma$ is a path from $x_0$ to $x_{m+1}$ in $\C - \{x_1, ..., x_m\}$, and the cycle of integration  $\Delta_\gamma$ 
consists of all $m$-tuples of points $(\gamma(t_1), ..., \gamma(t_m))$ with $t_i\leq t_j$ for $i<j$.

Multiple polylogarithms can be written as iterated integrals (cf.~loc.cit.). 
In particular, here is the
 iterated integral presentation of the multiple logarithm function:
\begin{equation} \label{ml}
Li_{1,\dots,1} (z_1,\dots,z_m)= (-1)^m I(0; x_1,\dots,x_m; 1) \,,
\end{equation}
 where we set
\begin{equation} \label{ml2}
x_1:= (z_1\cdots z_m)^{-1}, \quad x_2:= (z_2\cdots z_m)^{-1},\quad \dots\quad, x_m:= z_m^{-1}\,.
\end{equation}

Observe that in (\ref{ml2}) the parameters $x_1, \dots, x_m$ are non-zero. 
Moreover, many properties of the iterated  integrals will 
change if we put some of the $x_i$'s equal to zero. As a result, 
many features of the theory of multiple polylogarithms are different from the ones 
of multiple logarithms. We will use the notation $I_{1,\dots,1}(x_1, \dots, x_m)$ for
$I(0; x_1,\dots,x_m; 1)$ when the $x_i$'s are non-zero.

\subsection{From iterated integrals to Hopf algebras} 
In \cite{GonFund} it was shown that any iterated integral (\ref{int_def}) 
gives rise to an element 
\begin{equation}\label{int_defa}
I^{\cal H}(x_0; x_1,\dots,x_m; x_{m+1}) \in {\cal A}^{\cal H}_{m}, \qquad x_i \in \C,
\end{equation}
of a certain commutative graded  Hopf algebra ${\cal A}^{\cal H}_{\bullet}$, called 
the {\it fundamental Hodge-Tate Hopf algebra}. 
The category of graded finite-dimensional comodules over this Hopf algebra is canonically equivalent to the category of mixed $\Q$-Hodge-Tate structures. 
The element 
(\ref{int_defa}) is called a {\it Hodge iterated integral}.

It was shown in loc. cit. that this construction has 
$l$-adic and motivic counterparts. In particular, given a number field 
$F$, there are  the {\it motivic iterated integrals} 
\begin{equation}\label{int_defm}
I^{\cal M}(x_0; x_1,\dots,x_m; x_{m+1}) \in {\cal A}^{\cal M}_{m}, \qquad x_i \in F,
\end{equation}
which live in the {\it fundamental motivic Tate Hopf algebra} 
$ {\cal A}^{\cal M}_{\bullet}(F)$ of $F$. 

The benefits gained by  working 
with the Hodge/motivic iterated integrals are 
explained in the introduction to \cite{GonFund}. 
In particular the coproduct in the Hopf algebra is a new powerful structure, 
 which is not visible  
 on the level of numbers. It is conjectured that any relation of algebraic-geometric origin 
between the iterated integrals gives rise to a similar relation 
between the corresponding Hodge/motivic  iterated integrals. 
Therefore upgrading iterated integrals to elements of a Hopf algebra 
we see new structures, but conjecturally do not lose any information. 

\subsection{Multiple polylogarithms and the Hopf algebra of algebraic cycles} 
Given any field $F$, Bloch and Kriz \cite{BK} constructed,  
using algebraic cycles, yet another 
graded Hopf algebra $\chi_{\rm cycle}(F)$, denoted $\chi_\Mot(F)$ in loc.cit., and conjectured 
that it is isomorphic to the 
fundamental motivic Tate Hopf algebra of $F$.

 Thus it is natural to ask whether we can find elements 
\begin{equation} \label{1m}
{\Li}^{\rm cycle}_{n_1,\dots,n_m}(z_1,\dots,z_m) \in \chi_{\rm cycle}(F), \quad z_i \in \Fx\,,
\end{equation}
corresponding to multiple polylogarithms 
in the {\it  Hopf algebra of algebraic cycles} $\chi_{\rm cycle}(F)$ (defined in (\ref{BKchi}) below). 
The elements corresponding to the classical 
polylogarithms have been 
defined in \cite{BK} by using  the Bloch-Totaro cycles \cite{Bl1}. 
The main result of our paper is a 
construction of elements (\ref{1m}) for generic $z_i \in \Fx$, 
where generic means that the 
products of the arguments $\prod_{j=i}^k z_j$ for any pair $(i,k)$, where 
$1\leq i\leq k\leq n$,  are all different from~1.
We show that the $\Q$-linear combinations of these  elements 
form a Hopf subalgebra
\begin{equation}
\widetilde \chi^{\rm MP}_{\rm cycle}(F) \subset \chi_{\rm cycle}(F)\,.
\end{equation}
Presumably one can define the elements 
(\ref{1m}) for all 
$z_i \in \Fx$. They should generate the ``true''  
multiple polylogarithm Hopf subalgebra 
$
 \chi^{\rm MP}_{\rm cycle}(F) \subset \chi_{\rm cycle}(F)
$. 

It follows from  Conjecture 17a) in \cite{GonICM} that, for any field $F$, the 
$\Q$-linear combinations of the elements (\ref{1m}) should 
span over $\Q$
 the fundamental motivic Hopf algebra. In other words, 
one should have  
$$
\chi^{\rm MP}_{\rm cycle}(F) \stackrel{?}{=} \chi_{\rm cycle}(F)\,.
$$

\subsection{Strategy} Let ${\cal N}^{\bullet}(F, p)$ be Bloch's 
weight $p$ cubical cycle complex \cite{BlochAlgCyc}. 
The direct sum of these complexes 
$$
{\cal N}^{\bullet}(F, \ast):= \bigoplus_{p \geq 0} {\cal N}^{\bullet}(F, p)
$$ 
has a structure of graded-commutative augmented differential graded  algebra 
(DGA). Applying the bar construction functor $\B$ to it and then taking  
the ${\rm zero^{th}}$ cohomology of the resulting complex we get the Bloch-Kriz  
Hopf algebra: 
\begin{equation}\label{BKchi}
\chi_{\rm cycle}(F):= \H^0\B({\cal N}^{\bullet}(F, \ast))\,.
\end{equation}

Our Hopf subalgebra  $\widetilde\chi^{\rm MP}_{\rm cycle}(F)$ is 
also defined using the functor  $\H^0\B$. 
Its definition is organized into three steps outlined below. 

\begin{enumerate} \item[1.] {\bf Algebraic cycles from trees}.
Let $R$ be a set.  An {\em $R$-deco tree} is a plane tree 
with additional data:  
a root vertex of valency~1 and a decoration of the external vertices by elements $r_i$ of  
$R$, as in the picture (where the encircled vertex, decorated by $r_5$,
is the root vertex).

 \vskip 10pt \hskip 100pt
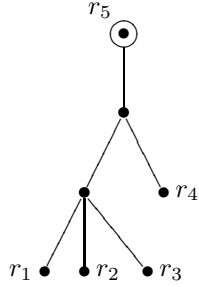
\hbox{
 \vbox{\xy 0;<-3pt,0pt>:
 \POS(0,0) *+{\bullet} *\cir{} 
 \ar @{-} +(0,10)*{\bullet}
  \POS(0,10) *{\bullet}
 \ar @{-} +(5,10)*{\bullet}
 \ar @{-} +(-5,10)*{\bullet}
 \POS(5,20) *{\bullet}
 \ar @{-} +(5,10)*{\bullet}
 \ar @{-} +(0,10)*{\bullet}
 \ar @{-} +(-8,10)*{\bullet}
 \POS(3,-3) *{r_5}
\POS(13,30) *{r_1}
\POS(2,30) *{r_2}
\POS(-6,30) *{r_3}
\POS(-8,20) *{r_4} \endxy}
}
\vskip 10pt

We define the grading of a tree as the number of its edges.   
Consider the graded $\Q$-vector space generated by 
$R$-deco trees. 
Let ${\CT}_{\bullet}(R)$ be the graded-commutative algebra 
generated by this vector space.  
A set of  linear generators  for ${\CT}_{\bullet}(R)$ 
is given by monomials in $R$-deco
trees, called {\it $R$-deco forests}. It has another grading, given by the number of leaves. 
So we get a bigraded algebra ${\CT}^{\ast}_{\bullet}(R)$. 
There is a differential 
${\CT}^{\ast}_{\bullet}(R) \to {\CT}^{\ast}_{\bullet-1}(R)$  
defined by contracting individual edges of a forest.
It equips ${\CT}^{\ast}_{\bullet}(R)$ with a structure of graded-commutative  
DGA with an extra grading. It is a rooted version of the DGA defined by the second author in \cite{GonArb}.

Now take $R = \Fx$, where $F$ is a field. 
We consider the subclass of ``generic'' decorations (as in Definition \ref{genericdeco})
on plane rooted trees. The corresponding 
DGA is denoted by $\widetilde {\CT}^{\ast}_{\bullet}(\Fx)$. 
We associate to such an $\Fx$-decorated tree with $m$ edges 
an element of the cubical cycle group ${\cal N}^{2p-m}(F, p)$ given by the 
codimension $p$ admissible algebraic 
cycles in the $m$-dimensional algebraic cube. We show 
that by extending this map to forests by multiplicativity 
we get  a homomorphism of DGAs 
$$
\Phi: \widetilde {\CT}^{\ast}_{\bullet}(\Fx) \longrightarrow 
{\cal N}^{2\ast - \bullet}(F, \ast)
$$
 which we call the {\sl forest cycling map}.
  
\ms
\item[2.] {\bf Trees from polygons.}
In a second step, 
we single out distinguished combinations of the above ($R$-deco) trees  
encoded by ($R$-deco) polygons. Namely, given a  polygon with
additional data (a root side and a decoration),  
we consider all its triangulations.
Each triangulation defines by duality a  plane trivalent 
tree. Let  $\Psi$ be the map 
which assigns  to the polygon the formal sum of these 
trees.  This map is 
extended to the graded-commutative algebra generated by the
$R$-deco polygons. 
The latter is equipped with a differential, which gives rise to a
graded-commutative DGA ${\CP}^{\ast}_{\bullet}(R)$. 
The DGA ${\CP}^{\ast}_{\bullet}(R)$ is identified with the standard cochain complex of 
 one of the  Lie coalgebras of iterated integrals defined in 
 \cite{GonFund}.
 
Let us consider its ``generic part'',  
a sub-DGA $\widetilde{\CP}^{\ast}_{\bullet}(R)$. 
The map  $\Psi$ provides a homomorphism of DGAs 
$$
\Psi: \widetilde{\CP}^{\ast}_{\bullet}(R) \longrightarrow \widetilde 
{\CT}^{\ast}_{\bullet}(R)\,.
$$
The {\em crucial fact} is the following. 
The image of the composition
\begin{equation} \label{2hs}
\widetilde 
{\CP}_{\bullet}^{\ast}(\Fx) \stackrel{\Psi}{\longrightarrow} 
\widetilde{\CT}_{\bullet}^{\ast}(\Fx) 
\stackrel{\Phi}{\longrightarrow} {\cal N}^{2\ast-\bullet}(F, \ast)\
\end{equation}
lies in the subspace  of elements with 
``decomposable''  differentials. 

\ms
\item[3.] {\bf The bar construction.}
Applying the $\H^0\B$ functor to (\ref{2hs}) we get the Hopf algebras 
$$
\chi_{\widetilde \CP}(\Fx):= \H^0\B (\widetilde {\CP}_{\bullet}^{\ast}(\Fx)), \quad 
\chi_{\widetilde \CT}(\Fx):= \H^0\B (\widetilde {\CT}_{\bullet}^{\ast}(\Fx)),
$$
and homomorphisms between them
\begin{equation} \label{3hs}
\chi_{\widetilde \CP}(\Fx) \stackrel{\H^0{\B}(\Psi)}{\longrightarrow} \chi_{\widetilde \CT}(\Fx) 
\stackrel{\H^0{\B}(\Phi)}{\longrightarrow} \chi_{\rm cycle}(F)\,.
\end{equation}

We prove that any $\Fx$-decorated polygon with the generic decoration gives rise 
to an element of $\chi_{\rm cycle}(F)$, called the {\it multiple logarithm cycle}. 
Its components in the bar complex are parametrized by 
dissections of the original polygon. 
By a slight generalization of this construction, we obtain 
our {\it multiple polylogarithm cycle}.

\ms
\item[4.] {\bf Comparison theorems.} The very name of the cycles suggests that 
they correspond to the motivic multiple polylogarithms from \cite{GonFund}. We give two
pieces of evidence for this: 

(i) The coproduct for the multiple polylogarithm cycles agrees with 
the one for the motivic multiple polylogarithms  
computed in \cite{GonFund}. 

(ii) Using  
the Hodge realization of the cycle Hopf algebra (\cite{BK}), we 
recover from our algebraic cycles the corresponding analytic functions.  

\end{enumerate}

\bigskip
\subsection{The structure of the paper}
In Section~2 we recall a few facts about 
Bloch's DGA of cubical algebraic cycles $\CN^\bullet(F, p)$, for a field $F$.  

In Section~3 we define, for a set  $R$, another 
DGA $\CT_\bullet^\ast(R)$, built from $R$-decorated rooted forests. 
A very similar differential graded algebra, for non-rooted forests, was introduced in \cite{GonArb}. 

In  Section~4 we relate the forest DGA to the cycle DGA in the case when $R = \Fx$. More precisely, we define 
a subalgebra $\widetilde \CT_\bullet^\ast(\Fx)$ of $\CT_\bullet^\ast(\Fx)$ 
by imposing some explicit genericity condition 
on the decoration.
Then we construct a map of graded DGA's 
$$
\Phi: \widetilde \CT_\bullet^\ast(\Fx) \longrightarrow \CN^{2\ast-\bullet}(F, \ast)\,.
$$

In  Section~5 we introduce the second important ingredient of our construction:  
given an ordered collection of elements $x_1, \dots, x_m \in \Fx$, we define an element 
$$
\tau(x_1, \ldots, x_m) \in \CT_\bullet^\ast(\Fx)\,
$$
with decomposable differential.
Under the genericity conditions on the $x_i$'s, it belongs to 
$\widetilde \CT_\bullet^\ast(\Fx)$,
and the algebraic cycle $\Phi\tau(x_1, \ldots, x_m)$ corresponds to the 
multiple logarithm (cf. \ref{mulpol}). 
Then we introduce $R$-deco 
polygons as more convenient combinatorial objects, and define a DGA 
of  $R$-deco polygons, which maps to the tree DGA from Section~4. 
 
We show in Section 8.2 that the DGA  of  $R$-deco polygons is identified with the standard cochain complex of 
 a  Lie coalgebra of (formal) iterated integrals defined in 
 \cite{GonFund}. 
There is a yet another differential $\overline \partial$ on the graded algebra of $R$-deco polygons, 
which plays only a technical role in the paper.  
In Section 5.4 the corresponding DGA is related to the rooted version of 
one of the graded Lie coalgebras introduced in \cite{GonDihedral} under 
the name  of the  dihedral Lie coalgebras.

In Section~6, we perform the bar construction on polygons. It gives rise to
a Hopf algebra with an induced map to $\chi_\cycle$.

In Section~7, an explicit formula for the coproduct in 
this Hopf algebra of polygons is given. 

In Section~8, we prove a coproduct comparison theorem: the coproduct 
of our multiple polylogarithm cycles agrees with   
 the coproduct of the corresponding 
motivic iterated integrals obtained in \cite{GonFund}. This ensures 
that our multiple polylogarithm cycles agree 
with the corresponding motivic iterated integrals.  We relate the formula for the coproduct in 
the Hopf algebra of polygons to 
the one from \cite{CK1}. 

In Section~9, we show how to get the original multivalued 
analytic functions (\ref{ml}) from the constructed cycles. 


In Section~10, we indicate how to modify our construction in order to obtain cycles corresponding to
multiple {\em poly}logarithms.

\vskip 2mm
See Sections 7-9 of \cite{GonECM} for connections between the sums of plane decorated 
trivalent trees and multiple polylogarithms closely related to the one discussed in our paper. 
\vskip 3mm
Our paper \cite{GGL} can serve as an introduction to the present paper.

\vskip 2mm
During the preparation of the final version of this paper, we learned that 
 H. Furusho and A. Jafari obtained 
 similar results, see \cite{FJ}.  

\medskip
{\bf Acknowledgements}. 
This work benefited tremendously from the hospitality of the MPI (Bonn). We 
are grateful to the MPI for providing ideal working conditions and  
support. H.G. is indebted to Spencer~Bloch and Dirk~Kreimer for enlightening discussions while being invited at the University 
of Chicago in April 2003 and at the IHES (Jan/Feb 2004 and April 2005).
He further acknowledges support for conferences at the ESI Vienna (Sept 2002), Les Houches (March 2003), Oxford Ohio (April 2003),
Morelia (June 2003) and Luminy (Feb 2004) where a preliminary version of this work (cf.~\cite{GGL}) has been presented.
A.G. was supported by the NSF grants DMS-0099390 and DMS-0400449. A.L. was partially supported by the grant RFFI 04-01-00642.

\section{Background on algebraic cycles}\label{preliminaries}

We review some of the main properties of cubical algebraic cycles. Our main references for the general set-up are 
the papers by Bloch \cite{Bl1}, \cite{Bl2}, \S2.4.1, and by Bloch and Kriz \cite{BK}, \S5. The reader can find more examples in \cite{GMS}.

Let $F$ be a field. Following  \cite{Bl2}, we define the algebraic $1$-cube $\cub_F$ as a pair
$$\cub_F = \Big(\P_F^1 \setminus \{1\} \simeq \A_F^1\,, (0) - (\infty)\Big).$$
Here we consider the 
 standard coordinate $z$ on the projective line $\P_F^1$ and remove from it the point 1. Furthermore, 
 $(0) - (\infty)$ denotes 
the divisor defined by the two points $0$ and $\infty$. The algebraic $n$-cube is defined by setting 
$\cub_F^n=(\cub_F)^n$.  

Bloch defined (we use the conventions from \cite{BK}, \S5), for $p,n\in \N$, the cycle groups
\begin{multline*}
\Cycle^p\langle F,n\rangle = \Z\big[\text{\{admissible closed irreducible subvarieties over $F$,}\\
 \text{\qquad \qquad of codimension $p$ in $\cub_F^n$\}}\big]\,.\qquad\qquad\qquad
\end{multline*}
Here a cycle is called {\bf admissible} if it intersects all the faces (of any 
codimension) of $\cub_F^n$ properly, i.e., in codimension~$p$ or not at all.

Let $G(k)$ be the semidirect product 
of the symmetric group $S_k$ and the group $(\Z/2\Z)^k$, acting by 
  permuting and inverting the 
coordinates in $\cub_F^k$. Let $\varepsilon_k$ 
be the sign representation of this group. 
We set in a slight variation to \cite{BK}, indicated by a subscript $\rmc$ for ``coinvariants'',
$$
\CN_\rmc^n(F,p) =  \big(\Cycle^p\langle F,2p-n\rangle\otimes 
\varepsilon_{2p-n}\big)_{G(2p-n)} \,.$$

\medskip\noindent
{\bf Remark:} The cycle groups were defined in \cite{Bl1} using invariants instead of coinvariants,
and therefore the groups above agree with those rationally.
We will refer to them as the {\em cubical algebraic cycle groups}.

\smallskip
Putting the group $\CN_\rmc^n(F,p)$ in 
degree $n$, we get, for a fixed ~$p$, a complex
$$ \dots \to \CN_\rmc^n(F,p) \ { \buildrel \partial \over \to}\ \CN_\rmc^{n+1}(F,p) \to \dots $$
where the differential $\partial$ is given by 
$$\partial = \sum_{i=1}^n (-1)^{i-1} (\partial_0^i - \partial_\infty^i)$$
and $\partial_\varepsilon^i$ denotes 
the operator given by the intersection with the coordinate hyperplane 
$\{z_i=\varepsilon\}$, $\varepsilon \in \{0,\infty\}$. 

The concatenation of coordinates, followed by 
taking the corresponding co\-inva\-riants, provides a product on algebraic cycles,
and together with the above one gets (\cite{Bl1}, Lemma 1.1, and \cite{BK}, Prop.~5.1):

\begin{Proposition} Given a field $F$, the algebraic cycle groups $\CN_\rmc^n(F,p)$ give rise to
a differential graded algebra $\CN_\rmc^\bullet = \bigoplus_{p\geq 0} \CN_\rmc^\bullet(F,p)$. 
Its cohomology groups are the higher Chow groups of $F$.
\end{Proposition}
The higher Chow groups of a field $F$ are known, by Bloch (\cite{BlochAlgCyc} and \cite{BlochMovLem}) and 
also by Levine \cite{LevHCG}, to be rationally isomorphic 
to its algebraic $K$-groups.

\medskip\noindent
\begin{Example} 
\begin{enumerate}
\item[a)] In $\cub_F^3$, a curve 
intersecting the line $z_1=z_2=0$ is not admissible, 
since an admissible curve has to meet
faces of dimension 1 in codimension 2, i.e. not at all. 

\item[b)] For any element $a$ in $F$, one can associate a cubical algebraic cycle which corresponds to
the dilogarithm $Li_2(a)$ if $F\hookrightarrow \C$. This cycle has been given by Totaro as 
the image of the map 
\begin{eqnarray*}
\varphi_a : \P_F^1 &\to& \big(\P_F^1\big)^3\,, \\
             t &\mapsto & (t,1-t, 1-\frac a t)\,,
\end{eqnarray*}
restricted to the algebraic cube $\cub_F^3$: we write
$$C_a := \Big[t,1-t,1-\frac a t\Big] := \varphi_a\big(\P_F^1 \big) \cap \cub_F^3\,.$$
The cycle $C_a$ belongs to the group $\CN_\rmc^1(F,2)$. One has 
$$ \partial C_a = [a,1-a]\in \cub_F^2$$
(only $\partial_0^3$ gives a non-empty contribution). 
The same computation shows that $C_a$ is in fact admissible. Observe the apparent similarity with 
the 
formula $d\, Li_2(a) = -\log(1-a)\,d\log(a)$ for the differential of  the dilogarithm. 

\item[c)]
Recall that the double logarithm can be defined by an iterated integral
$$ I_{1,1}(x_1,x_2)= \int_{\Delta_\gamma} \frac{dt_1}{t_1-x_1} \wedge \frac{dt_2}{t_2-x_2}\,,$$
in the notation of (\ref{int_def}) in the Introduction.
Its differential can be easily computed as
\begin{equation} \label{diffLi}
\qquad d\,I_{1,1}(x_1,x_2) = I_1(x_1)\;d\,I_1(x_2) - I_1(x_1)\;d\,I_1\Big(\frac{x_2}{x_1}\Big)  + I_1(x_2)\;d\,I_1
\Big(\frac{x_1}{x_2}\Big)\,.
\end{equation}

The cycle
\begin{equation} \label{Z2}
Z_{x_1,x_2}:=  \Big[1-\frac {1}t, 1-\frac t {x_1},1-\frac t{x_2}\Big] \,,
\end{equation}
will play the role of the double logarithm among the algebraic cycles. Its boundary is readily evaluated as
\begin{equation}\label{bdry}
\qquad \partial Z_{x_1,x_2} = \underbrace{\Big[1-\frac 1{x_1}, 1-\frac 1{x_2}\Big]}_{\text{from }z_1=0} - \underbrace{\Big[1-\frac 1{x_1}, 1-\frac {x_1}{x_2}\Big]}_{\text{from }z_2=0} + \underbrace{\Big[1-\frac1{x_2},1-\frac{x_2}{x_1}\Big]}_{\text{from }z_3=0} \ \in 
\CN_\rmc^2(F,2)
\end{equation}
whose individual terms are already very reminiscent of the three terms in (\ref{diffLi}). 
\end{enumerate}
\end{Example}
The cycle associated to the triple logarithm $I_{1,1,1}(x_1,x_2,x_3)$ is given in Example \ref{triplelog} below.

\section{The differential graded algebra of $R$-deco forests}
In this paper, a {\em plane tree} 
is a finite tree whose internal vertices are of valency~$\geq 3$, and where at each vertex a cyclic ordering 
of the incident edges is given. We assume that 
all the other vertices are of valency~1, and call them {\em external} vertices. A plane tree is {\em planted}  if it
has a distinguished external vertex of valency~1, called its {\em root}; in particular, 
a planted tree has at least one edge. 
A {\em forest} is a disjoint union of trees. 

Below we work with algebras over $\Q$, although we may replace it by any field.

\subsection{The orientation torsor}
Recall that a {\it torsor} under a group $G$ is a set on which $G$ acts freely transitively.

Let $S$ be a finite (non-empty) set. We impose on the set of orderings of $S$ an equivalence relation, given by
even permutations of the elements. The  equivalence classes form a 2-element set ${\rm Or}_S$ of {\em orientation classes},
for short also {\em orientations}. It has an obvious $\Z/2\Z$-torsor structure and is called the {\it orientation
torsor of $S$}. 

\begin{Definition} 
The {\bf orientation torsor} of  a plane forest is the orientation
torsor of the set of its edges. 
\end{Definition}

\subsection{The algebra of $R$-deco forests}\label{forests}
\begin{Definition}\label{Rdeco} Let $R$ be a (non-empty) set. An {\bf $R$-deco tree} is a planted plane tree with a map, 
called {\bf $R$-decoration}, from its external vertices to $R$. An {\bf $R$-deco forest} is a disjoint 
union of $R$-deco trees.  
\end{Definition}
\begin{Remark}\label{rem1}
\begin{enumerate}
\item There is a canonical  direction for each edge in an $R$-deco tree, away from the root.
\item An edge ordering in a tree $\tau$ provides an orientation on $\tau$.
For a planted plane tree there is a 
canonical linear order of the edges, starting from the root edge, which is 
induced by the cyclic order of edges at internal vertices.
\end{enumerate}
\end{Remark}
We draw the trees so that the cyclic order of edges around
the vertices is displayed in counterclockwise direction, and the root vertex is at the top.

\begin{Example} We draw an $R$-deco tree $\tau$ with root vertex decorated by $x_4\in R$; its other external vertices are decorated
by $x_1$, $x_2$, $x_3\in R$. The above-mentioned ordering of the edges $e_i$ coincides
with the natural ordering of their indices, while the direction
of the edges (away from the root) is indicated by small arrows along the edges. The root vertex is marked by a circle 
around it.
\end{Example}
\vskip 15pt
\hskip 100pt 
\xy  0;<-3pt,0pt>:
(5,0)*+{};(-5,-1)*+{}
 **\crv{(6,5)&(7,20)&(20,30)&(10,45)&(6,25)&(5,25)&(4,30)&(0,42)&(-6,35)&(-3,25)&(4,18)&(-1,15)&(-2,22)&(-7,28)&(-13,20)&(-8,15)&(-5,5)}
?(.03)*\dir{>} ?(.35)*\dir{>} ?(.93)*\dir{>}
 \POS(0,0) *+{\bullet} *\cir{} 
 \ar @{-} +(0,10)*{\bullet}^{e_1} \POS(0,6)*\dir3{>}
  \POS(0,10) *{\bullet}
 \ar @{-} +(5,10)*{\bullet}_{e_2} \POS(3,16)*\dir3{>}\POS(0,10)
 \ar @{-} +(-5,10)*{\bullet}^{e_5} \POS(-3,16)*\dir3{>}
 \POS(5,20) *{\bullet}
 \ar @{-} +(5,10)*{\bullet}_{e_3} \POS(8,26)*\dir3{>}\POS(5,20)
 \ar @{-} +(-5,10)*{\bullet}^{e_4} \POS(2,26)*\dir3{>}
\POS(3,-3) *{x_4}
\POS(13,32) *{x_1}
\POS(-2,32) *{x_2}
\POS(-8,20) *{x_3}
\POS(-13,10) *{\tau}
\endxy

\medskip\noindent
\centerline{\small Figure 1: An $R$-deco tree $\tau$ with root vertex decorated by $x_4$.}
\vskip 15pt

We define the {\bf grading} of an $R$-deco tree $\tau$  by
$$e(\tau) = \#\{\text{edges of }\tau\}\,$$ 
and extend it to forests by linearity: $e(\phi_1 \disjoint \phi_2): = e(\phi_1) + e(\phi_2)$. (Here 
$\disjoint$ denotes the disjoint union.)

Let $V^\tr_{\bullet}({R})$ be the graded vector space where $V^\tr_n(R)$ for $n>0$ has a basis given by 
$R$-deco trees of grading $n$ and where $V_0^\tr(R):=\Q\cdot\fatone$
for an extra generator $\fatone$.

\begin{Definition} 
The  algebra $\CT_{\bullet}(R)$ is the free graded-commutative algebra with unit $\fatone$,
generated by the graded vector space $V^\tr_{\bullet}({R})$.
\end{Definition}

Here is a definition which does not use a choice of orientation.

\begin{Definition} 
For a set $R$,  ${\CT_{\bullet}^\Or}(R)$ is generated as a vector space by pairs $(\phi, \omega)$, where $\phi$ 
is an $R$-deco forest and $\omega$ an orientation on it, subject to the relation $(\phi, -\omega) = -(\phi, \omega)$. 
\end{Definition}

We define an algebra structure $\star$ on ${\CT}_{\bullet}^\Or(R)$ by letting $\fatone$ be the neutral element with respect 
to $\star$ and by setting ($\disjoint$ denotes the disjoint union)
$$
(\phi_1, \omega_1)\star (\phi_2, \omega_2) := (\phi_1\disjoint \phi_2, \omega_1 \otimes \omega_2)\,.
$$
It makes ${\CT}_{\bullet}^\Or(R)$ into a graded-commutative algebra,  
so the basis elements of $V^\tr_{\bullet}({R})$ commute in the $R$-deco tree algebra via the rule
$$ (\tau_1,\omega_{1})\star (\tau_2,\omega_{2}) = (-1)^{e(\tau_1)e(\tau_2)} (\tau_2,\omega_{2})\star (\tau_1,\omega_{1})\,.$$

Since each $R$-deco tree comes equipped with a canonical orientation (cf.~Remark \ref{rem1}), we can---and in the 
following will---identify $\CT_{\bullet}(R)$ with $\CT_{\bullet}^\Or(R)$.
This identification will be useful in describing the differential, since the orientation
torsor takes care of the signs which are more complicated to describe for $\CT_{\bullet}(R)$.
We will refer to either algebra as the {\bf $R$-deco tree algebra}.

\subsection{The differential}
A differential on ${\CT}_{\bullet}(R)\cong\CT_{\bullet}^\Or(R) $ is a map
$$
d: {\CT}_{\bullet}(R) \longrightarrow {\CT}_{\bullet-1}(R) 
$$
satisfying $d^2=0$ and the Leibniz rule 
$$d\big((\tau_1,\omega_{1})\star (\tau_2,\omega_{2})\big) = d\big((\tau_1,\omega_{1})\big)\star (\tau_2,\omega_{2}) +
(-1)^{e(\tau_1)}(\tau_1,\omega_{1})\star d\big((\tau_2,\omega_{2})\big)\,.$$
Since ${\CT}_{\bullet}(R)$ is a free graded-commutative algebra, 
it is sufficient to define the differential on the algebra generators, that is on the elements $(\tau, \omega)$, where $\tau$ is an $R$-deco tree
and $\omega$ is an orientation of $\tau$.

The terms in the differential of a tree $\tau$ arise by contracting an edge of~$\tau$---they fall into two types,
according to whether the edge is internal or external. We will need the notion of a splitting.

\begin{Definition} A {\bf splitting} of a tree $\tau$ at an internal vertex $v$ is the disjoint union of 
the trees which arise as $\tau_i\cup v$ where the $\tau_i$ are the connected components of $\tau\setminus v$.

The following further structures on $\tau$  are inherited for each $\tau_i\cup v$: a decoration at $v$, planarity of $\tau$ 
and an ordering of its edges.
Also, if $\tau$ has a root $r$, then $v$ plays the role of the root for all $\tau_i\cup v$ which do not contain $r$.
\end{Definition} 

\begin{Definition} Let $e$ be an edge of a tree $\tau$. The {\bf contraction} of $\tau$
along $e$, denoted $\tau/e$, is given as follows:
\begin{enumerate}
\item If $e$ is an internal edge, then $\tau/e$ is again a tree: it is the same tree as $\tau$ 
except that $e$ is contracted and the incident vertices $v$ and $v'$ of $e$ are identified 
to a single vertex.
\item If $e$ is an external edge not containing the root vertex, then $\tau/e$ is obtained as follows: first we contract
the edge $e$ to a vertex $w$ and then we perform a splitting at $w$.
\item If the tree consists of a single edge, its contraction gives the empty tree.
\end{enumerate}
\end{Definition} 

Two typical examples are given below: in 
Figure~2 we contract a {\em leaf}, i.e.~an external vertex which is not the root vertex,
and in Figure~3 the root edge is contracted.
\vskip 10pt

\vbox{\hskip 0pt
 \hbox{
\vtop{\xy  0;<-3pt,0pt>: 
 \POS(0,0) *+{\bullet} *\cir{} 
 \ar @{-} +(0,10)*{\bullet}
  \POS(0,10) *{\bullet}
 \ar @{-} +(8,10)*{\bullet}
 \ar @{-} +(0,10)*{\bullet}
 \ar @{-} +(-8,10)*{\bullet}
\POS(3,-3) *{p}
\POS(-11,20) *{s} 
\POS(3,20) *{r}
\POS(-2,17) *{e}
\POS(12,20) *{q}
\endxy}
 \vtop{\hbox{$\quad$}\hbox{\small contract}\xymatrix{\ar@{~>}[r]&}\hbox{\small along $e$}}
\vbox{\xy  0;<-3pt,0pt>: 
 \POS(0,0) *+{\bullet} *\cir{} 
 \ar @{-} +(0,10)*{\bullet}
  \POS(0,10) *{\bullet}
 \ar @{-} +(8,10)*{\bullet}
 \ar @{-} +(-8,10)*{\bullet}
\POS(3,-3) *{p}
\POS(-11,20) *{s} 
\POS(3,9) *{r}
\POS(12,20) *{q}
\endxy}
 \vtop{\hbox{$\quad$}\hbox{\small split at}\xymatrix{\ar@{~>}[r]&}\hbox{\small internal vertex}\hbox{$\quad$}}
\vbox{\xy  0;<-3pt,0pt>: 
 \POS(0,0) *+{\bullet} *\cir{} 
 \ar @{-} +(0,10)*{\bullet}
  \POS(3,13) *+{\bullet} *\cir{}
 \ar @{-} +(8,10)*{\bullet}
  \POS(-3,13) *+{\bullet} *\cir{}
 \ar @{-} +(-8,10)*{\bullet}
\POS(3,-3) *{p}
\POS(-11,20) *{s} 
\POS(0,12) *{r}
\POS(12,20) *{q}
\endxy}
 }
}
\vskip 10pt
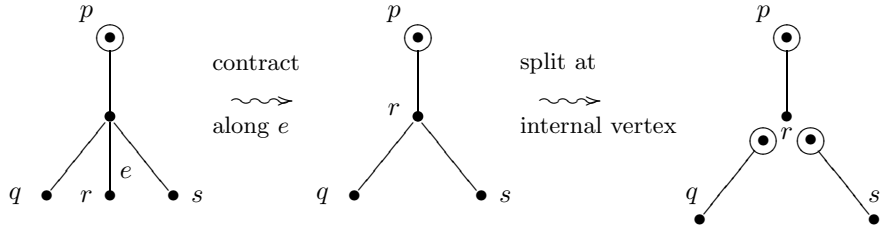
\centerline{\small Figure 2: Contracting a leaf.}

\vskip 15pt

\vbox{\hskip 0pt
 \hbox{
\vbox{\xy  0;<-3pt,0pt>: 
 \POS(0,0) *+{\bullet} *\cir{} 
 \ar @{-} +(0,10)*{\bullet}
  \POS(0,10) *{\bullet}
 \ar @{-} +(8,10)*{\bullet}
 \ar @{-} +(0,10)*{\bullet}
 \ar @{-} +(-8,10)*{\bullet}
\POS(3,-3) *{p}
\POS(-11,20) *{s} 
\POS(3,20) *{r}
\POS(-2,5) *{e'}
\POS(12,20) *{q}
\endxy}
 \vtop{\hbox{$\quad$}\hbox{\small contract}\xymatrix{\ar@{~>}[r]&}\hbox{\small along $e'$}\hbox{$\quad$}}
\vbox{\xy  0;<-3pt,0pt>: 
  \POS(0,10) *+{\bullet} *\cir{} 
 \ar @{-} +(8,10)*{\bullet}
 \ar @{-} +(0,10)*{\bullet}
 \ar @{-} +(-8,10)*{\bullet}
\POS(-11,20) *{s} 
\POS(3,20) *{r}
\POS(0,5) *{p}
\POS(12,20) *{q}
\endxy}

 \vtop{\hbox{$\quad$}\hbox{\small split at}\xymatrix{\ar@{~>}[r]&}\hbox{\small root vertex}\hbox{$\quad$}}
\vbox{\xy  0;<-3pt,0pt>: 
  \POS(9,10) *+{\bullet} *\cir{}
 \ar @{-} +(8,10)*{\bullet}
  \POS(-3,10) *+{\bullet} *\cir{}
 \ar @{-} +(-8,10)*{\bullet}
  \POS(3,10) *+{\bullet} *\cir{}
 \ar @{-} +(0,10)*{\bullet}
\POS(3,5) *{p}
\POS(-3,5) *{p}
\POS(9,5) *{p}
\POS(-8,20) *{s} 
\POS(0,20) *{r}
\POS(12,20) *{q}
\endxy}
 }
}
\vskip 10pt
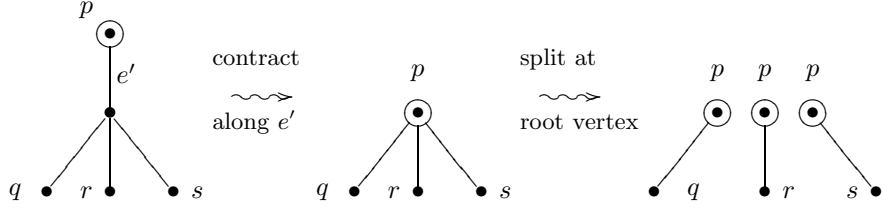
\centerline{\small Figure 3: Contracting the root edge.}
\vskip 15pt

\begin{Remark} The linear ordering on the edges of $\tau$ induces a linear ordering on each of
the trees which result from a splitting, and it furthermore induces a linear ordering on the
trees in the resulting forest. The notion of an {\em unshuffle} is appropriate to describe
this.
\end{Remark}

\medskip
\def \sgn {{\rm sgn}}
Let $S$ be a finite set and ${\rm Or}_S$ the orientation torsor of $S$. 
We present its elements as $s_1 \wedge ... \wedge s_n$, where $n = |S|$, with the relation 
$s_{\sigma(1)}\wedge\dots\wedge s_{\sigma(n)} = \sgn(\sigma) \big(s_1\wedge\dots\wedge s_n\big)$ for
any permutation $\sigma$ on $n$ letters. For $|S|=1$ there is a unique orientation.

Given $s \in S$ and 
$\omega\in {\rm Or}_S$, we define 
an element $i_s\omega \in {\rm Or}_{S-s}$ as follows:
$$
i_s\omega:= s_2 \wedge ... \wedge s_n \quad \mbox{if $\omega =  s \wedge s_2 \wedge ... \wedge s_n$}\,.
$$

\begin{Definition} \label{defdiff} Let $\tau$ be a finite tree with set of edges $E(\tau)$, and let $\om$ be an orientation of $\tau$. 
The {\bf differential} on $(\tau,\om)$ is defined as
$$d: (\tau,\om) \mapsto \sum_{e\in E(\tau)}(\tau/e, i_e\om)\,.$$
In particular, $d$ maps a tree with at most one edge to zero (which corresponds to the empty tree).
\end{Definition}

\begin{Proposition} \label{showdiff} The map $d$ just defined is in fact a differential, i.e. $d^2=0$.
\end{Proposition}

\medskip \noindent
{\bf Proof.} For a tree $\tau$ with one edge only we already have $d(\tau)=0$.
Applying $d$ to a generator $(\tau,\omega)$ with $\geq 2$ edges, we obtain $\quad \sum_{e\in E(\tau)} (\tau/e,i_e\omega)$,
and a second application of $d$ gives (with the obvious notation $\tau/(e\disjoint e')$ for the 
independent contraction of the edges $e$ and $e'$ to a point each)
\begin{eqnarray*}
 \sum_{e'\in E(\tau/e)}
\big((\tau/e)/e', i_{e'}(i_e(\omega))\big)  &=& \sum_{{\buildrel \scriptstyle{(e,e')\in E(\tau)^2}\over { e\not= e'}}} \big(\tau\big/(e\disjoint e'), i_{e'}
(i_e \omega)\big) \\
&=& \sum_{e<e'} \big(\tau\big/(e\disjoint e'), i_{e'}
(i_e \omega)\big) + \big(\tau\big/(e\disjoint e'), i_e
(i_{e'} \omega)\big)
\end{eqnarray*}
which vanishes since $i_{e'} i_e \omega = - i_e i_{e'}\omega$.
Note that this argument works also if $e$ and $e'$ are contracted to the same point.
 \qed

\begin{Proposition}
The differential of Definition \ref{defdiff} is compatible with the $R$-deco structure.
\end{Proposition}
{\bf Proof.} We can identify $R$-deco forests which arise from successively contracting edges in an $R$-deco
tree with trees having external and possibly internal decoration. Thus, as long as the two contracted
edges are not adjacent leaves, the same proof as for the previous proposition applies.

If the two contracted edges are adjacent leaves, then the contraction of the second one gives zero,
since it corresponds to contracting a tree with a single edge. \qed

\begin{Example}
The simplest non-trivial example for the differential of an $\Fx$-deco tree, where $\Fx$ is the multiplicative group of a field $F$, can be seen on a tree with one internal vertex, 
as given in Figure~4. Here we choose the $\Fx$-decoration $(x_1,x_2)$ with $x_1,x_2\in \Fx$ for the leaves and the
decoration $1$ for the root:
\vskip 3pt
\hskip 10pt
\hbox{
 \vbox{\xy  0;<-3pt,0pt>: 
 \POS(0,0) *+{\bullet} *\cir{} 
 \ar @{-} +(0,10)*{\bullet}
  \POS(0,10) *{\bullet}
 \ar @{-} +(8,10)*{\bullet}
 \ar @{-} +(-8,10)*{\bullet}
 \POS(8,20) *{\bullet}
 \POS(-8,20) *{\bullet}
\POS(3,-3) *{1}
\POS(8,23) *{x_1}
\POS(-8,23) *{x_2}
\endxy}
\vtop{\hbox{$\quad$}\hbox{$\quad$}\hbox{$\quad$}\hbox{$\displaystyle{\buildrel d\over \longmapsto}$\quad }
\hbox{$\quad$}\hbox{$\quad$}}
 \vbox{\xy  0;<-3pt,0pt>: 
 \POS(0,5) *+{\bullet} *\cir{} 
 \ar @{-} +(0,10)*{\bullet}
\POS(3,2) *{1}
\POS(0,18) *{x_1}
\endxy}
\vtop{\hbox{$\quad$}\hbox{$\quad$}\hbox{$\quad$}\hbox{\ $\displaystyle{\star}$ }\hbox{$\quad$}\hbox{$\quad$}}
 \vbox{\xy  0;<-3pt,0pt>: 
 \POS(0,5) *+{\bullet} *\cir{} 
 \ar @{-} +(0,10)*{\bullet}
  \POS(0,15) *{\bullet}
\POS(3,2) *{1}
\POS(0,18) *{x_2}
\endxy}
\vtop{\hbox{$\quad$}\hbox{$\quad$}\hbox{$\quad$}\hbox{$\displaystyle{\quad - \ }$ }\hbox{$\quad$}\hbox{$\quad$}}
 \vbox{\xy  0;<-3pt,0pt>: 
 \POS(0,5) *+{\bullet} *\cir{} 
 \ar @{-} +(0,10)*{\bullet}
  \POS(0,15) *{\bullet}
\POS(3,2) *{1}
\POS(0,18) *{x_1}
\endxy}
\vtop{\hbox{$\quad$}\hbox{$\quad$}\hbox{$\quad$}\hbox{\ $\displaystyle{\star}$ }\hbox{$\quad$}\hbox{$\quad$}}
 \vbox{\xy  0;<-3pt,0pt>: 
 \POS(0,5) *+{\bullet} *\cir{} 
 \ar @{-} +(0,10)*{\bullet}
  \POS(0,15) *{\bullet}
\POS(0,1) *{x_1}
\POS(0,18) *{x_2}
\endxy}
\vtop{\hbox{$\quad$}\hbox{$\quad$}\hbox{$\quad$}\hbox{$\displaystyle{\quad +\ }$ }\hbox{$\quad$}\hbox{$\quad$}}
 \vbox{\xy  0;<-3pt,0pt>: 
 \POS(0,5) *+{\bullet} *\cir{} 
 \ar @{-} +(0,10)*{\bullet}
  \POS(0,15) *{\bullet}
\POS(3,2) *{1}
\POS(0,18) *{x_2}
\endxy}
\vtop{\hbox{$\quad$}\hbox{$\quad$}\hbox{$\quad$}\hbox{\ $\displaystyle{\star}$ }\hbox{$\quad$}\hbox{$\quad$}}
 \vbox{\xy  0;<-3pt,0pt>: 
 \POS(0,5) *+{\bullet} *\cir{} 
 \ar @{-} +(0,10)*{\bullet}
  \POS(0,15) *{\bullet}
\POS(0,1) *{x_2}
\POS(0,18) *{x_1}
\endxy}
}

\vskip 7pt

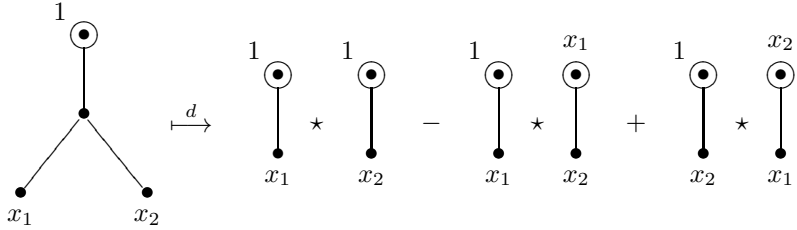
\centerline{\small Figure 4: The differential on a tree with one internal vertex.}

\vskip 7pt
\noindent For the drawings, we use the canonical ordering of edges (cf.~Remark \ref{rem1}) for $\Fx$-deco 
trees and the induced
ordering of edges for forests which arise from a splitting. The edge ordering of the forest 
produces its orientation, i.e. the choice of an element in the orientation torsor. 
\end{Example}

\begin{Remark}
There is a $\Z$-bigrading on $\CT_{\bullet}(R)$ given by 
$$
\CT_n^p(R) = \Z[\{\text{$R$-deco forests with $n$ edges and $p$ leaves}\}]\,.
$$
For $n\leq 0$ or $p\leq 0$ it is zero, except for $p=n=0$ where it is $\simeq \Q$. 
It will correspond to the bigrading on the cycle groups $\CN_\rmc^{2p-n}(F,p)$. 
\end{Remark}

For a set $R$, put
$$
\CT_{\bullet}^\ast(R):= \bigoplus_{p\geq 0}\ \bigoplus_{0\leq n\leq p}\CT_n^p(R)\,.
$$
We have an obvious augmentation map $\varepsilon:\CT_\bullet^\ast(R)\to \Q$ defined by taking the zero${}^\th$ component.

With the above definitions, we have:
\begin{Proposition} The algebra $\CT_{\bullet}^\ast(R)$, together with the differential $d$ above,
forms a bigraded DGA with augmentation. The differential lowers
the grading indicated with subscripts by~1 and leaves the grading indicated by superscripts unchanged.
\end{Proposition}

\medskip\noindent
{\bf Proof.} We have already proved in Proposition \ref{showdiff} that $d^2=0$.
It remains to show that the differential is compatible with the bigrading. In fact,
we have more precisely
$$d:\CT_n^p(R) \to \CT_{n-1}^p(R)\,,$$
where $\CT_0^p(R) = 0$ for any $p>0$.

We distinguish three cases, dealing first with a tree $\tau$ in $\CT_n^p(R)$:
\begin{enumerate} 
\item If we contract an internal edge in $\tau$, we obtain a tree with one edge less (so the subscript grading
is lowered by~1) and with the same root and decorations of external edges.
\item For contractions of an external non-root edge, the vertex to which the edge is contracted
serves as a root vertex for all but one tree in the resulting forest. Therefore the number
of external non-root vertices does not change under the contraction while the number of edges obviously does. 
\item If the root edge is contracted, none of the external non-root edges changes its properties,
except when the tree consists of a single edge in which case the contraction gives the empty tree. In both cases
the claim follows immediately.
\end{enumerate}
Using the Leibniz rule, we get the corresponding statement for forests, and by linearity it follows for any
element in $\CT_n^p(R)$. \qed


\section{Mapping forests to algebraic cycles}

In the special case where $R=\Fx$, the multiplicative group of a field $F$, we can 
establish the connection between the two differential graded algebras
above, and Theorem \ref{dga_morphism} below gives the first main result of this paper.

It turns out that the admissibility condition on algebraic cycles mentioned above forces us to restrict
to a subalgebra of $\CT^\ast_\bullet$, which we now describe.

\begin{Definition} \label{genericdeco}
We call an $R$-deco tree  {\bf generic} if all the individual 
decorations of external vertices are different.

We denote the subalgebra of $\CT^\ast_\bullet(R)$ generated by generic $R$-deco trees by $\CTT^\ast_\bullet(R)$.
\end{Definition}

One of our key results is the following statement:
\begin{Theorem}\label{dga_morphism} For a field $F$, there is a natural map of differential graded algebras
$$\CTT^*_\bullet(\Fx) \to \CN_\rmc^{2*-\bullet}(F,*)\,.$$
It is given by the map in Definition \ref{fc} below.
\end{Theorem}

A proof of this statement is given at the end of this section.

\begin{Definition} \label{fc}
The {\bf forest cycling map} for a field $F$ is the map $\Phi$ from $\CT^\ast_\bullet(\Fx)$ to (not necessarily admissible) 
cubical algebraic cycles over $F$ given on generators, i.e. $\Fx$-deco trees $\tau$ with orientation $\omega$, as follows:
\begin{enumerate}
\item to each internal vertex $v$ of $\tau$ we associate a decoration consisting of an independent (``parametrizing'') variable;
\item to each oriented edge from a vertex $v$ to another one $w$ equipped with respective decorations $y_v$ and $y_w$
(variables or constants), we associate the expression $[1-y_v/y_w]$ as a parametrized coordinate in $\P_F^1$;
\item choosing an ordering of edges of $\tau$ corresponding to $\om$, we concatenate all the respective coordinates 
produced in the previous step.
\end{enumerate}
\end{Definition}
The  last step in the definition is well-defined due to fact that we take the coinvariants of algebraic cycles with respect
to permutation of coordinates: changing the orientation of $\tau$ amounts to performing an odd permutation of the
coordinates of $\Phi(\tau)$, so in both settings we multiply by $-1$.

The above somewhat lengthy description of the map $\Phi$ is easily understood by looking at an example. 
We denote the concatenation product for algebraic cycles by~$*$, and
we encode the expression $1-\frac x y$, for $x,y$ in a field, by the following picture
\vskip 5pt 
\vbox{\hskip 150pt \xy  0;<-3pt,0pt>: \POS(0,0) *{\bullet} 
 \ar @{-} +(0,10)*{\bullet} \POS(0,6)*\dir3{>}
  \POS(-3,10) *{y}
  \POS(-3,0) *{x}
\endxy }\vskip 5pt \noindent

\begin{Example} Let us consider the forest cycling map $\Phi$ for the following $R$-deco tree $(\tau,\om)$, where the 
orientation $\om$ is given by $e_1\wedge e_2\wedge e_3$ (we leave out the arrows 
since the edges are understood to be directed away from the root):
\vskip 10pt
\hskip -13pt
\hbox{
 \vtop{\xy  0;<-3pt,0pt>: 
 \POS(0,0) *+{\bullet} *\cir{} 
 \ar @{-} +(0,10)*{\bullet}
  \POS(0,10) *{\bullet}
 \ar @{-} +(8,10)*{\bullet}
 \ar @{-} +(-8,10)*{\bullet}
 \POS(8,20) *{\bullet}
 \POS(-8,20) *{\bullet}
\POS(3,-3) *{1}
\POS(11,23) *{x_1}
\POS(-8,23) *{x_2}
\POS(-2,5) *{e_1}
\POS(7,15) *{e_2}
\POS(-2,16) *{e_3}
\POS(3,10) *{u}
\endxy}\hskip 5pt
\hskip 5pt\vtop{\hbox{$\quad$}\hbox{$\quad$}\hbox{$\buildrel \Phi \over \longmapsto $}\hbox{$\quad$}\hbox{$\quad$}} \ 
\vtop{\hbox{$\quad$}\hbox{$\quad$}\hbox{$\displaystyle{\Big[1-\frac {1}u\Big]*\Big[1-\frac u{x_1}\Big]*
\Big[1-\frac u{x_2}\Big] =\Big[1-\frac {1}u\,,1-\frac u{x_1}\,,1-\frac u{x_2}\Big]\,.}$}}
}

\vskip 5pt
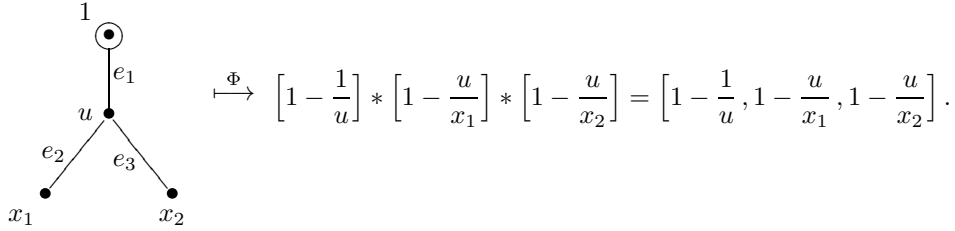
\centerline{\small Figure 5: The forest cycling map on a tree with one internal vertex.}
 \vskip 10pt\noindent
\end{Example}
\noindent
This cycle, as already mentioned, corresponds to the double logarithm $I_{1,1}(x_1,x_2)$, as we will see in Section~6.

\medskip\noindent
{\bf Proof} (of Theorem \ref{dga_morphism}).
The bigrading of a tree $\tau$ in $\CTT^p_n(\Fx)$ translates to the respective bigrading of the cycle $\Phi(\tau)$:
since $\tau$ has $n$ edges, its image under $\Phi$ lands in $\cub_F^n$. Furthermore, $\tau$ has $n+1$ vertices,
$p+1$ of which are external, and so maps to a cycle with $n-p$ independent parameters, i.e., $\Phi(\tau)$ is of codimension $p$.

It remains to check the admissibility of $\Phi(\tau)$ and the compatibility of the differentials.
We prove inductively the following two statements:

\smallskip
{\em Claim 1:} \label{admiss} Each generic $\Fx$-deco tree in $\CTT^p_n(\Fx)$ maps under $\Phi$ to an admissible cycle 
in $\CN_\rmc^{2p-n}(F,p)$.

{\em Claim 2:} The differential on (not necessarily admissible) cycle groups, restricted to $\Phi\big(\CTT^p_n(\Fx)\big)$,
can be written as $\sum_{j=1}^n (-1)^{j-1}\partial_0^j$.

\smallskip
The case $n=1$: Concerning Claim~1, note that for $p=1$ the only obstruction for a cycle $[a]$ in $\CN_\rmc^1(F,1)$ 
to be admissible is that the
(constant) coordinate $a$ is not allowed to be equal to $0$ or to $\infty$.
Translated back to the $\Fx$-deco trees, this condition means that $1-\frac {y_v}{y_w}\not\in \{0,
\infty\}$, where $y_v$ and $y_w$ are the elements in $\Fx$ decorating root and non-root vertex,
respectively. The second possibility ($1-\frac {y_v}{y_w}=\infty$) is excluded since $y_w=0$ is not allowed as a decoration, 
and the first possibility ($1-\frac {y_v}{y_w}=0$)
entails $y_v=y_w$ which is excluded by the genericity condition for $\CTT^1_1(\Fx)$.
For $p\not=1$, there is nothing to show.

Concerning Claim~2, note that the differential for $n=1$ (and any $p$) on both sides turns out to be zero.

\smallskip
The case $n\geq 2$: using induction, we show both claims at level $n$, assuming the validity of both statements for $n-1$.

Concerning Claim~1, note that a tree $\tau$ in $\CTT^p_n(\Fx)$ produces a cycle $\Phi(\tau)$ which satisfies the admissibility condition
everywhere, except
possibly at the zeros and poles of its coordinates. Therefore testing for admissibility amounts
to testing each term in the differential for admissibility separately.

In case we contract an internal edge, we obtain a single tree in which $n$ is reduced by~1, and which is 
still generic since the decoration does not change, so that induction applies. 

In case we contract an external edge, there are at least two trees (or none) produced,
and either tree lies in  $\CTT^p_\nu(\Fx)$ for some $\nu<n$. The admissibility condition
imposes that all decorations for each individual tree are mutually different.
Since decorations are inherited by such a contraction, it is necessary and sufficient for $\tau$
to have itself mutually different decorations in order to make $\Phi(\tau)$ admissible:
otherwise there would be some sequence of contractions of the tree producing an edge with two identical decorations
(in $\Fx$); since this produces under $\Phi$ a coordinate equal to the constant 0, admissibility would be violated.
This proves Claim~1.

\smallskip
Concerning Claim~2, note that each partial differential $\partial_\infty^j$ for a typical coordinate $1-\frac{y_v} 
{y_w}$ in some cycle $Z$ produces a non-empty cycle only if $y_v=\infty$ or $y_w=0$. 
If $y_v$ and $y_w$ are both constants, i.e., lie in $\Fx$, this is not possible.
If, on the other hand, $y_v$ or $y_w$ are variables, 
i.e., correspond to an internal vertex each, 
they can attain those values in $\P^1_F$, but for each such instance the tree structure ensures that there is a further
coordinate in $Z$ of the form $1-\frac x{y_v}$ or $1-\frac{y_w}x$, respectively (with some $x$ either constant or
variable). This will cause the resulting cycle for this partial differential to vanish, since one coordinate
becomes equal to $1$. Thus any $\partial_\infty^j$ produces only the empty cycle, which proves Claim~2.

\smallskip
Now each of the $\partial_0^j$ corresponds to the contraction $\big(\tau/e_j, i_{e_j}(\omega)\big)$ of the $j$th edge $e_j$, 
in the tree differential of the generator $(\tau,\omega)$, where the ordering of edges
corresponds to the ordering of coordinates induced via $\Phi$ from $\omega$. This shows the compatibility
of differentials.
\qed

\section{Tree sums with decomposable boundary and polygons}\label{tree_sums}

In this section, we single out a number of formal linear 
combinations over $R$-deco trees whose boundary consists of decomposable terms only. 
Each such combination can be encoded by an $R$-decorated polygon, which suggests to pull back the DGA structure on trees
to polygons. 

It is convenient to introduce a second differential $\overline \partial$ on polygons. 
Its tree realization  coincides with the main one. 
We relate the Lie coalgebra corresponding to the 
differential $\overline \partial$  to the dihedral Lie coalgebra from \cite{GonDihedral}. 

\subsection{Tree sums with decomposable boundary}
\begin{Definition} \label{alltrivalent}
Let $\{x_1, \dots, x_{m+1}\}$ be a collection of distinct elements of $R$. 
Then $\tau(x_1, \dots, x_{m+1})\in \CT_{2m-1}^m(R)$ is the sum of all trivalent $R$-deco trees with $m$ leaves whose $R$-decoration is
given by $(x_1, x_2,\dots,  x_{m+1})$, the last one decorating the root.  
\end{Definition} 
Recall that the number of such trees is given by the Catalan number $\frac 1{m}{2(m-1) \choose m-1}$.


\begin{Example} \label{triplelog}
1. For $m=2$, the tree $\tau(x_1,x_2,1)$ is given by the tree in Figure~5.

2. For $m=3$, the sum of trees has already more than a single  term: 

\vskip 10pt
\hskip 5pt
 \vtop{\hbox{$\quad$}\hbox{$\quad$}\hbox{$\quad$}\hbox{$\tau(x_1,x_2,x_3,x_4)\ =\ $}\hbox{$\quad$}\hbox{$\quad$}}
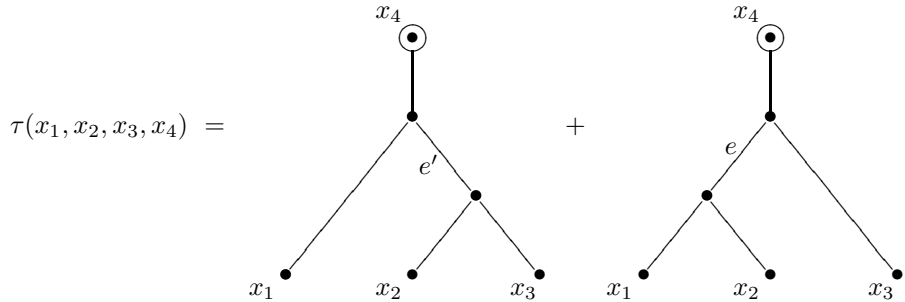
\hbox{
 \vbox{\xy  0;<-3pt,0pt>: 
 \POS(0,0) *+{\bullet} *\cir{} 
 \ar @{-} +(0,10)*{\bullet}
  \POS(0,10) *{\bullet}
 \ar @{-} +(16,20)*{\bullet}
 \ar @{-} +(-8,10)*{\bullet}
 \POS(16,30) *{\bullet}
 \POS(-8,20) *{\bullet}
 \ar @{-} +(8,10)*{\bullet}
 \ar @{-} +(-8,10)*{\bullet}
\POS(3,-3) *{x_4}
\POS(19,32) *{x_1}
\POS(3,32) *{x_2}
\POS(-14,32) *{x_3}
\POS(-2,16) *{e'}
\endxy}\hskip 5pt}
 \vtop{\hbox{$\quad$}\hbox{$\quad$}\hbox{$\quad$}\hbox{$+$}\hbox{$\quad$}\hbox{$\quad$}}
\hbox{
 \vbox{\xy  0;<-3pt,0pt>: 
 \POS(0,0) *+{\bullet} *\cir{} 
 \ar @{-} +(0,10)*{\bullet}
  \POS(0,10) *{\bullet}
 \ar @{-} +(8,10)*{\bullet}
 \ar @{-} +(-16,20)*{\bullet}
 \POS(8,20) *{\bullet}
 \ar @{-} +(8,10)*{\bullet}
 \ar @{-} +(-8,10)*{\bullet}
\POS(3,-3) *{x_4}
\POS(3,32) *{x_2}
\POS(19,32) *{x_1}
\POS(-14,32) *{x_3}
\POS(5,14) *{e}
\endxy}\hskip 5pt}
\vskip 8pt

\centerline{\small Figure 6: The sum of trees corresponding to the quadruple  $(x_1, x_2, x_3, x_4)$.}
\vskip 8pt
\end{Example}

Applying $\Phi$ to the sum of trees $\tau(x_1,x_2,x_3,1)$ with mutually different $x_i\in \Fx\setminus\{1\}$ for a field $F$, 
we get the following cycle:
\begin{equation} \label{Z3}
Z_{x_1,x_2,x_3} = \Big[1-\frac 1 t,1-\frac t{x_1},1-\frac t u,1-\frac u{x_2} ,1-\frac u {x_3}\Big]
\ + \ \Big[1-\frac 1 t, 1-\frac t u, 1-\frac u{x_1}, 1-\frac u{x_2} ,1-\frac t {x_3}\Big]. \,
\end{equation}
Here we have two parametrizing variables, $t$ and $u$, and $Z_{x_1,x_2,x_3}\in \CN_\rmc^1(F,3)$. 

\begin{Proposition} \label{decomp_boundary}
Let $\{x_1, \dots, x_m\}$ be a collection of distinct elements of $\Fx$. 
The image of the tree sum $\tau(x_1, \dots, x_m)$ under the forest cycling map is an admissible cycle with decomposable boundary.
\end{Proposition}

\noindent {\bf Proof.} The admissibility of the resulting cycle follows immediately from Lemma \ref{admiss}, since the
decorations $x_1,\dots,x_m$ are mutually different.

Concerning the statement about the boundary, note first that the contraction of an {\em external edge} in a trivalent tree produces 
a forest which, if non-empty, is not a tree, and thus maps to a decomposable cycle under $\Phi$. 

The contributions of {\em internal edges} in the tree sum cancel.
Indeed, there are exactly two trivalent trees, $\tau$ and $\tau'$, 
containing the same subtrees $\tau_i$, $i=1,\ldots,4$, but joined
in a different way, via edges $e$ and $e'$, as shown in the picture: 
\vskip 10pt \hskip 100pt
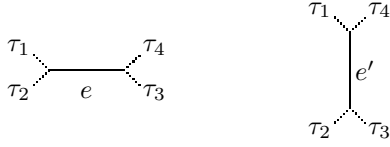
\xy   \POS(0,5) \ar@{-} +(10,0)\ar@{..} +(-2,2) \ar@{..} +(-2,-2) \POS(5,2) *{e}  
\POS(-4,8) *{\tau_1} \POS(-4,2) *{\tau_2}
\POS(10,5)\ar@{..} +(2,2) \ar@{..} +(2,-2) \POS(14,8) *{\tau_4} \POS(14,2) *{\tau_3}
\POS(40,0) \ar@{-} +(0,10)\ar@{..} +(2,-2) \ar@{..} +(-2,-2) \POS(42,5) *{e'}  \POS(36,-3) 
*{\tau_2} \POS(44,-3) *{\tau_3}
\POS(40,10)\ar@{..} +(-2,2) \ar@{..} +(2,2) \POS(36,13) *{\tau_1} \POS(44,13) *{\tau_4}
\endxy
\vskip 10pt
We can  write the canonical orientations $\omega(\tau)$ of $\tau$ and $\omega(\tau')$ of $\tau'$ 
in terms of the subtrees $\tau_i$ ($i=1,\dots,4$), where $\tau_1$ contains the root of $\tau$, as follows: 
$$
\omega_{\tau}=\omega_{\tau_1}\,\wedge\,\omega_{\tau_2}\,\wedge\, e\,\wedge\, \omega_{\tau_3}\,\wedge\,\omega_{\tau_4}, \qquad 
\omega_{\tau'}=\omega_{\tau_1}\,\wedge\, e\,\wedge\,\omega_{\tau_2}\,\wedge\, \omega_{\tau_3}\,\wedge\,\omega_{\tau_4}
$$
 where $\omega_{\tau_i}$
denotes the canonical orientation for $\tau_i$. 
The essential point is that each of the trees $\tau_i$ have an odd number of edges, and thus the corresponding orientations anticommute. 
Thus $\omega_{\tau'}=
-\,\omega_\tau$. 
Therefore the terms in the differential corresponding to the contraction of $e$ in $\tau$ and of $e'$ in $\tau'$ produce the same tree with
opposite orientation torsor and thus cancel. \qed

\subsection{The algebra of $R$-deco polygons}\label{polygonalgebra}
A polygon is an {\bf oriented}
 convex polygon with $N\geq2$ sides. It inherits 
a cyclic order of the sides.
Let $R$ be a set.

\begin{Definition} 
An {\bf $R$-deco polygon} is a polygon whose sides are decorated by elements of $R$, with a distinguished side, called its
{\bf root side}.
\end{Definition}

The orientation of the polygon induces a linear  order of its sides such that the root side is the last one. 
The {\bf first vertex} of the polygon is the common vertex of the root side and the first side.
Once the root side is determined, the orientation determines the first vertex, and vice versa.
So we indicate the orientation in pictures by marking the first vertex by a bullet.

We use the notation $\lceilx a_1,\dots,a_N\rfloorx$  for an $R$-deco $N$-gon 
with linearly ordered sides decorated by $a_1, \ldots, a_N$, where $a_N$ is the root decoration.

\begin{Example} Here is an $R$-deco 6-gon $\pi = \lceilx a_1,\dots,a_{6}\rfloorx$.
The root side is drawn by a double line, the first vertex is marked by a bullet, and 
the orientation is counterclockwise. 

\hskip 120pt
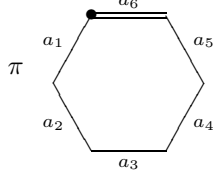
\xy
\POS(0,0)  \ar@{-} +(10,0)_{a_3} \ar@{-} +(-5,9)^{a_2} 
\POS(10,0) \ar@{-} +(5,9)_{a_4}
\POS(15,9) \ar@{-} +(-5,9)_{a_5}
\POS(10,18) \ar@{=} +(-10,0)_{a_6}
\POS(0,18) *+{\bullet}
\POS(-5,9) \ar@{-} +(5,9)^{a_1}\POS(-10,11) *{\pi}

\endxy

\end{Example}
The {\bf weight} $\chi(\pi)$ of an $R$-deco polygon $\pi$ is the number of its non-root sides.

We define the graded vector space $V^\pg_\smallbullet(R)=\bigoplus_{n=0}^\infty V^\pg_{n}(R)$, 
where $V^\pg_n(R)$, 
for $n\geq 1$, denotes the vector space of $R$-deco ($n+1$)-gons, and 
$V^\pg_0(R):= \Q$.

\begin{Definition}
The {\bf polygon algebra} $\CP_\smallbulletP^{(\smallbullet)}=\CP_\smallbulletP^{(\smallbullet)}(R)$
is the exterior algebra of the graded vector space $V^\pg_\smallbullet(R)$.  It is bigraded: the 
second grading, denoted via superscripts, comes from the exterior power.
\end{Definition}
As we will often ignore the second grading, 
we mostly denote this algebra by $\CP_\smallbulletP$. 

To define a differential on $\CP_\smallbulletP$, we
 introduce arrows and dissections:

\begin{Definition}
An {\bf arrow in an $R$-deco polygon} $\pi$ is a line segment beginning at a vertex and ending at the interior of a side of $\pi$. An arrow connecting a vertex with an incident side is trivial.

An arrow in  $\pi$ is a {\bf backward arrow} if, 
in the linear order of the sides of $\pi$, its  end
is before its beginning. Otherwise, it is a 
{\bf forward arrow}.
\end{Definition}

We identify arrows with the same vertex ending at the same side. 


\begin{Example} 
The following picture shows two forward arrows $\alpha$ and $\gamma$ 
and a backward arrow $\beta$ in the $R$-deco 6-gon of the previous example. 

\hskip 90pt
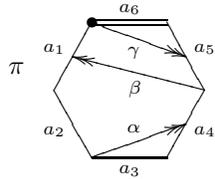
\xy
\POS(0,0)  \ar@{-} +(10,0)_{a_3} \ar@{-} +(-5,9)^{a_2}  \ar@{->>} +(12.5,4.5)^\alpha
\POS(15,9) \ar@{->>} -(17.5,-4.5)^\beta
\POS(0,18) \ar@{->>} +(12.5,-4.5)_\gamma
\POS(10,0) \ar@{-} +(5,9)_{a_4}
\POS(15,9) \ar@{-} +(-5,9)_{a_5}
\POS(10,18) \ar@{=} +(-10,0)_{a_6}
\POS(0,18) *+{\bullet}
\POS(-5,9) \ar@{-} +(5,9)^{a_1} \POS(-10,12) *{\pi}

\endxy
\vskip10pt
\end{Example}

Any non-trivial arrow dissects a polygon into two regions. 

\def\rt{{\bullet\hskip-2pt=}}
\def\rst{\disjoint}

\begin{Definition} \label{defposreg} Let $\pi$ be an $R$-deco polygon with $\geq3$ sides, and $\alpha$ 
a non-trivial arrow in $\pi$.
A {\bf dissection} of $\pi$ by $\alpha$ 
is the ordered pair
$(\rho^\rt_\alpha ,\rho^\rst_\alpha)$ of a {\bf root region} $\rho^\rt_\alpha$ 
and a {\bf cut-off region} 
$\rho^\rst_\alpha$, where 
$\rho^\rt_\alpha\cup \rho^\rst_\alpha=\pi$. 
The root region is the one containing the first vertex and part of the root side.

We assign to these regions $R$-deco polygons $\pi^\rt_\alpha$ and $\pi^\rst_\alpha$ by contracting 
the arrow $\alpha$ in either region.
The polygon $\pi^\rt_\alpha$ inherits the root and the orientation from $\pi$. 
In the polygon $\pi^\rst_\alpha$, the root side is the
one on which $\alpha$ ends. 
Its orientation is inherited from the one of $\pi$ 
if the dissecting arrow is forward, and is the opposite one otherwise. 
\end{Definition}

\begin{Example} The arrow $\alpha$ dissects the hexagon $\pi$ in the picture into the regions $\rho^\rt_\alpha$ and 
$\rho^\rst_\alpha$. On the right hand side, we draw the $R$-deco polygons associated to the two regions.

\vskip 10pt
\hskip 30pt
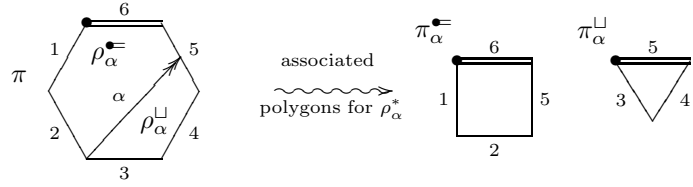
\xy
\POS(0,0)  \ar@{-} +(10,0)_3  \ar@{-} +(-5,9)^2  \ar@{->>} +(12.5,13.5)^\alpha
\POS(10,0) \ar@{-} +(5,9)_4
\POS(15,9) \ar@{-} +(-5,9)_5
\POS(10,18) \ar@{=} +(-10,0)_6
\POS(0,18) *+{\bullet}
\POS(-5,9) \ar@{-} +(5,9)^1
\POS(3,14) *+{\rho^\rt_\alpha}
\POS(9,5) *+{\rho^\rst_\alpha}
\POS(-9,11) *+{\pi}
\hskip 70pt  \POS(0,9) \ar@{~>}^{\text{associated}\phantom{\Big|}}_{\text{polygons for }\rho_\alpha^*} +(16,0) \hskip 70pt
\POS(0,13) *+{\bullet}
\POS(10,13) \ar@{=} +(-10,0)_6
\ar@{-} +(0,-10)^5
\POS(10,3) \ar@{-} +(-10,0)^2
\POS(0,13) \ar@{-} +(0,-10)_1
\POS(-3,17) *+{\pi^\rt_\alpha}
\hskip 60pt
\POS(0,13) *+{\bullet}
\POS(10,13) \ar@{=} +(-10,0)_5
\ar@{-} +(-5,-8)^4
\POS(0,13) \ar@{-} +(5,-8)_3
\POS(-3,17) *+{\pi^\rst_\alpha}

\endxy

\end{Example}

\vskip 10pt Let 
$\chi(\alpha) := \chi(\pi_{\alpha}^\rst)$ be the number of the non-root edges of the polygon $\pi_{\alpha}^\rst$. 
\begin{Definition} \label{sign}
To an arrow $\alpha$ in an $R$-deco polygon $\pi$ we assign its sign:
$$\sgn(\alpha) = \begin{cases} (-1)^{\chi(\alpha)} &\text{if }\alpha\text{ is backward,}\\ 1&\text{otherwise.}
                 \end{cases}
$$
\end{Definition}

The sum over all the dissections of $\pi$ provides a differential $\partial\,\pi$:

\begin{Definition} \label{ggt}
Given a non-trivial arrow $\alpha$ of an $R$-deco polygon $\pi$,  set 
$$\partial_\alpha \pi = \pi^\rt_\alpha\wedge \pi^\rst_\alpha \in \CP_\smallbulletP^{(2)},$$ 
where $\pi^\rt_\alpha$ and $\pi^\rst_\alpha$ are as in Definition \ref{defposreg}.

The {\bf differential}  $\partial\, \pi$ of $\pi$ is the sum $\sum_\alpha \sgn(\alpha) \partial_\alpha \pi$ 
over the arrows $\alpha$ in $\pi$.

It extends to a {\bf differential on} $\CP_\smallbulletP^{(\smallbullet)}$ via
the Leibniz rule $\partial(a\wedge b)=(\partial a)\wedge b \, +(-1)^{\deg(a)}\, a\wedge (\partial b)$. Here $\deg$ is the degree in the exterior algebra.

\end{Definition}

\begin{Example} \label{diffpol}
The differential applied to a triangle gives three terms:
$$\partial\Big(\threegon{1}{2}{3}
\Big) = \twogon{1}{3}\wedge  \twogon{2}{3}\ + \
\twogon{2}{3} \wedge  \twogon{1}{2}\ -\ \twogon{1}{3} \wedge \twogonright{2}{1}  \,.$$

The summands correspond respectively to the following three dissections:
$$\threegonarrowc{1}{2}{3}\quad ,\quad \threegonarrowb{1}{2}{3}\quad ,\quad \threegonarrowa{1}{2}{3}\,.$$
\end{Example}

\vskip 10pt

\smallskip
The justification of the name ``differential'' is given by  Proposition \ref{just} below. 

\vskip 3mm

{\bf Another differential on $R$-deco polygons}. Given an arrow $\alpha$ in an $R$-deco polygon, set 
$$\varepsilon_\alpha = \begin{cases} -1 &\text{if }\alpha\text{ is backward,}\\ 1&\text{otherwise.}
                 \end{cases}
$$
\begin{Definition} \label{tilde}
In the set-up of  Definition \ref{defposreg}, a dissecting arrow $\alpha$ gives rise 
to two 
$R$-deco polygons $\overline \pi^\rt_\alpha$ and $\overline \pi^\rst_\alpha$, 
whose orientations are inherited from the one of $\pi$, and whose decorations are as in Definition \ref{defposreg}.

A differential $\overline \partial$ is defined as in Definition \ref{ggt}, with
$\overline \partial \pi= \sum_\alpha\varepsilon_\alpha \overline \pi^\rt_\alpha \wedge \overline \pi^\rst_\alpha$.
\end{Definition}

Note that $\overline \pi^\rt_\alpha = \pi^\rt_\alpha$, and that 
 $\overline \pi^\rst_\alpha = \pi^\rst_\alpha$ for a forward arrow $\alpha$. If $\alpha$ is backward, then $\overline \pi^\rst_\alpha$ has the opposite orientation.

\begin{Proposition} \label{just} 1. One has $\partial^2=0$. The polygon algebra $\CP_\smallbullet^{(\smallbullet)}(R)$ is an augmented 
 DGA 
with the differential  $\partial$. 

2. One has $\overline \partial^2=0$, providing the same polygon algebra with another augmented DGA structure, 
denoted $\overline \CP_\smallbullet^{(\smallbullet)}(R)$. 
\end{Proposition}

 \begin{Remark} 
Proposition \ref{just} is equivalent to the following statement: 
Each of the maps  $\overline \partial$ and $\partial$ provides a Lie cobracket 
\begin{equation} \label{cobracket}
V_{\smallbullet}^\pg(R) \to V_{\smallbullet}^\pg(R) \wedge V_{\smallbullet}^\pg(R)
\end{equation}
So  there are two graded Lie coalgebra structures on the graded vector space ${V}^\pg_{\smallbullet}(R)$. 

The graded-commutative DGA $\overline {\cal P}_{\smallbullet}^{(\smallbullet)}(R)$ (respectively ${\cal P}_{\smallbullet}^{(\smallbullet)}(R)$)
is just the standard cochain 
complex of the Lie coalgebra $(V_{\smallbullet}^\pg(R), \overline \partial)$ (respectively $(V_{\smallbullet}^\pg(R), \partial)$). 
\end{Remark}

We will give a direct proof of Proposition \ref{just} in Section 6. 

On the other hand, there is another proof. Namely, we identify in 
Proposition \ref{hihiks} the map $\partial $  with 
the cobracket in a Lie coalgebra of formal iterated integrals defined in \cite{GonFund}. 
The latter Lie coalgebra was interpreted in loc. cit. as the dual to the Lie algebra of all derivations of a certain structure. 
So the property $\partial^2=0$ is valid on the nose. 
Similarly,  we identify in Subsection 5.4 the map $\overline \partial $ with 
the cobracket in a rooted version of one of the dihedral Lie coalgebras defined in \cite{GonDihedral}.

\subsection{Relating the polygon algebra to the tree algebra}
The standard duality between triangulated polygons and trivalent trees carries over to the respective
$R$-deco objects (to each side we first associate a region outside of the polygon, and then each of its sides and diagonals 
corresponds to an edge of the tree, while each region---inside or outside---corresponds to a vertex):

\begin{Example} We give a picture of an $\N$-deco hexagon with triangulation, on the left (the dotted lines indicate the 
outside regions of the polygon), and of its dual $\N$-deco tree, on the right. Note the little arrow
attached to the root edge which indicates the orientation (counterclockwise for the given example).

\vskip10pt
\hskip 15pt
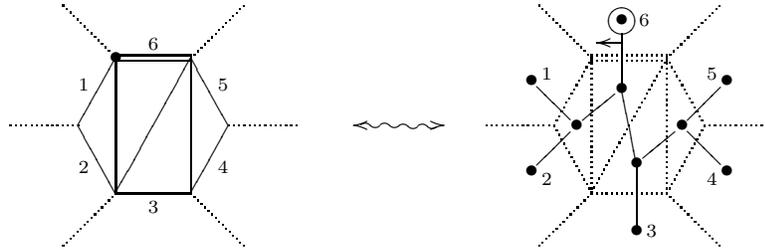
\xy 
\POS(0,0)  \ar@{-} +(10,0)_{3}  \ar@{-} +(-5,9)^{2}   \ar@{-} +(0,18)
\POS(10,0) \ar@{-} +(5,9)_4  \ar@{-} +(0,18) \ar@{.} +(7,-7) 
\POS(15,9) \ar@{-} +(-5,9)_5 \ar@{.} +(9,0)
\POS(10,18) \ar@{=} +(-10,0)_{6} \ar@{.} +(7,7) 
\POS(-5,9) \ar@{-} +(5,9)^1 \ar@{.} +(-9,0)
\POS(0,0) \ar@{-} +(10,18) \ar@{.} +(-7,-7)
\POS(0,18) \ar@{.} +(-7,7)  \POS(0,18)*+{\bullet}

\hskip 90pt 
\POS(0,9) \ar@{<~>} +(12,0)\hskip 60pt

\hskip 30pt
\POS(0,0)  \ar@{.} +(10,0)  \ar@{.} +(-5,9)   \ar@{.} +(0,18)
\POS(10,0) \ar@{.} +(5,9) \ar@{.} +(0,18) \ar@{.} +(7,-7) 
\POS(15,9) \ar@{.} +(-5,9) \ar@{.} +(9,0)
\POS(10,18) \ar@{:} +(-10,0) \ar@{.} +(7,7) 
\POS(-5,9) \ar@{.} +(5,9) \ar@{.} +(-9,0)
\POS(0,0) \ar@{.} +(10,18) \ar@{.} +(-7,-7)
\POS(0,18) \ar@{.} +(-7,7)

\POS(-2,9) *{\bullet} \ar@{-} +(-6,6)  *{\bullet} \ar@{-} +(-6,-6)  *{\bullet} \ar@{-} +(6,5)  *{\bullet}
\POS(12,9) *{\bullet} \ar@{-} +(6,-6)  *{\bullet} \ar@{-} +(6,6)  *{\bullet}\ar@{-} +(-6,-5)  *{\bullet}
\POS(4,14) \ar@{-} +(0,9) *+{\bullet}*\cir{}
\POS(4,20) \ar@{->} +(-3.3,0)   
\POS(6,4) \ar@{-} +(0,-9)  *{\bullet} \ar@{-} +(-2,10)
\POS(8,-5) *{\scriptstyle 3}
\POS(16,2) *{\scriptstyle 4}
\POS(16,16) *{\scriptstyle 5}
\POS(-6,2) *{\scriptstyle 2}
\POS(-6,16) *{\scriptstyle 1}
\POS(7,23) *{\scriptstyle 6}

\endxy
\end{Example}

\vskip 15pt

\begin{Proposition} \label{inj}
The standard duality induces a map from $\CP_\smallbulletP(R)$ to the subalgebra $\CT_{\rm 3-val}(R)$ of the tree algebra $\CT_\smallbullet(R)$ which is generated by {\em trivalent} 
trees. 
\end{Proposition}
\noindent{\bf Proof.} 
Let $\nabla$ be a triangulation of an $R$-deco polygon $\pi$ with dual tree $\tau$.
For each internal vertex of $\tau$, the orientation of $\pi$ induces a cyclic ordering on the incident edges, 
and the root of $\pi$ induces a root of $\tau$. Putting both data together, $\tau$
inherits an orientation (in fact, even a linear ordering on its edges).

We associate to $\pi$ the following sum, where $\nabla$ runs through all triangulations of $\pi$ and 
$\tau(\pi,\nabla)$ denotes the (oriented) 
$R$-deco tree dual to the triangulation $\nabla$ of $\pi$,
\begin{equation}\label{nabla}
\pi \mapsto \sum_{\nabla\text{ triangulation}} \tau(\pi,\nabla)\,.
\end{equation}
By the duality above, we identify each term in the sum with the corresponding oriented $R$-deco tree. \qed


\begin{Theorem} The map from (\ref{nabla}) induces a map of differential graded algebras $\CP_\smallbulletP(R) \to 
\CT_{\smallbullet}(R)$ as well as $\overline \CP_\smallbulletP(R) \to 
\CT_{\smallbullet}(R)$. 
\end{Theorem}
\noindent{\bf Proof.} We start from the second claim, which is a bit more straightforward. 
Let $\pi$ be an $R$-deco polygon with triangulation $\nabla$, 
and let $\tau=\tau(\pi,\nabla)$ be the dual $R$-deco tree.
Let us calculate the effect of contracting an edge $e$. 

{\em Internal edges:} The proof of Proposition \ref{decomp_boundary} shows that
the terms in the differential for $\tau$ arising from internal edges cancel.

{\em External edges:} Let $e$ be an external edge. Recall that 
the dual tree for a triangulation of an oriented polygon has a canonical orientation 
(that is, an element 
of the orientation torsor of the tree).

\begin{Remark} \label{s}
Changing the orientation of an $N$-gon $\pi$ amounts to changing the corresponding orientation of the dual 
tree of any triangulation of $\pi$ by $(-1)^{N}$. (Flip the two branches at each of the $N-2$
internal vertices.)
\end{Remark}

The edge $e$ is dual to a side $s_e$ of the polygon $\pi$. 
Let us cut out the triangle $t_e$ of the triangulation containing the side $s_e$ 
from the triangulation of the polygon $\pi$. We obtain two polygons, $\pi_R$ and $\pi_O$, where 
$\pi_R$ is the one containing the root side if $s_e$ is not a root, and
 the first vertex in the latter case. Let  $\omega_R$ and $\omega_O$ be 
the canonical orientations of the dual trees assigned to the induced triangulations of these polygons. 

Let $\alpha$ be the dissecting arrow corresponding to the triangle $t_e$: it ends at the side $s_e$, 
and its vertex coincides with the one of $t_e$ opposite to the side $s_e$. 
Then, if $\alpha$ is a forward (respectively a backward) arrow, 
the canonical orientation of the tree $\tau$ is 
 $\omega_R \wedge \omega_O\wedge e = e  \wedge \omega_R \wedge \omega_O$ (respectively $\omega_R \wedge e \wedge \omega_O = -e \wedge \omega_R \wedge \omega_O$). 
Here we used the fact that the dual tree always has an odd number of edges. 
In particular, the root edge corresponds to a forward arrow and produces the orientation $e\wedge \omega_R \wedge \omega_O$. 
In each case this matches the definition $\varepsilon_\alpha$. 

A proof of the first claim is deduced from this and Remark \ref{s}. Indeed, this is obvious for a forward dissecting arrow. 
For a backward dissecting 
arrow $\alpha$ we have a different orientation of the cut-off polygon, which thanks to Remark \ref{s} amounts to 
the sign $(-1)^{\chi(\alpha)+1}$, and $\sgn(\alpha) = (-1)^{\chi(\alpha)}$. So the total sign difference will be  
$-1$, matching $\varepsilon_\alpha$. 
\qed

\bigskip
An $R$-deco polygon is {\bf generic} if its decorations are all different. Now we can 
connect the polygons to algebraic cycles.

\begin{Corollary}
There is a map of differential graded algebras from the generic part $\widetilde{\CP}_\smallbullet(\Fx)$ of $\CP_\smallbulletP(\Fx)$
to the cubical algebraic cycles $\CN_\rmc^{2*-\smallbullet}(*)$ over $F$.

More precisely, $\widetilde{\CP}_{N+1}(\Fx)$ gets mapped to $\CN_\rmc^{2*-N}(*)$.

\end{Corollary}
\noindent{\bf Proof.} We just need to combine the above theorem with Theorem \ref{dga_morphism}. \qed

\subsection{Appendix: Relation with the dihedral Lie coalgebra from \cite{GonDihedral}} \label{dihedral}
Let us identify 
the graded Lie coalgebra of $R$-deco polygons $({V}^\pg_{\smallbullet}(R), \overline \partial)$ with a rooted version of  one of the Lie coalgebras 
introduced in \cite{GonDihedral} under the name of dihedral Lie coalgebras. 
Let $G$ be a commutative group. Recall the bigraded dihedral Lie coalgebra 
$\widetilde {\cal D}_{\smallbullet, \smallbullet}(G) = \oplus_{w,m}\widetilde {\cal D}_{w,m}(G)$ defined in Section 4.5 of loc.~cit. 
Consider the corresponding graded Lie sub-coalgebra $\widetilde {\cal D}_{\smallbullet}(G)$:
\begin{equation} \label{dih}
\widetilde {\cal D}_{\smallbullet}(G):= \bigoplus_{w>0}\widetilde {\cal D}_{w,w}(G), 
\end{equation} 
called the diagonal part of the dihedral Lie coalgebra $\widetilde {\cal D}_{\smallbullet, \smallbullet}(G)$. 
We relate it to ${V}^{\pg}_{\smallbullet}(R)$ as follows. 
Put $R=G$. Then the group $G$ acts on ${V}^{\pg}_{\smallbullet}(G)$ by simultaneous multiplication on the $G$-decoration, so we can 
pass to the Lie coalgebra of $G$-coinvariants. Killing 
in the latter coalgebra, for each integer $n>0$,  the 
generator $(\underbrace{e: \dots :e}_{n+1})$, where $e$ is the unit in $G$,  we get 
the Lie coalgebra (\ref{dih}).

\medskip
{\bf The correspondence on pictures}. In \cite{GonDihedral} a 
typical generator of (\ref{dih}) is given as a 
cyclically ordered set
of $n$ ``black'' points  on the circle. The complement to these points is the union of $n$ arcs. 
The cobracket is the sum of several terms, each of which is
obtained via the following three step procedure: 

\begin{enumerate}

\item Cut the circle at a black  point and at an arc 
 not adjacent to that point.

\item Complete each of the two resulting half-circles
to a circle by gluing its ends. 

\item Take the wedge product of the above two circles, so that in the
left factor the path from the black point cut to the arc cut follows the cyclic order.
\end{enumerate}

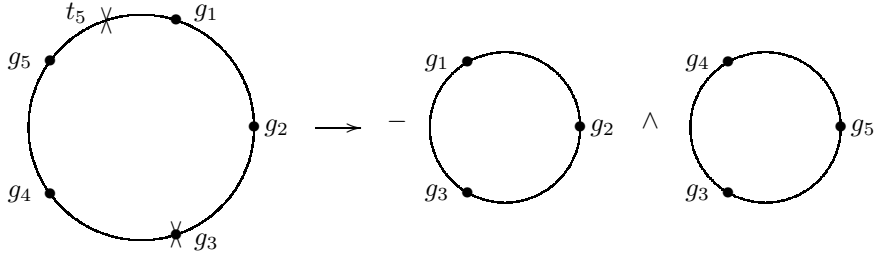
\xyoption{arc}

\vskip10pt
\hskip -5pt\xy
*\xycircle(15,15){} 
\POS(4.63,14.26) *{\bullet} \POS(8.63,15.26)   *{g_1}
\POS(-12.13,8.82) *{\bullet} \POS(-16.13,8.82) *{g_5}
\POS(-12.13,-8.82) *{\bullet} \POS(-16.13,-8.82)*{g_4}
\POS(4.63,-14.26) *{\bullet}*{\setminus} *{/} \POS(8.63,-15.26)*{g_3} 
\POS(15,0) *{\bullet}\POS(18,0)*{g_2}
\POS(-4.63,14.26)  *{\setminus} *{/} \POS(-8.63,15.26) *{t_5}
 \endxy 
 \xy\hskip 10pt \POS(0,0) \ar@{->} +(6,0) \endxy \hskip 20pt $-\ $
\xy  
*\xycircle(10,10){}
\POS(-5.0,8.66) *{\bullet} \POS(-9.0,8.66) *{g_1} 
\POS(-5.0,-8.66) *{\bullet} \POS(-9.0,-8.66) *{g_3} 
\POS(10,0) *{\bullet} \POS(13,0) *{g_2} 
 \endxy \hskip 10pt$\wedge$\hskip 10pt
\xy
*\xycircle(10,10){}
\POS(-5.0,8.66) *{\bullet} \POS(-9.0,8.66) *{g_4} 
\POS(-5.0,-8.66) *{\bullet} \POS(-9.0,-8.66) *{g_3} 
\POS(10,0) *{\bullet} \POS(13,0) *{g_5} 
\endxy
\vskip10pt\noindent

The corresponding term in the differential on polygons looks as follows: 
we assign to the circle with black points a polygon whose vertices are centers of the arcs, and whose 
sides are labelled by the black points. 
Then the cutting of the circle corresponds to an arrow from the middle of an arc 
to the respective black point.
The dissection induced by this arrow gives a term in the differential and corresponds precisely to the
wedge product of the two smaller circles above. 

\vskip10pt
\hskip 5pt
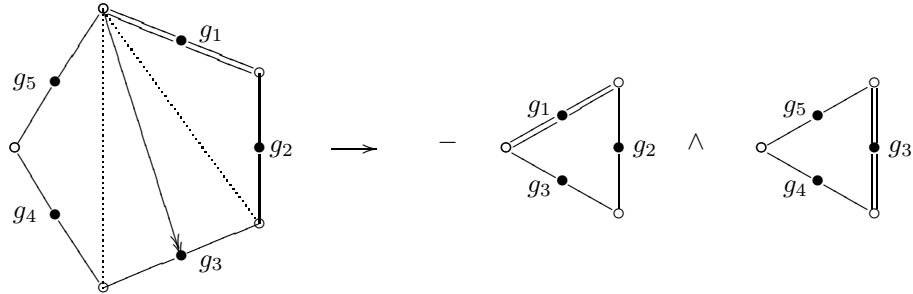
\xy {}
\POS(4.63,14.26) *{\bullet}  \ar@{=} +(10.37,-4.26) *{\circ} \ar@{=} +(-10.37,4.26) *{\circ}\POS(8.63,15.26)   *{g_1}
\POS(-12.13,8.82) *{\bullet} \ar@{-} +(6.39,9.7)*{\circ} \ar@{-} +(-5.39,-8.82)*{\circ}
\POS(-5.74,18.52) \ar@{->>} +(10.37,-32.78) \ar@{.} +(0,-37.04)\ar@{.} +(20.74,-28.52) *{\circ}
 \POS(-16.13,8.82) *{g_5}
\POS(-12.13,-8.82) *{\bullet} \ar@{-} +(6.39,-9.7)*{\circ} \ar@{-} +(-5.39,8.82)*{\circ}\POS(-16.13,-8.82)*{g_4}
\POS(4.63,-14.26) *{\bullet}
\ar@{-} +(10.37,4.26) \ar@{-} +(-10.37,-4.26)
\POS(8.63,-15.26)*{g_3} 
\POS(15,0) *{\bullet} \ar@{-} +(0,10) \ar@{-} +(0,-10) \POS(18,0)*{g_2}
\hskip 70pt \POS(0,0) \ar@{->} +(6,0) \hskip 40pt $-\ $
\hskip 40pt
\POS(5,0) *{\bullet}\ar@{-}+(0,-8.) \ar@{-}+(0,8.) \POS(8.5,0) *{g_2} 
\POS(-2.5,4.33) *{\bullet}\ar@{=}+(7.5,4.33)*{\circ} \ar@{=}-(7.5,4.33) *{\circ}
\POS(-5.5,5.33) *{g_1} 
\POS(-2.5,-4.33) *{\bullet} \ar@{-}+(7.5,-4.33)*{\circ} \ar@{-}-(7.5,-4.33) *{\circ}\POS(-5.5,-5.33) *{g_3} 
\hskip 40pt$\wedge$\hskip 20pt
\hskip 30pt
\POS(5,0) *{\bullet}\ar@{=}+(0,-8.) \ar@{=}+(0,8.) \POS(8.5,0) *{g_3} 
\POS(-2.5,4.33) *{\bullet}\ar@{-}+(7.5,4.33)*{\circ} \ar@{-}-(7.5,4.33) *{\circ}
\POS(-5.5,5.33) *{g_5} 
\POS(-2.5,-4.33) *{\bullet} \ar@{-}+(7.5,-4.33)*{\circ} \ar@{-}-(7.5,-4.33) *{\circ}\POS(-5.5,-5.33) *{g_4} 
\endxy
\vskip 3mm
We will see in Section 8 that the other graded Lie coalgebra $({V}^\pg_{\smallbullet}(R), \partial)$ 
is related to the Lie coalgebra of motivic iterated integrals.

\section{Cycles in the bar construction of the polygon algebra} \label{bar}
\subsection{Higher dissections of polygons}
To deal with the bar construction, we introduce $n$-fold dissections.

\begin{Definition}
We say that two arrows in a polygon {\bf do not intersect} if they have no point in common, except possibly the starting point.
\end{Definition}

\begin{Example} Two intersecting arrows on the left, four non-intersecting arrows on the right.
\end{Example}

\hskip 50pt
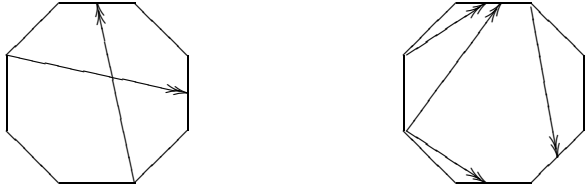
\xy (0,0) 
\POS(0,0)  \ar@{-} +(10,0)
\POS(10,0) \ar@{-} +(7,7)
\POS(17,7) \ar@{-} +(0,10)
\POS(17,17) \ar@{-} +(-7,7)
\POS(10,24) \ar@{-} +(-10,0)
\POS(0,24) \ar@{-} +(-7,-7)
\POS(-7,17) \ar@{-} +(0,-10)
\POS(-7,7) \ar@{-} +(7,-7)

\POS(5,24) \ar@{<<-} +(5,-24)
\POS(17,12)\ar@{<<-} +(-24,5)

\hskip 150pt 

\POS(0,0)  \ar@{-} +(10,0)
\POS(10,0) \ar@{-} +(7,7)
\POS(17,7) \ar@{-} +(0,10)
\POS(17,17) \ar@{-} +(-7,7)
\POS(10,24) \ar@{-} +(-10,0)
\POS(0,24) \ar@{-} +(-7,-7)
\POS(-7,17) \ar@{-} +(0,-10)
\POS(-7,7) \ar@{-} +(7,-7)

\POS(4,24)  \ar@{<<-} +(-10.5,-7.)
\POS(6,24)  \ar@{<<-} +(-12.5,-17.)
\POS(13.5,3.5) \ar@{<<-} +(-3.5,20)

\POS(4,0)  \ar@{<<-} +(-10.5,7.)

\endxy
\vskip 10pt

\begin{Definition}
Let $\pi$ be an $R$-deco polygon and $n$ a positive integer. An {\bf $n$-fold dissection $D$ of $\pi$} is the 
set of $n$ 
regions provided by  $n-1$ non-intersecting arrows
in $\pi$: the root region of $D$ is defined in the same way as for 2-fold dissections. 
\end{Definition}

There is a partial order on the set of regions of a dissection, induced by the root region and by 
the adjacency of regions. A 1-fold dissection of $\pi$ is equal to $\pi$ itself.

An  $n$-fold dissection $D$ of an $R$-deco polygon $\pi$ gives rise to its {\bf dual tree} $\tau(D)$ whose 
vertices 
correspond to
the regions of the dissection, and whose $(n-1)$ edges are dual to the arrows of $D$. 
The root region 
of $D$ determines the root of $\tau(D)$.

\medskip
\begin{Example} \label{dual_tree_example} 
Here is a polygon with a 5-fold dissection, and its dual tree.
\end{Example}
\hskip 120pt
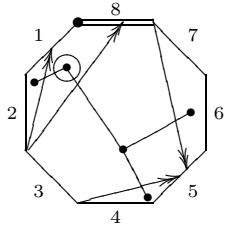
\xy (0,0) 
\POS(-1.5,18) *+{\scriptstyle{\bullet}} *{\cir{}}
\POS(6,7) *{\scriptstyle{\bullet}}
\POS(9.3,.7) *{\scriptstyle{\bullet}}
\POS(-5.8,16) *{\scriptstyle{\bullet}}
\POS(15,12) *{\scriptstyle{\bullet}}

\POS(-1.5,18) \ar@{-} +(7.5,-11)
\POS(6,7) \ar@{-} +(3.3,-6.3)
\POS(6,7) \ar@{-} +(9,5)
\POS(-5.8,16) \ar@{-} +(4.3,2)


\POS(0,0)  \ar@{-} +(10,0)_4 
\POS(10,0) \ar@{-} +(7,7)_5
\POS(17,7) \ar@{-} +(0,10)_6
\POS(17,17) \ar@{-} +(-7,7)_7
\POS(10,24) \ar@{=} +(-10,0)_8
\POS(0,24) \ar@{-} +(-7,-7)_1 \POS(0,24)*+{\bullet}
\POS(-7,17) \ar@{-} +(0,-10)_2 
\POS(-7,7) \ar@{-} +(7,-7)_3

\POS(-3.5,20.5) \ar@{<<-} +(-3.2,-13.5)
\POS(6,24)  \ar@{<<-} +(-12.5,-17.)
\POS(14.5,4.5) \ar@{<<-} +(-4.5,20)
\POS(13.5,3.5)  \ar@{<<-} +(-13.5,-3.5)

\endxy

\vskip 10pt

Each $n$-fold dissection $D$ of $\pi$ gives rise to
$n$ \  $R$-deco polygons:
the root region inherits the root of $\pi$, and for any other region, there exists a unique 
arrow in its boundary which is the closest one to the root of $\pi$. Contracting 
arrows in the boundary of a region, we get an $R$-deco polygon. Its root side is the one where 
the arrow specified above ends.

\medskip
\begin{Example} \label{octagon} We illustrate the above on an $\N$-deco octagon $\pi$, where the sides are
decorated according to their place in the linear ordering of $\pi$.
The regions $r_1$, $r_2$ and $r_3$ below get equipped with the root sides 8, 5 and 8, respectively. 
To each of the first two regions we associate
an $R$-deco 3-gon, decorated by $\lceilx 3,5,8\rfloorx$ and $\lceilx 7,6,5\rfloorx$, respectively, while to region 
$r_3$ there is the 2-gon $\lceilx 1,8\rfloorx$ associated.
Finally, there are two more regions $r_4$ and $r_5$ leading to the 2-gons $\lceilx2,1\rfloorx$ and $\lceilx4,5\rfloorx$, respectively.
We illustrate the above by giving the polygons associated to $r_i$, $i=1,2,3$.

\end{Example}

\vskip 10pt
\hskip 50pt
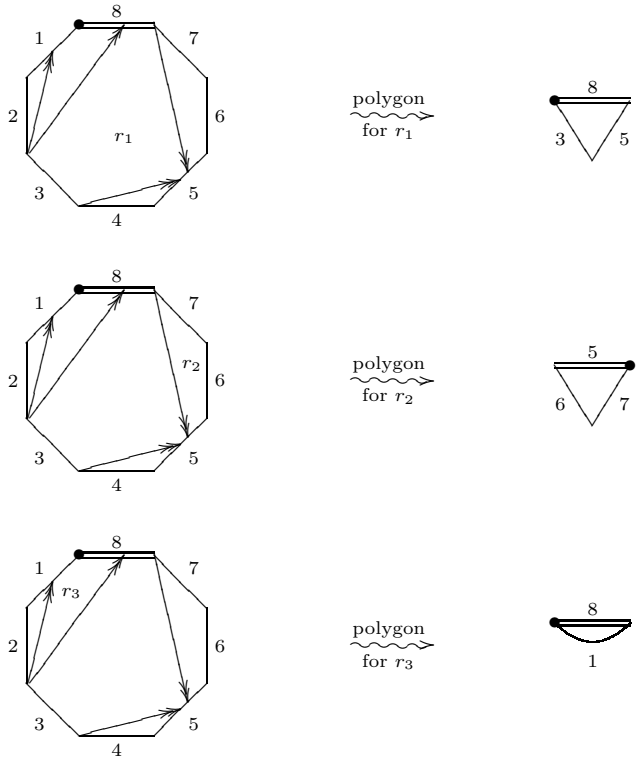
\xy (0,0) 
\POS(6,9) *{\scriptstyle r_1}
\POS(0,0)  \ar@{-} +(10,0)_4 
\POS(10,0) \ar@{-} +(7,7)_5
\POS(17,7) \ar@{-} +(0,10)_6
\POS(17,17) \ar@{-} +(-7,7)_7
\POS(10,24) \ar@{=} +(-10,0)_8
\POS(0,24) \ar@{-} +(-7,-7)_1 \POS(0,24)*+{\bullet}
\POS(-7,17) \ar@{-} +(0,-10)_2 
\POS(-7,7) \ar@{-} +(7,-7)_3 

\POS(-3.5,20.5) \ar@{<<-} +(-3.2,-13.5)
\POS(6,24)  \ar@{<<-} +(-12.5,-17.)
\POS(14.5,4.5) \ar@{<<-} +(-4.5,20)
\POS(13.5,3.5)  \ar@{<<-} +(-13.5,-3.5)

\hskip 100pt  \POS(0,12) \ar@{~>}^{\text{polygon}}_{\text{for }r_1} +(12,0) \hskip 80pt

\POS(10,14) \ar@{=} +(-10,0)_8 
\ar@{-} +(-5,-8)^5
\POS(0,14) \ar@{-} +(5,-8)_3
\POS(0,14)*+{\bullet}
\endxy

\vskip 10pt
\hskip 50pt
\xy (0,0) 

\POS(15,14) *{\scriptstyle r_2}
\POS(0,0)  \ar@{-} +(10,0)_4 
\POS(10,0) \ar@{-} +(7,7)_5
\POS(17,7) \ar@{-} +(0,10)_6
\POS(17,17) \ar@{-} +(-7,7)_7
\POS(10,24) \ar@{=} +(-10,0)_8
\POS(0,24) \ar@{-} +(-7,-7)_1 \POS(0,24)*+{\bullet}
\POS(-7,17) \ar@{-} +(0,-10)_2 
\POS(-7,7) \ar@{-} +(7,-7)_3 

\POS(-3.5,20.5) \ar@{<<-} +(-3.2,-13.5)
\POS(6,24)  \ar@{<<-} +(-12.5,-17.)
\POS(14.5,4.5) \ar@{<<-} +(-4.5,20)
\POS(13.5,3.5)  \ar@{<<-} +(-13.5,-3.5)

\hskip 100pt  \POS(0,12) \ar@{~>}^{\text{polygon}}_{\text{for }r_2} +(12,0) \hskip 80pt

\POS(10,14) \ar@{=} +(-10,0)_5 
\ar@{-} +(-5,-8)^7
\POS(0,14) \ar@{-} +(5,-8)_6
\POS(10,14)*+{\bullet}
\endxy

\vskip 10pt
\hskip 50pt
\xy (0,0)

\POS(-1,19) *{\scriptstyle r_3}

\POS(0,0)  \ar@{-} +(10,0)_4 
\POS(10,0) \ar@{-} +(7,7)_5
\POS(17,7) \ar@{-} +(0,10)_6
\POS(17,17) \ar@{-} +(-7,7)_7
\POS(10,24) \ar@{=} +(-10,0)_8
\POS(0,24) \ar@{-} +(-7,-7)_1 \POS(0,24)*+{\bullet}
\POS(-7,17) \ar@{-} +(0,-10)_2 
\POS(-7,7) \ar@{-} +(7,-7)_3 

\POS(-3.5,20.5) \ar@{<<-} +(-3.2,-13.5)
\POS(6,24)  \ar@{<<-} +(-12.5,-17.)
\POS(14.5,4.5) \ar@{<<-} +(-4.5,20)
\POS(13.5,3.5)  \ar@{<<-} +(-13.5,-3.5)

\hskip 100pt  \POS(0,12) \ar@{~>}^{\text{polygon}}_{\text{for }r_3} +(12,0) \hskip 80pt
\POS(10,15) \ar@{=} +(-10,0)_8 \POS(0,15)*+{\bullet}
\POS(0,15)*{}="A";
\POS(10,15)*{}="B";
"A";"B"; **\crv{(5,10)};
\POS(5,10) *+{\scriptstyle 1}
\endxy
\vskip 10pt












\subsection{Proof of Proposition \ref{just}}
We first prove the second claim. 
A typical term in $\overline \partial^2 \pi$ for an $R$-deco polygon $\pi$ arises from two non-intersecting arrows $\alpha$
and $\beta$.
There are two cases to consider, depending on whether the dual tree $\tau$ of the 3-fold dissection associated to
$(\alpha,\beta)$ is linear or not. 

In case $\tau$ is non-linear, it has the form
\vskip 5pt
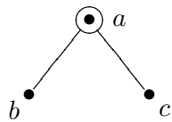
\xy \hskip 100pt \POS(0,0) *+{\bullet} *\cir{}  \ar@{-} +(-8,-10) *{\bullet} \ar@{-} +(8,-10) *{\bullet} \POS(4,0) *{a} \POS(-10,-12) *{b} \POS(10,-12) *{c} \endxy
\vskip 5pt
\noindent
for some subpolygons $a$, $b$ and $c$, where $a$ is the root polygon. In $\overline \partial^2\pi$ both terms $(a\wedge b)
\wedge c$ and $(a\wedge c)\wedge b$ occur with the same sign, regardless of the type of the arrows
$\alpha$ and $\beta$. So they  cancel in the wedge product. 

If $\tau$ is linear, it has the form 
\vskip 5pt
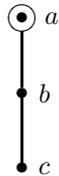
\xy \hskip 140pt \POS(0,0) *+{\bullet} *\cir{} \ar@{-} +(0,-10) *{\bullet} \POS(0,-10) \ar@{-} +(0,-10) *{\bullet}
\POS(4,0) *{a} \POS(3,-10) *{b} \POS(3,-20) *{c} \endxy
\vskip 5pt
\noindent
Let us first exclude the case when both arrows $\alpha$ and $\beta$ are backward, and end at the same side. Then,
by the Leibniz rule, cutting first the upper edge and then the lower one we obtain the same term
$a\wedge b\wedge c$ with the opposite sign compared to cutting the lower edge first. This 
does not depend on the directions of the arrows. In the excluded case we get $a\wedge c \wedge b + a\wedge b \wedge c =0$.
Finally, the augmentation is defined by taking the zero${}^\th$ component. 
This proves the second claim.

Let us prove the first claim. Given  an arrow  $\alpha$ in an $R$-deco polygon $\pi$, we use the notation  
$\varepsilon_{\alpha, \pi} :=  \varepsilon_\alpha$ to emphasize $\pi$. 
Let $\alpha$ and $\beta$ be non-intersecting arrows in an $R$-deco polygon $\pi$, and 
$\alpha$ is closer to the root than $\beta$. 
From now on, given a collection $D$ of non-intersecting arrows in $\pi$ including an arrow $\alpha$, 
we denote by $\pi_{\alpha}^\rst$ the unique polygon sitting right under the arrow $\alpha$, i.e., the one corresponding
 to the bottom vertex of the edge of the tree dual to $\alpha$. 
Denote by $a$ and $b$ the corresponding cut-off polygons, $a=\pi_{\alpha}^{\rst}$, 
$b=\pi_{\beta}^{\rst}$, and let $a \uparrow b$ be  the cut-off polygon for the single dissecting arrow $\alpha$ in $\pi$. 

\vskip 10pt
\hskip 80pt
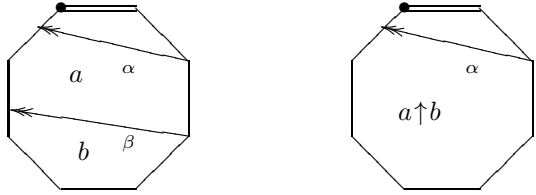
\xy (0,0) 
\POS(9,6) *{\scriptstyle \beta}
\POS(9,16) *{\scriptstyle \alpha}
\POS(3,5) *{ b}
\POS(2,15) *{ a}
\POS(0,0)  \ar@{-} +(10,0) 
\POS(10,0) \ar@{-} +(7,7)
\POS(17,7) \ar@{-} +(0,10)
\POS(17,17) \ar@{-} +(-7,7)
\POS(10,24) \ar@{=} +(-10,0)
\POS(0,24) \ar@{-} +(-7,-7) \POS(0,24)*+{\bullet}
\POS(-7,17) \ar@{-} +(0,-10) 
\POS(-7,7) \ar@{-} +(7,-7)

\POS(17,7) \ar@{->>} +(-23.6,3.5)
\POS(17,17) \ar@{->>} +(-20,4.5)

\hskip 130pt

\POS(9,16) *{\scriptstyle \alpha}
\POS(2,10) *{ a\hskip-2pt\uparrow\hskip-2pt b}
\POS(0,0)  \ar@{-} +(10,0) 
\POS(10,0) \ar@{-} +(7,7)
\POS(17,7) \ar@{-} +(0,10)
\POS(17,17) \ar@{-} +(-7,7)
\POS(10,24) \ar@{=} +(-10,0)
\POS(0,24) \ar@{-} +(-7,-7) \POS(0,24)*+{\bullet}
\POS(-7,17) \ar@{-} +(0,-10) 
\POS(-7,7) \ar@{-} +(7,-7)

\POS(17,17) \ar@{->>} +(-20,4.5)

\endxy

\vskip 20pt

In the above notation, one has the following obvious lemma. 

\begin{Lemma} \label{23we}
(1). $\chi(a)+ \chi(b) = \chi(a\hskip-2pt\uparrow\hskip-2pt b)$. 

(2). $
\varepsilon_{\alpha, \pi} \varepsilon_{\beta, \pi} = \varepsilon_{\beta, a\uparrow  b}
$. 
\end{Lemma}
Note that, if $\alpha$ and $\beta$ are backward in $\pi$, 
 then $\beta$ is not necessarily forward in $\overline \pi_\alpha^\rst$.
Thus the straightforward analog of (2) may not be correct for $\overline \pi_\alpha^\rst$ polygons. 

\vskip 3mm 
We compare the same two triple wedge products as above. 
The corresponding $R$-deco polygons match each other, so we have to worry only about the extra signs which we pick up using 
$\partial$ instead of $\overline \partial$. After these remarks, the proof is identical to the one above in the first case (i.e., for $\tau$ non-linear), and 
in the second case (i.e., if $\tau$ is linear) if the top arrow $\alpha$ is forward. Let us prove it in the second case, 
when  $\alpha$ is  backward. 

Assume that $\beta$ is forward. Then cutting out $\alpha$ first and $\beta$ second we pick up an extra sign 
$(-1)^{\chi(\alpha)+\chi(\beta)}$: indeed, by Lemma \ref{23we}(2) after the first cut $\beta$ becomes a backward arrow in $a\hskip-2pt\uparrow\hskip-2pt b$.  
By Lemma \ref{23we}(1) the weight of the middle polygon is $\chi(\alpha) - \chi(\beta)$, so cutting in the other order we pick up
 the same sign. 

Now assume $\beta$ is backward. Then cutting out $\alpha$ and then $\beta$ we pick up an extra sign 
$(-1)^{\chi(\alpha)}$: indeed, after the first cut $\beta$ becomes a forward
arrow in $a\hskip-2pt\uparrow\hskip-2pt b$.  Cutting in the other order we pick up 
the  sign $(-1)^{\chi(\beta)} (-1)^{\chi(\alpha)-\chi(\beta)} = (-1)^{\chi(\alpha)}$. 
 \qed

\subsection{The sign of a dissection}

The last ingredient needed for description of the cocycles in the bar construction of polygons is the sign of a dissection.
Given a dissection $D$ of an $R$-deco polygon $\pi$, and a dissecting arrow $\alpha$, there are two regions defined by the dissection $D$ containing $\alpha$. Take the one which is further away from the root on the dual tree, and 
denote by $\pi^\rst_{D, \alpha}$ the $R$-deco polygon obtained by shrinking its  arrow sides. 
Recall  the weight $\chi(\pi)$ of an $R$-deco polygon $\pi$.

\begin{Definition} Let $D$ be a dissection of an $R$-deco polygon $\pi$.
The {\bf sign}  of $D$ is 
$$
\sgn(D):= \prod_{\text{\rm backward arrows $\alpha$ of }D}  (-1)^{\chi(\pi^\rst_{D, \alpha})} = 
\prod_{\text{\rm arrows $\alpha$ of }D}\varepsilon_\alpha^{\chi(\pi^\rst_{D, \alpha})}\,.
$$
\end{Definition}

\medskip
\begin{Example} The dissection in Example \ref{octagon} has two backward arrows: the arrow ending in 1 cuts off a polygon
of weight~1 and thus contributes a
minus sign, while the arrow starting from between 7 and 8 contributes a plus sign: the cut-off polygon has weight~2. 
Thus the sign of the dissection is equal to $\,-1$.
\end{Example}

\subsection{The cocycle attached to a polygon in the bar construction} \label{dualtree}

We define the {\bf Adams grading} of an $R$-deco $N$-gon as  $2(N-1)$.
Since it is even, the shuffle product in the bar construction below is commutative.

Recall that the {\bf bar construction} $\B(\CA)$ associated to an augmented DGA $\CA_\smallbullet=\oplus_{m\geq 1} \CA_m$, 
with graded-symmetric multiplication $\wedge$ and differential $\partial$ of degree $+1$, 
is the tensor coalgebra $\oplus_i \CA_\smallbullet^{\otimes i}$ with differential $D_1+D_2$, where for homogeneous elements $a_1$, \dots, $a_N$ in $\CA_\smallbullet$ one defines, denoting any tensor sign by a $\Bar$, 

$$ D_1\big([a_1|a_2| \cdots| a_N]\big)=\displaystyle
\sum_{j=1}^{N-1}
(-1)^{\sum_{i\le j}({\rm deg}(a_i)-1)}[a_1| a_2|
\cdots| a_j \wedge a_{j+1}| \cdots
| a_N]\,,
$$
$$D_2\big([a_1\Bar\dots\Bar a_N]\big)=\displaystyle
 \sum_{j=1}^{N} (-1)^{\sum_{i<j}({\rm deg}({a_i})-1)}
 [a_1| a_2| \cdots|
\partial(a_j)| \cdots | a_N]\,.
$$


Furthermore, $\B(\CA_\smallbullet)$ carries an algebra structure, given by the graded-commutative shuffle product $\sha$.

Considering $\B(\CA_\smallbullet)$ as a double complex with respect to $D_1$ and $D_2$, we are particularly interested
in the cohomology of the ``main diagonal'', i.e., in $\H^0\B(\CA_\smallbullet)$, where the terms live which have only bars and 
no wedges.

The bar construction $\B(\CP_\smallbulletP)$ for the polygon algebra
$\CP_\smallbulletP=\CP_\smallbulletP(R)$ is a commutative Hopf algebra.
The unit in $\B(\CP_\smallbulletP)$
is given by $\Q\to \Q\cdot \fatone \hookrightarrow \B(\CP_\smallbulletP)$, the counit is the augmentation,
and the antipode is induced by the remaining structures. 
By the Milnor-Moore theorem, there exists a Lie algebra whose universal enveloping algebra is dual to
$\B(\CP_\smallbulletP)$. 


\begin{Definition}\label{bardef}
Let $\pi$ be an $R$-deco polygon. We associate to it an element $\B(\pi)$ in the bar construction $\B(\CP_\smallbulletP)$. 
Its component in $\CP_\smallbulletP^{| n}$ is the sum
$$\sum_D \sgn(D) \sum_{\lambda} \lceilx \pi_D^{\lambda(1)}\Bar \dots \Bar \pi_D^{\lambda(n)} \rfloorx$$
 over all $n$-fold dissections $D$ of $\pi$, where the inner sum runs through all linear orders $\lambda$
of the associated subpolygons $\pi_D^i$ compatible with the
partial order on $\tau(D)$.
\end{Definition}

\begin{Example} \label{barthree} The case $n=2$: the cocycle $\B(\lceilx 1\,2\,3\rfloorx)$ associated to a triangle 
$$\lceilx 1\,2\,3\rfloorx= \threegon{1}{2}{3} $$ 
is given as
\begin{equation*} 
 \Bigg[ \threegon{1}{2}{3}\quad {\large{,}} \quad \twogon{1}{3}\Barr  \twogon{2}{3}\ + \
\twogon{2}{3} \Barr  \twogon{1}{2}\ -\ \twogon{1}{3} \Barr  \twogonright{2}{1}\Bigg]\,.
\end{equation*}

Note that the terms of the second component of $\B(\lceilx 1\,2\,3\rfloorx)$  
cancel pairwise under $D_1$
with the terms of  the differential on polygons $\partial\big(\lceilx 1\,2\,3\rfloorx\big)$.
\end{Example}

\begin{Example} \label{barfour} The case $n=3$:  the cocycle associated to a quadrangle 
$$\lceilx 1\,2\,3\,4\rfloorx= \fourgon{1}{2}{3}{4}$$ 
has the following component in $\CP_\smallbulletP\Bar\CP_\smallbulletP$ (we omit $\lceilx $ and $\rfloorx$)
\begin{align*}
 & \quad\quad 14\Bar 234 \ +\ 34\Bar 123\ -\ 124\Bar 32\ -\ 134\Bar 21 \\
 & \,\ +\ 124\Bar 34 \ +\ 134\Bar 23\ +\ 234\Bar 12\ +\ 14\Bar321\,,
\end{align*}
corresponding to the following respective 2-fold dissections:
\vskip10pt
\xy 
\hskip 50pt \fourgona \POS(5,5) \ar@{<<-} +(-5,-10) \hskip 60pt \fourgona \POS(10,0) \ar@{<<-} +(-10,5)
\hskip 60pt \fourgona \POS(5,-5) \ar@{<<-} +(5,10) \hskip 60pt \fourgona \POS(0,0) \ar@{<<-} +(10,-5)
\endxy
\vskip 20pt
\xy 
\hskip 50pt  $14\Bar 234 \quad +\quad 34\Bar 123\quad -\quad 124\Bar 32\quad -\quad 134\Bar 21 $
\endxy
\vskip20pt
\xy 
\hskip 50pt \fourgona \POS(5,5) \ar@{<<-} +(5,-10) \hskip 60pt \fourgona \POS(10,0) \ar@{<<-} +(-10,-5)
\hskip 60pt \fourgona \POS(5,-5) \ar@{<<-} +(-5,10) \hskip 60pt \fourgona \POS(0,0) \ar@{<<-} +(10,5)
\endxy
\vskip 20pt
\xy 
\hskip 37pt  $\ +\ 124\Bar 34 \quad  +\quad 134\Bar 23\quad  +\quad 234\Bar 12\quad  +\quad  14\Bar321$
\endxy
\vskip 20pt
\noindent
The component of the cocycle in $\CP_\smallbulletP\Bar\CP_\smallbulletP\Bar\CP_\smallbulletP$ (here we use, e.g.,  
the shorthand $\,12\sha 34\,$ for $\,12\Bar 34\,+\,34\Bar 12\,$, the shuffle product of $\,12\Bar 34\,$ and $\,34\Bar 12\,$) is
\begin{align*}
& \ \ \ \ \ \, 14\Bar24\Bar34\ +\ 34\Bar13\Bar23\ -\ 24\Bar 12\sha 32 \ +\ 14\Bar 31\Bar 21\\
&\ +\ 14\Bar 34\Bar 23 \ +\ 34\Bar23\Bar 12 \ +\ 14\Bar21\Bar32 \ -\ 14\Bar 21\sha 34\\
& \ -\ 14\Bar 24\Bar 32 \ -\ 34\Bar13\Bar 21\ +\ 24\Bar 12\sha 34\ -\ 14\Bar31\Bar23\,.
\end{align*}
On pictures, the above twelve summands look as follows:
\vskip10pt
\xy 
\hskip 50pt \fourgona \POS(4,5) \ar@{<<-} +(-4,-10)  \POS(6,5) \ar@{<<-} +(4,-10) 
\hskip 60pt \fourgona \POS(10,1) \ar@{<<-} +(-10,4) \POS(10,-1) \ar@{<<-} +(-10,-4)
\hskip 60pt \fourgona \POS(6,-5) \ar@{<<-} +(4,10) \POS(4,-5) \ar@{<<-} +(-4,10)
\hskip 60pt \fourgona \POS(0,-1) \ar@{<<-} +(10,-4) \POS(0,1) \ar@{<<-} +(10,4)
\endxy
\vskip 20pt
\xy 
\hskip 50pt $14\Bar24\Bar34\ +\ 34\Bar13\Bar23\ -\ 24\Bar 12\sha 32 \ +\ 14\Bar 31\Bar 21$
\endxy
\vskip 20pt

\xy 
\hskip 50pt \fourgona \POS(5,5) \ar@{<<-} +(-5,-10) \POS(10,0) \ar@{<<-} +(-10,-5)
\hskip 60pt \fourgona \POS(10,0) \ar@{<<-} +(-10,5)\POS(5,-5) \ar@{<<-} +(-5,10)
\hskip 60pt \fourgona \POS(5,-5) \ar@{<<-} +(5,10) \POS(0,0) \ar@{<<-} +(10,5)
\hskip 60pt \fourgona \POS(0,0) \ar@{<<-} +(10,-5) \POS(5,5) \ar@{<<-} +(5,-10) 
\endxy
\vskip 20pt
\xy 
\hskip 36pt $\ +\ 14\Bar 34\Bar 23 \ +\ 34\Bar23\Bar 12 \ +\ 14\Bar21\Bar32 \ -\ 14\Bar 21\sha 34$
\endxy
\vskip 20pt
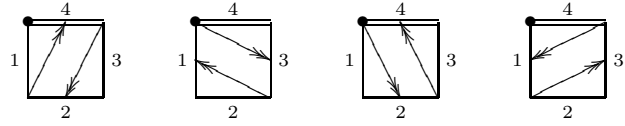
\xy 
\hskip 50pt \fourgona \POS(5,5) \ar@{<<-} +(-5,-10) \POS(5,-5) \ar@{<<-} +(5,10) 
\hskip 60pt \fourgona \POS(10,0) \ar@{<<-} +(-10,5) \POS(0,0) \ar@{<<-} +(10,-5)
\hskip 60pt \fourgona \POS(5,5) \ar@{<<-} +(5,-10) \POS(5,-5) \ar@{<<-} +(-5,10)
\hskip 60pt \fourgona \POS(10,0) \ar@{<<-} +(-10,-5)\POS(0,0) \ar@{<<-} +(10,5)
\endxy
\vskip 20pt
\xy 
\hskip 36pt $\ -\ 14\Bar 24\Bar 32 \ -\ 34\Bar13\Bar 21\ +\ 24\Bar 12\sha 34\ -\ 14\Bar31\Bar23$
\endxy
\vskip 20pt

Again, all the terms in the differential cancel
because they appear twice with different signs; we give two examples:
\begin{enumerate}
\item
$\partial(14\Bar 234)$ gives a term $-\,14\Bar24\wedge 34$ (the sign comes from the fact that we take the 
differential in the second factor) from taking the differential in $\CP_\smallbulletP$, while $14\Bar24
\Bar34$ gives a term $14\Bar24\wedge 34$ (the sign is $+$ since we replace the second $\Bar$ by a
$\wedge $).
\item 
$\partial(134\Bar 23)$ gives a term $+\,14\wedge  34\Bar23$ (no sign since we take the differential
in the first factor), while $14\Bar34\Bar23$ gives $- \,14\wedge  34\Bar23$ (the sign is $-$ 
since we replace the first $\Bar$ by a $\wedge $).
\end{enumerate}
\end{Example}
\vskip10pt

We need the multiplicativity of the sign of a dissection, formulated as follows. 
Let $D'\supset D$ be an overdissection of a dissection $D$. So for each arrow $\alpha$ of $D$ there is 
the cut-off polygon $\pi_{\alpha}^{\rst}$. We include here a ``root  arrow'' 
 providing the root polygon of the dissection as a cut-off polygon. Then there is a dissection $D_\alpha$ of the polygon 
$\pi_{\alpha}^{\rst}$ induced by $D'$. 
 
\begin{Lemma} \label{rty} One has 
\begin{equation} \label{rty1}
\sgn(D') = \sgn(D) \prod_{\alpha}\sgn(D_{\alpha}, \pi_{\alpha}^{\rst})
\end{equation}
where the product is over all arrows of $D$, including the root arrow. 

Further, let $\alpha$ be the unique predecessor of $\beta$ in the dual tree of $D$. Then 
\begin{equation} \label{rty2}
\sgn(D- \{\beta\}) = \sgn(D) (\eps_{\alpha}\eps_{\beta})^{\chi(\pi_\beta^\rst)},
\end{equation}
\end{Lemma}

\noindent{\bf Proof.} We first prove the second claim. Denote by $a$ and $b$ the subpolygons of
$D$ whose root arrows are $\alpha$ and $\beta$, respectively.
Then one has 
$$
\sgn(D) = \frac
{\sgn(D - \{\beta\})}{\varepsilon_{\alpha}^{\chi(a\uparrow b)}}\varepsilon_{\alpha}^{\chi(a)}\varepsilon_{\beta}^{\chi(b)}\,.
$$
So we have to check that 
$$
\frac
{\varepsilon_{\alpha}^{\chi(a)}\varepsilon_{\beta}^{\chi(b)}}
{\varepsilon_{\alpha}^{\chi(a\uparrow b)}} = (\varepsilon_{\alpha}\varepsilon_{\beta})^{\chi(b)}
$$
which follows from Lemma \ref{23we}(1). This proves the second claim.

We pass to the first claim. Given a dissection $D'$, we 
proceed by the induction on the number of arrows in $D$. If $D$ has just one arrow, it is the root arrow, and it is 
a forward arrow. So $\sgn (D) =1$, and 
(\ref{rty1}) is just the definition of $\sgn(D')$.

\hskip 140pt
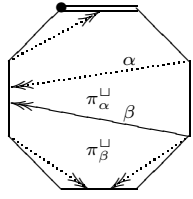
\xy (0,0)

\POS(9,10) *{\scriptstyle \beta}
\POS(9,17.) *{\scriptstyle \alpha}
\POS(0,0)  \ar@{-} +(10,0)
\POS(10,0) \ar@{-} +(7,7)
\POS(17,7) \ar@{-} +(0,10)
\POS(17,17) \ar@{-} +(-7,7)
\POS(10,24) \ar@{=} +(-10,0)
\POS(0,24) \ar@{-} +(-7,-7) \POS(0,24)*+{\bullet}
\POS(-7,17) \ar@{-} +(0,-10)
\POS(-7,7) \ar@{-} +(7,-7)
\POS(5,5) *{\scriptstyle \pi_\beta^\rst}
\POS(5,12) *{\scriptstyle \pi_\alpha^\rst}

\POS(17,17) \ar@{.>>} +(-23.6,-3.5)
\POS(17,7) \ar@{->>} +(-23.6,4.5)
\POS(-7,7) \ar@{.>>} +(10.5,-7)
\POS(17,7) \ar@{.>>} +(-10.5,-7)
\POS(-7,17) \ar@{.>>} +(12,6.5)
\endxy

Let us declare an arrow $\beta$ of $D'-D$ as a new arrow 
of the dissection $D$. Denote by $\alpha$ the dissecting arrow of $D$ just above $\beta$, i.e., preceding it in the dual tree. 
Then $\sgn(D')$ stays the same. By the second claim, $\sgn(D)$ is multiplied by 
$$
\varepsilon_\beta^{\chi(\pi_\beta^\rst)}\varepsilon_\alpha^{\chi(\pi_\beta^\rst)}\,.
$$

The $\sgn(D_{\alpha}, \pi_{\alpha}^{\rst})$ is multiplied by 
$\prod_{\gamma}(\varepsilon_\alpha\varepsilon_\gamma)^{\chi(\varphi_\gamma^\rst)}$, where $\gamma$ runs through all 
arrows in $\pi_{\alpha}^\rst$, and $\varphi_\gamma^\rst$ denotes the polygon cut-off by the arrow $\gamma$ in 
$\pi_{\beta}^{\rst}$. 
The $\sgn(D_{\beta}, \pi_{\beta}^{\rst})$ is multiplied by 
$\prod_{\gamma}(\varepsilon_\beta\varepsilon_\gamma)^{\chi(\varphi_\gamma^\rst)}$. Multiplying these two products, we see that 
the total change of 
$\prod_{\alpha}\sgn(D_{\alpha}, \pi_{\alpha}^{\rst})$ in (\ref{rty1}) is 
$$
\prod_{\gamma}(\varepsilon_\alpha\varepsilon_\beta)^{\chi(\varphi_\gamma^\rst)} = 
(\varepsilon_\alpha\varepsilon_\beta)^{\sum_{\gamma}\chi(\varphi_\gamma^\rst)}= 
(\varepsilon_\alpha\varepsilon_\beta)^{\chi(\pi_\beta^\rst)}\,.
$$
Here in the last step we used Lemma \ref{23we}(1). The lemma is proved. 
\qed

\begin{Proposition}\label{cocycle}
The element $\B(\pi)$ associated to the $R$-deco polygon $\pi$ is a 0-cocycle in $\B(\CP_\smallbulletP)$.
\end{Proposition}

\noindent{\bf Proof.} 
We check that the terms in the differential $(D_1+D_2)\B(\pi)$ cancel pairwise.

A typical term of $\B(\pi)$ is of the form $\sgn(D) \lceilx a_1\Bar\dots\Bar a_n\rfloorx$ ($n\geq 1$) where the order 
$(a_1,\dots,a_n)$ on the regions $a_n$ of an $n$-fold 
dissection $D$ of $\pi$ is compatible with the partial order $\prec\,$ on the dissection.

Therefore, a typical term of $D_1\big(\B(\pi)\big)$ is of the form ($\deg(a_i)$ is even)
\begin{equation}\label{typical}
(-1)^j\,\sgn(D)\,\lceilx a_1\Bar\dots\Bar a_j\wedge a_{j+1}\Bar\dots\Bar a_n\rfloorx
\end{equation}
with the same requirements on the $a_i$ (i.e., the compatibility of orders).
Since the differential $\partial$ respects the partial order of a dissection, the terms of $D_2\big(\B(\pi)\big)$  
are also ordered compatibly with $\prec$.

\smallskip
We distinguish two cases:

\smallskip
{\bf Case 1:} The regions $a_j$ and $a_{j+1}$ are comparable with respect to $\prec$. 
Then the same term arises, up to sign, from $D_2\big(\B(\pi)\big)$: we must have $a_j\prec a_{j+1}$, 
and they have to be neighbouring regions (any $\gamma$
with $a_j\prec \gamma\prec a_{j+1}$ would have to occur between $a_j$ and $a_{j+1}$ in any compatible linear order which would exclude a term 
with component $a_j\wedge a_{j+1}$). 

Let $\beta$ be the arrow shared by $a_j$ and $a_{j+1}$, and let $\alpha$ be its predecessor
in the dual tree, which is also the root arrow of $a_j$. Set $D'=D-\{\alpha_{j+1}\}$.
The term $\ \sgn(D') \lceilx a_1\Bar\dots\Bar a_n\rfloorx\ $ of $\B(\pi)$ gives, via $D_2$, a term
\begin{equation} \label{splitsign}
(-1)^{j-1}\sgn(D') (\varepsilon_{\alpha}\varepsilon_{\beta})^{\chi(a_{j+1})}\lceilx a_1\Bar\dots \Bar a_j\wedge a_{j+1}\Bar \dots \Bar a_n\rfloorx\,.
\end{equation}
Equation (\ref{rty2}) of the above lemma shows that
\begin{equation} \label{splitsign1}
\sgn(D) = \sgn(D')(\varepsilon_{\alpha}\varepsilon_{\beta})^{\chi(a_{j+1})}\,.
\end{equation}
Therefore the terms (\ref{splitsign}) and (\ref{typical})  cancel.

Note that all terms in $D_2\big(\B(\pi)\big)$ are thus taken care of.

\smallskip
{\bf Case 2:} Suppose the regions $a_j$ and $a_{j+1}$ are {\em not} comparable with respect to $\prec$. Then 
there is a 
summand $\ \sgn(D)\lceilx a_1\Bar\dots\Bar a_{j+1}\Bar a_j\Bar \dots\Bar a_n\rfloorx\ $ whose $D_1$-image
gives a term (\ref{typical}) with opposite sign, due to antisymmetry. This proves Case~2. 

\medskip
In both cases the terms in $(D_1+D_2)\big(\B(\pi)\big)$ cancel pairwise, so $\B(\pi)$ is a cocycle.\qed

\section{Recognizing the coproduct for $\B(\pi,\partial)$} \label{recognize}
In this section, we show that the algebra generated by the $\B(\pi)$ for $R$-deco polygons $\pi$ forms
a sub-Hopf algebra of $\B(\CP_\smallbullet)$. 
Let us formulate the main result. 
\vskip 2mm
Recall (\cite{CK1}) that an {\bf admissible cut} in a rooted tree $\tau$ is given by a collection $C$ of edges of $\tau$ such that
each simple path from the root to a leaf vertex contains at most one element of $C$. 
We use the same notion for the rooted {\em plane} trees. We call a dissection $D$
{\bf admissible} if the set of those edges of $\tau(D)$ corresponding to the arrows in $D$ form an admissible cut of $\tau(D)$.

An admissible dissection $D$ gives rise to  a region $R^D$ containing the root of
$\pi$, and the remaining regions $P_i^{D}$.
There are two particular cases: the empty cut of $\tau(D)$ corresponding to
$R^D=\emptyset$, and the "full cut" corresponding to $R^D=\pi$.
\begin{Theorem} \label{MT} The coproduct $\Delta\big(\B(\pi)\big)$ for a polygon $\pi$ has the following expression
\begin{equation}\label{CKlike}
\Delta\big(\B(\pi)\big) =\sum_{D\rm{\ adm}} \sgn(D) \ 
\B(R^D) \otimes \Sha{\ i\  }{\phantom{i=1}} \B\big(P_i^{D}\big)\,,
\end{equation}
where $D$ runs through all admissible dissections in $\pi$. 

Thus the subalgebra generated by the $\B(\pi)$ is a Hopf subalgebra of $\B(\CP_\smallbullet)$. 
\end{Theorem}

Here is a rough idea of the proof. 
The terms of $\B(\pi)$ for an $R$-deco polygon 
$\pi$ are parametrized by dissections, together with a
linear ordering on the ensuing smaller polygons.
The terms in the coproduct arise from splitting 
each such linearly ordered set into two connected parts, the left and the right. 
The left one---if non-empty---contains the root region of the underlying dissection.

We regroup those terms in the form $\sum_i \pm \B(\pi_i^1)\otimes \Shaonly_{j>1} \B(\pi_i^j)$ for certain polygons 
$\pi_i^j$, where on the right  we have a (commutative) shuffle product.

\subsection{Level structure on plane trees} \label{leveltrees}
Recall that a dissection of an $R$-deco polygon $\pi$ provides a rooted {\em plane} tree whose vertices 
are labeled by the regions. 

Let $\B(\pi)$ be the element in the bar construction of $\B(\CP_\smallbullet)$, let $D$ be a dissection of $\pi$ and $\tau=\tau(D)$
the associated rooted plane tree.
For the terms in $\B(\pi)$ we need to control the linear orders compatible with this tree. It is convenient to replace 
the notion of ``linear orders compatible with the partial order on $\tau$'' by the equivalent one of a ``tree level structure''
(this notion was used by Loday in \cite{Loday}): 
each compatible linear order corresponds to a (metrically distorted) version of the tree where each vertex lies on a different
level. A linear order on the edges of $\tau$ arises from a level structure by identifying the vertex at level $i$ with the unique 
edge ending at that vertex. For example, we give a simple tree $\tau$ together with two of its eight different level structures $\tau_i$:
\vskip7pt
\xy 
\POS(0,0) \ar@{-} +(-3,-5) \ar@{-} +(6,-10) 
\POS(-3,-5) \ar@{-} +(-3,-5) \ar@{-} +(3,-5) 
\POS(0,0) *+{\bullet} \POS(-3,-5) *+{\bullet} \POS(0,-10) *+{\bullet} \POS(-6,-10) *+{\bullet} \POS(6,-10) *+{\bullet} 
\POS(-5,0) *+{\tau}

\hskip 80pt \POS(-10,-5) \ar@{~>} +(6,0) \hskip40pt
\POS(0,0) \ar@{-} +(-3,-5) \ar@{-} +(12,-20) \POS(-3.5,-2) *{_{e_1}}
\POS(-3,-5) \ar@{-} +(-3,-5) \ar@{-} +(6,-10)  \POS(-7,-7) *{_{e_2}} \POS(-.5,-13) *{_{e_3}} \POS(8.5,-18) *{_{e_4}}
\POS(0,0) *+{\bullet} \POS(-3,-5) *+{\bullet} \POS(3,-15) *+{\bullet} \POS(-6,-10) *+{\bullet} \POS(12,-20) *+{\bullet} 
\POS(-8,2) *+{\tau_1}
\POS(-3,0) \ar@{.} +(6,0) \POS(10,0)  *+{\scriptstyle{{\rm level\ }0}}
\POS(-8,-5) \ar@{.} +(16,0) \POS(12,-5)  *+{\scriptstyle{1}}
\POS(-10,-10) \ar@{.} +(20,0) \POS(14,-10)  *+{\scriptstyle{2}}
\POS(-10,-15) \ar@{.} +(22,0) \POS(16,-15)  *+{\scriptstyle{3}}
\POS(-10,-20) \ar@{.} +(25,0) \POS(18,-20)  *+{\scriptstyle{4}}

 \hskip120pt
\POS(0,0) \ar@{-} +(6,-10) \ar@{-} +(-12,-20) \POS(-3.5,-2) *{_{e_1}}
\POS(-3,-5) \ar@{-} +(-3,-5) \ar@{-} +(6,-10) \POS(-8,-18) *{_{e_4}} \POS(-.5,-13) *{_{e_3}} \POS(6.5,-7) *{_{e_2}}
\POS(0,0) *+{\bullet} \POS(-3,-5) *+{\bullet} \POS(3,-15) *+{\bullet} \POS(6,-10) *+{\bullet} \POS(-12,-20) *+{\bullet} 
\POS(-8,2) *+{\tau_2}
\POS(-3,0) \ar@{.} +(6,0) \POS(11,0)  *+{\scriptstyle{{\rm level\ }0}}
\POS(-8,-5) \ar@{.} +(16,0) \POS(14,-5)  *+{\scriptstyle{1}}
\POS(-10,-10) \ar@{.} +(20,0) \POS(14,-10)  *+{\scriptstyle{2}}
\POS(-12.5,-15) \ar@{.} +(22,0) \POS(14,-15)  *+{\scriptstyle{3}}
\POS(-15,-20) \ar@{.} +(25,0) \POS(14,-20)  *+{\scriptstyle{4}}
\endxy
\vskip7pt
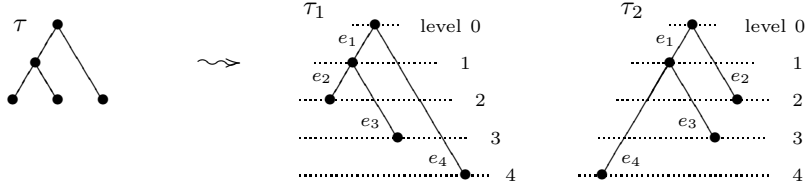
\centerline{\small A tree $\tau$, and two of its eight level trees, denoted $\tau_1$ and $\tau_2$.}

\subsection{Regrouping the terms in the coproduct $\Delta\big(\B(\pi)\big)$}\label{regrouping}
Let $\pi$ be an $R$-deco polygon. Its associated element $\B(\pi)$ is given in terms of a double sum $\displaystyle{ \sum_{D\text{ diss.}}
\sum_{\lambda\text{ lin.order}}}$, where the first sum runs through all dissections $D$ of $\pi$ while the second sum runs 
through all linear orders $\lambda$ compatible with the induced partial order on (the  regions $\pi_i$ of) $D$. 
The coproduct of $\B(\pi)$ is given by a triple sum
\begin{equation}\label{triplesum}
\ \sum_{D\text{ diss.}}\ \sum_{\lambda\text{ lin.order}}\ \sum_{\text{splitting of }\lambda}\,, 
\end{equation}
where the first two sums 
are as above and the third sum runs through all the $|D|+1$ splittings of the sequence of polygons $\pi=(\pi_1,\dots,\pi_{|D|})$ 
into a ``left'' part $(\pi_1,\dots,\pi_n)$ and  ``right'' part  $(\pi_{n+1},\dots,\pi_{|D|})$.
The empty parts are treated as the unit $\fatone\in \B(\CP_\smallbullet)$.

Now replace each compatible linear order $\lambda$ by a level structure of $\tau(D)$. 
The levels 
are ordered from top to bottom (as in the pictures above). Then a splitting of $\lambda$ corresponds
to a horizontal cut (up to isotopy) of $\tau(D)$ avoiding the vertices.
\vskip10pt
\xy 
\hskip 120pt 
\POS(0,0) \ar@{-} +(-3,-5) \ar@{-} +(12,-20) \POS(-3.5,-2) *{_{e_1}}
\POS(-3,-5) \ar@{-} +(-3,-5) \ar@{-} +(6,-10)  \POS(-7,-7) *{_{e_2}} \POS(-.5,-13) *{_{e_3}} \POS(8.5,-18) *{_{e_4}}
\POS(0,0) *+{\bullet} \POS(-3,-5) *+{\bullet} \POS(3,-15) *+{\bullet} \POS(-6,-10) *+{\bullet} \POS(12,-20) *+{\bullet} 
\POS(-8,2) *+{\tau_1}
\POS(-3,0) \ar@{.} +(6,0) \POS(10,0)  *+{\scriptstyle{{\rm level\ }0}}
\POS(-8,-5) \ar@{.} +(16,0) \POS(12,-5)  *+{\scriptstyle{1}}
\POS(-10,-10) \ar@{.} +(20,0) \POS(14,-10)  *+{\scriptstyle{2}}
\POS(-10,-15) \ar@{.} +(22,0) \POS(16,-15)  *+{\scriptstyle{3}}
\POS(-10,-20) \ar@{.} +(25,0) \POS(18,-20)  *+{\scriptstyle{4}}
\POS(-13,-12.5) \ar@{--} +(33,0) \POS(30,-12.5)  *+{\scriptstyle{\text{horizontal cut}}}
\endxy
\vskip5pt
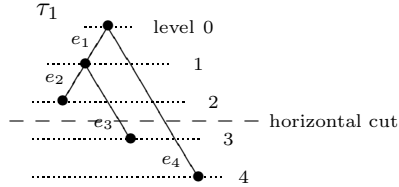
\centerline{\small A horizontal cut of the level tree $\tau_1$ with edge sequence $(e_3,e_4)$.}
\vskip10pt
\medskip
The ``left'' part of the splitting of $\lambda$ corresponds to the upper (=root) part, while the ``right'' part
corresponds to the lower part of the horizontal cut.
Therefore the above triple sum (\ref{triplesum}) is rewritten as
\begin{equation}\label{triplesum1}
  \sum_{D\rm{\ diss.}} \quad \sum_{\rm{\ level\ tree\ }\tau'\rm{\ for\ \tau(}D\rm{)}} 
\quad\sum_{{\rm horizontal\ cut\ of\ }\tau'}\,. 
\end{equation}

Each horizontal cut is described by the sequence
$E(\lambda)$ of edges in $\tau(D)$ which it meets, 
together with the induced level structure on the remaining parts of 
$\tau(D)- E(\lambda)$. The set $E(\lambda)$ gives rise to an admissible dissection of $D$. 

Horizontal cuts and admissible dissections are related as follows: 

(i) Each horizontal cut of a level tree for $\tau$
obviously produces an admissible cut for $\tau$ by forgetting the level structure. 

(ii) Each admissible dissection $\overline D$ of $\tau$ can be connected by a path $\beta$ from ``left of tree'' to ``right 
of tree'' intersecting the edges in $\tau$ precisely once and avoiding the remaining parts of the tree. 
The pair $(\tau,\beta)$ can be ``straightened out'', 
not necessarily uniquely, via a plane isotopy such that $\beta$ becomes a 
horizontal cut and $\tau$ obtains a level structure.

\vskip10pt
\hskip50pt
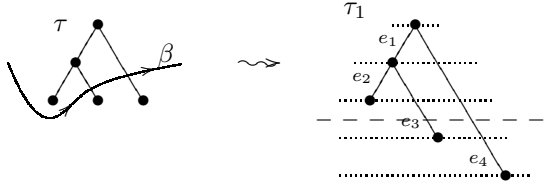
\xy   
(-12,-5)*+{}; (12,-5)*+{}
 **\crv{(-8,-15)&(-3,-10)&(0,-7)} ?(.5)*\dir{>} ?(.93)*\dir{>} 
\POS(9,-4) *+{\beta}
\POS(0,0) \ar@{-} +(-3,-5) \ar@{-} +(6,-10) 
\POS(-3,-5) \ar@{-} +(-3,-5) \ar@{-} +(3,-5) 
\POS(0,0) *+{\bullet} \POS(-3,-5) *+{\bullet} \POS(0,-10) *+{\bullet} \POS(-6,-10) *+{\bullet} \POS(6,-10) *+{\bullet} 
\POS(-5,0) *+{\tau}

\hskip 80pt \POS(-10,-5) \ar@{~>} +(6,0) \hskip40pt
\POS(0,0) \ar@{-} +(-3,-5) \ar@{-} +(12,-20) \POS(-3.5,-2) *{_{e_1}}
\POS(-3,-5) \ar@{-} +(-3,-5) \ar@{-} +(6,-10)  \POS(-7,-7) *{_{e_2}} \POS(-.5,-13) *{_{e_3}} \POS(8.5,-18) *{_{e_4}}
\POS(0,0) *+{\bullet} \POS(-3,-5) *+{\bullet} \POS(3,-15) *+{\bullet} \POS(-6,-10) *+{\bullet} \POS(12,-20) *+{\bullet} 
\POS(-8,2) *+{\tau_1}
\POS(-3,0) \ar@{.} +(6,0) \POS(10,0)  
\POS(-8,-5) \ar@{.} +(16,0) \POS(12,-5)  
\POS(-10,-10) \ar@{.} +(20,0) \POS(14,-10) 
\POS(-10,-15) \ar@{.} +(22,0) \POS(16,-15) 
\POS(-10,-20) \ar@{.} +(25,0) \POS(18,-20) 
\POS(-13,-12.5) \ar@{--} +(30,0) \POS(20,-12.5) 
\endxy
\vskip10pt
\centerline{\small An admissible cut of $\tau$ and an associated horizontal cut of a level tree $\tau_1$ for $\tau$.}
\medskip

\def \Dbar{{\overline{D}}}

Using these remarks, we can rewrite (\ref{triplesum1}) as follows: 
$$  \sum_{D\rm{\ diss.}} \quad \sum_{\rm{\ level\ tree\ }\tau'{\rm\ for\ \tau(}D\rm{)}} 
\quad\sum_{{\rm adm\ diss\ }\Dbar \subset D}\,. $$

\medskip
We can reverse the order of the sum in this triple sum, first summing over the {\em admissible} 
dissections $\Dbar$, then summing over all dissections $D$ of $\pi$ 
{\em containing} $\Dbar$, and finally summing over all level trees for $\tau(D)$, i.e.
$$  \sum_{\Dbar\rm{\ adm.\ diss.}} \quad\sum_{{\rm diss\ } D \supset \Dbar}
 \quad \sum_{{\rm \ level\ tree\  for\ \tau(}D{\rm)}} 
\,. $$

Recall that an admissible dissection $D$ separates a 
``root piece'' $R^{D}$ from the remaining  regions $P_{i}^D$. 
Each overdissection $D'\supset D$ gives 
both a dissection of the root piece $R_{D'- D}^D$ and a dissection of the remaining  regions $ P^{D}_{i,D'-  D}$. 
The separated parts can be resummed independently: more precisely, a plane tree $\tau$ with horizontal cut $C$ cuts $\tau$
into a root piece $\tau_R$ and a forest $\phi=\disjoint_i \tau_{P_i}$ of remaining pieces, which provides a bijection between
\begin{enumerate}
\item the level structures of the pair $(\tau, C)$ and
\item the product \{level structures on $\tau_R\} \times \{$level structures on $\phi\}$ .
\end{enumerate}

\medskip
Decomposing $D'- D= D_1 \disjoint  D_2$ into the 
dissections $D_1$ and $D_2$ of the root piece and of the remaining pieces, respectively, and using Lemma \ref{rty}, which implies the full multiplicativity for an {\em admissible} dissection,
we get
\begin{eqnarray*}
 \Delta\big(\B(\pi)\big)&=& \sgn(D')\sum_{D\rm{\ adm}} \quad \sum_{D'\supset D} \quad\sum_{{\rm level\ str.\ on\ }\tau(D')}
 R_{D'- D}^D \otimes \Sha{ i\ }{\phantom{i=1}} P^{D}_{i,D'-  D} \\
&=& \sum_{D{\rm \ adm}} \sgn(D) 
\Bigg(\sum_{D_1}\sgn(D_1)\sum_{{\rm level\ str.\ for\ }\tau(D_1)} R_{D_1}^D  \Bigg) \otimes \\
&& \qquad \otimes \Bigg(\sum_{D_2} \sum_{{\rm level\ str.\ 
for\ }\tau(D_2)}\Sha{\ i\  }{\phantom{i=1}} \sgn(D_{2,i}) P^{D}_{i,D_2}\Bigg)  \\
&=& \sum_{D{\rm \ adm}} \sgn(D)\ \B(R^D) \otimes \Sha{i  }{\phantom{i=1}} \B\big(P_i^{D}\big)\,.
\end{eqnarray*}
We need the  shuffle product $\Shaonly$ since precisely those linear orders on the remaining
 regions which are induced from such a shuffle do occur.
Theorem \ref{MT} is proved. \qed

\vskip 3mm
The vector space $V^\pg$ maps to $\H^0\B(\CP)$ via $\pi\to \B(\pi)$ 
for a polygon $\pi$. Theorem \ref{MT} implies that this map extends to a map of the symmetric algebra $\SS^\bullet V^\pg$ of $V^\pg$ 
 to $\H^0\B(\CP)$ sending 
the product in $\SS^\bullet V^\pg$ to the shuffle product in $\H^0\B(\CP)$.

We can view $\SS^\bullet V^\pg$ as a Hopf algebra via pullback of the 
coproduct formula (\ref{CKlike}).
So there is a Hopf algebra map $\SS^\bullet V^\pg \to \chi_\cycle=\H^0\B(\CN)$, defined via $\H^0\B(\CP)$. 

\section{Comparing the coproduct for $\B(\pi)$ and the coproduct for the Hopf algebra of iterated integrals
from \cite{GonFund}}

\subsection{The Hopf algebra $\CI_\smallbullet(R)$ of iterated integrals.}  Let $R$ be a set. 
In \cite{GonFund}, the second author described a commutative, graded Hopf 
algebra $\CI_\smallbullet(R)$. 
It is generated, as a graded commutative algebra, by elements $\I(a_0;a_{1},\dots, a_{n};a_{n+1})$ of degree $n$, $a_i\in R$. 
Their properties resemble the ones of  the iterated integrals. 

The coproduct is given via polygons inscribed in a semicircle. 
For the above generator it is written as a sum indexed by all---possibly empty---subsequences of $(a_1,\dots,a_n)$.
A picture describing the combinatorics of the terms arises when we place the $a_i$ on a semicircle, in the order 
dictated by their index. In particular, $a_0$ and $a_{n+1}$ are located at the endpoints of the semicircle.

Precisely, each subsequence $(a_{i_1},\dots,a_{i_r})$ of $(a_1,\dots,a_n)$ gives rise to a single ``main'' polygon with vertex set
$\{a_0,a_{i_1},\dots,a_{i_r},a_{n+1}\}$ and to a set of $r+1$ remaining polygons with $\geq2$ sides and set of vertices
$\{a_{i_0},a_{i_0+1},\dots,a_{i_1}\}$, $\{a_{i_1},a_{i_1+1},\dots,a_{i_2}\}$, \dots, $\{a_{i_r},a_{i_r+1},\dots,a_{i_{r+1}}\}$, 
where $a_{i_0}=a_0$ and $a_{i_{r+1}}=a_{n+1}$. We give an example for $n=4, r=2$.
\vskip10pt
\xy
\hskip 90pt
\POS(36,0) *{a_{5}}
\POS(-34,0) *{a_{0}}
\POS(-20,22.36) *{\bullet}
\POS(0,30) *{\bullet}
\POS(-26,14.97) *{\bullet}
\POS(23,19.26) *{\bullet}
\POS(-23.5,23.5) *{a_2}
\POS(-29.5,15) *{a_1}
\POS(-2,33) *{a_3}
\POS(27,20) *{a_4}

\POS(-30,0) *{\bullet}; \POS(30,0)*{\bullet},
{\ellipse_{}} 
\POS(-30,0) \ar@{-} +{(10,22.36)} \ar@{.} +{(60,0)}
\POS(-20,22.36) \ar@{-} +{(43,-3.1)}
\POS(23,19.26) \ar@{-} +{(6.81,-19.26)}
\endxy 
\vskip5pt
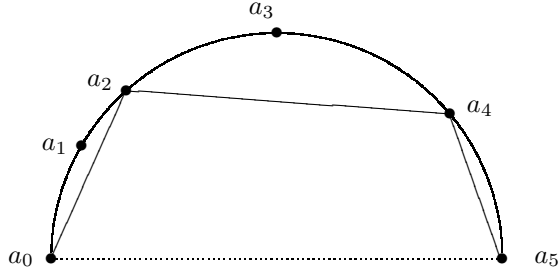
\centerline{\small A 4-gon for the subsequence $(a_2,a_4)$ of $(a_1,\dots,a_4)$ inscribed in the semicircle}
\vskip10pt
The corresponding term in the coproduct is then 
\begin{equation} \label{semicircle}
\I(a_0;a_{i_1},\dots,a_{i_r};a_{n+1}) \ \otimes \ \prod_{k=0}^r \I(a_{i_k};a_{i_k+1},\dots,a_{i_{k+1}-1};a_{i_{k+1}})\,,
\end{equation}

These generators are subject to the following relations:

\noindent
the {\em path composition formula} with respect to an element
$x\in R$:
\begin{equation} \label{sumx}
\I(a_0;a_1,\dots,a_n;a_{n+1}) = \sum_{k=0}^{n} \I(a_0;a_1,\dots,a_k;x) \,\I(x;a_{k+1},\dots,a_n;a_{n+1}) \,,
\end{equation}
the {\em inversion formula}
\begin{equation} \label{inversion}
\I(a_{n+1};a_n,\dots,a_1;a_{0}) =(-1)^n \I(a_0;a_1,\dots,a_n;a_{n+1})\,.
\end{equation}
and the {\it unit identity} 
\begin{equation} \label{unitx}
\I(a;b) = 1 \quad \mbox{for any $a,b \in R$}\,.
\end{equation}

Let us replace the $k$-th factor of the right hand product in (\ref{semicircle}) by the sum 
 arising from the path 
composition formula for $x=0$:
\begin{equation} \label{semicircle2}
\sum_{j=i_k}^{i_{k+1}} \I(a_{i_k};a_{i_k+1},\dots,a_j;0) \,
\I(0;a_{j+1},\dots,a_{i_{k+1}-1};a_{i_{k+1}})\,.
\end{equation}

\medskip
Using this and the inversion relation (\ref{inversion}), we write (\ref{semicircle}) as 
$$
\I(a_0;a_{i_1},\dots,a_{i_r};a_{n+1}) \ \otimes 
$$
\begin{equation} \label{semicircle3}
\ \prod_{k=0}^r \sum_{j=i_k}^{i_{k+1}} (-1)^{j-i_k}
\I(0;a_j,\dots,a_{i_k+1};a_{i_k}) \,
\I(0;a_{j+1},\dots,a_{i_{k+1}-1};a_{i_{k+1}})\,.
\end{equation}


We denote by $\CI_\smallbullet(R)$ the so obtained Hopf algebra. 
In \cite{GonFund} several similar Hopf algebras were defined: 
they differ by the relations added to the basic relations (\ref{sumx})-(\ref{unitx}). 
The Hopf algebra $\CI_\smallbullet(R)$ is the biggest one of them. 

\subsection{The Lie coalgebra of iterated integrals and polygons.} Recall that given a commutative Hopf algebra ${\cal I}_{\bullet}$, graded by 
the integers $n=0, 1, ...$, with $ {\cal I}_{0}=\Q$, the graded $\Q$-vector space of indecomposables
$$
Q({\cal I}_{\bullet}):= \frac{{\cal I}_{>0}}{{\cal I}_{>0}^2}
$$
has a natural graded Lie coalgebra structure with the cobracket 
$\delta: Q({\cal I}_{\bullet}) \longrightarrow \Lambda^2Q({\cal I}_{\bullet})$ induced by the coproduct in the Hopf algebra. 
Applying this to the Hopf algebra $\CI_\smallbullet(R)$, we arrive at the graded Lie coalgebra 
$Q(\CI_\smallbullet(R))$. We denote by $\overline \I(a_0;a_{i_1},\dots,a_{i_r};a_{n+1}) $ the projection of the generator 
$\I(a_0;a_{i_1},\dots,a_{i_r};a_{n+1}) $ to it. 

Recall the  $R$-deco $(n+1)$-gon  $[a_1, ..., a_{n+1}]$ decorated by the set 
$(a_1, ..., a_{n+1})$ so that the root side is decorated by $a_{n+1}$. 

\begin{Proposition} \label{hihiks} Let $R$ be a set. The map 
$$
\overline \I(0;a_{1},\dots,a_{n};a_{n+1}) \to [a_1, ..., a_{n+1}]
$$
gives rise to an isomorphism of the graded Lie coalgebras $$
Q(\CI_\smallbullet(R)) \stackrel{\sim}{\longrightarrow} 
(V_{\smallbullet}^\pg(R), \partial)\,.$$
\end{Proposition}

\noindent{\bf Proof.} It follows from (\ref{sumx}) - (\ref{inversion}) that the elements 
$\overline \I(0;a_{1},\dots,a_{n};a_{n+1})$ form a basis of the vector space $Q(\CI_\smallbullet(R))$.

\begin{lemma} \label{forx}
The Lie cobracket is given by 
\begin{equation} \label{semi3}
\delta\overline \I(0;a_{1},\dots,a_{n};a_{n+1}) = 
\sum_{0 \leq k < l \leq n+1}\overline \I(0;a_{1},\dots,a_{k}, a_l, \ldots, a_n;a_{n+1})\wedge 
\end{equation}
\begin{equation} \label{semi4}
\Bigl( 
\overline \I(0;a_{k+1},\dots, a_{l-1};a_{l}) + (-1)^{l-k-1} \overline \I(0;a_{l-1},\dots, a_{k+1};a_{k}) \Bigr)\,.
\end{equation}
\end{lemma}

\noindent{\bf Proof.} To get a non-zero term, the sequence of elements 
$\{0,a_{i_1},\dots,a_{i_r},a_{n+1}\}$ determining a term in the coproduct must have exactly one gap, i.e. all sides 
of the corresponding inscribed polygon except just one must be of length $1$. Thus only sequences 
$\{0, a_1, ..., a_k, a_l, \ldots, a_n, a_{n+1}\}$ can appear where $0 \leq k < l \leq n+1$. Given such a sequence, 
the corresponding term of the coproduct is
$$
\overline \I(0;a_{1},\dots,a_{k}, a_l, \ldots, a_n;a_{n+1})\wedge 
\overline \I(a_k;a_{k+1},\dots, a_{l-1};a_{l}) \,.
$$
Using (\ref{sumx}) and (\ref{inversion}), we get the required formula. \qed

The formula (\ref{semi3}) can be interpreted via $R$-deco polygons as follows. 
The two terms in (\ref{semi4}) correspond to two dissections of the polygon $[a_1, ..., a_{n+1}]$, given by the forward arrow 
starting from the vertex sharing $a_{k}$ and $a_{k+1}$ and ending at the side $a_l$, and, if $k>0$,  the backward arrow 
starting at the vertex sharing $a_{l-1}$ and $a_l$ and ending at $a_k$. The sign of this dissection is 
$(-1)^{l-k-1}$ for the backward arrow. \qed. 

\subsection{Comparing the two coproducts.} 
\begin{Theorem}
There is a map of coalgebras from $\langle\B(\pi)\mid \pi\in \CP\rangle$ to ${\cal I}_\smallbullet(R)$.
\end{Theorem}

\noindent{\bf Proof.}  The map is given by sending a 
generator $\B(\pi)$, assigned to an $R$-deco polygon $\pi$ 
decorated by  $(a_1,\dots,a_{n+1})$, the generator $\I(0;a_1,\dots,a_{n};a_{n+1})$.

To show that the coproduct on $\B(\pi)$ is compatible with $\Delta \I(0;a_1,\dots,a_{n};a_{n+1})$, note that   
any factor $\I(a_{i_0};a_{i_1},\dots,a_{i_r};a_{i_{r+1}})$ which occurs in the latter 
on the right of the tensor product can be written, using path composition and inversion formula, as a sum of products of 
terms  $\I(0;a_{i_1},\dots,a_{i_k}; a_{i_{k+1}})$. The factors $\I(0;b) =1$ are left out. 
We will show that after these reductions the coproducts are the same.

Let us establish a bijection of terms in (\ref{semicircle3}) with the terms in the coproduct for polygons.
We identify each semicircle polygon with an $R$-deco polygon. 
Let $P$ be an $(n+2)$-gon ($n\geq 1$) with 
 ``black'' vertices $a_0,\dots,a_{n+1}$ and a base side $a_0 a_{n+1}$, 
 inscribed in a semicircle. We put a ``white'' vertex, denoted by a $\circ$, 
on each of the $n+2$ arcs between the $a_i$'s. 
We order the points {\em counterclockwise}. 
The white vertices determine an $R$-deco $(n+2)$-gon $P^\vee$ with decorations 
given by the labels of the intermediate black vertices. From this, we obtain
an $R$-deco polygon $\pi$: 
the root arrow of $\pi$ emanates from the white vertex between the sides $a_0$ and $a_1$, and it ends
in the side $a_{n+1}$ (so that the side $a_0$ in $P^\vee$ is cut off).

\smallskip
Let us identify semicircle subpolygons with admissible dissections of $\pi$. 
Each subpolygon $P'$ of $P$ with $r+2$ sides, one of which must be  $a_0 a_{n+1}$, 
gives rise to a number of 
circumscribing $2(r+2)$-gons $Q$ with alternating black and white vertices. 
Both $P'$ and $Q$ inherit the orientation,  and $Q$ inherits
a direction of its edges $\alpha_1,\dots,\alpha_{2(r+2)}$ compatible with the orientation.

Such pairs $(P',Q)$ are in bijection with the terms in (\ref{semicircle2}). 
Precisely, the $2r+2$ directed sides $\alpha_i$ correspond to the factors on the right
of the tensor product.

Every second directed side goes from black to white vertex and thus does not correspond to the direction of arrows in an $R$-deco polygon
(going from vertex to side in an $R$-deco polygon corresponds to going from white vertex to black vertex in the polygon). Thus we invert
each such side $\alpha_i$ which becomes a backward arrow in the $R$-deco polygon. The sign which is assigned to such a backward arrow
in the bar construction for the $R$-deco polygon coincides with the sign which results from applying the inversion relation to 
the term for $\alpha_i$ in (\ref{semicircle3}).

Finally, we use
further relations for the integrals 
$\I(\dots)$: the directed sides $\alpha_i$ which connect a black vertex with an adjacent white one correspond to trivial arrows
in $\pi$ and can thus be dropped---in accordance with the fact that $\I(a;b)=1$. 
Furthermore, the leftover directed sides $\alpha_i$ which end in $a_0$ are set to zero---they correspond to terms of the form 
$\I(0;\dots;0)$ in (\ref{semicircle3}) which are set to zero if the sequence between the two zeros is not empty.

\vskip 3mm
Let us show how such a pair $(P',Q)$ arises from a term
$$
\I(a_0; a_{i_1},\dots,a_{i_r};a_{n+1})\ \otimes \ 
\prod_k \I(a_{i_k};a_{i_k+1},\dots, \dots ;0)\, \I(0;\dots ,\dots,a_{i_{k+1}-1}; 
a_{i_{k+1}})
$$
The left tensor factor of this term is pictured by a sub-(r+2)-gon $P'$ 
with vertices $a_0,a_{i_1},\dots,a_{i_r},a_{n+1}$, while the right hand factor, multiplied with the unit $\I(a_{n+1};0) \I(0;a_0)$,
is encoded by a sequence of $2(r+2)$ ``arrows'' 
which go from a black point $a_{i_k}$ to a white point $0$ or from a white point $0$ to a black point 
$a_{i_{k+1}}$. These arrows  form a polygon $Q$ with alternating black and white 
vertices. 
\qed

\begin{Example} 
The situation corresponding to the term 
$$\I(a_0; a_2,a_4;a_5)\otimes \I(a_0;a_1;a_2)\,\I(a_2;a_3;a_4)\,\I(a_4;a_5)$$
is illustrated by four pictures below.  

\vskip20pt
\xy
\hskip 90pt
\POS(36,0) *{a_{5}}
\POS(-34,0) *{a_{0}}
\POS(-20,-22.36) *{\bullet}
\POS(0,-30) *{\bullet}
\POS(-26,-14.97) *{\bullet}
\POS(23,-19.26) *{\bullet}
\POS(-23.5,-23.5) *{a_2}
\POS(-29.5,-15) *{a_1}
\POS(-2,-33) *{a_3}
\POS(27,-20) *{a_4}
\POS(0,30) *{\circ}
\POS(-28,-10.77) *{\circ}
\POS(10,-28.28) *{\circ}
\POS(-10,-28.28) *{\circ}
\POS(28,-10.77) *{\circ}
\POS(-23,-19.26) *{\circ}
\POS(-30,0) *{\bullet}; \POS(30,0)*{\bullet},
{\ellipse_{.}} 
\POS(30,0);\POS(-30,0), {\ellipse_{}} 
\POS(-30,0) \ar@{-} +{(10,-22.36)} \ar@{.} +{(60,0)}
\POS(-20,-22.36) \ar@{-} +{(43,3.1)}
\POS(23,-19.26) \ar@{-} +{(6.81,19.26)}
\POS(0,30) \ar@{-->>} +(30,-30)\ar@{-->>} +(-30,-30)
\POS(-28,-10.77) \ar@{-->>} +(-2,10.77)\ar@{-->>} +(8,-11.66)
\POS(10,-28.28) \ar@{-->>} +(-30,5.92)\ar@{-->>} +(13,8.49)
\POS(28,-10.77) \ar@{-->>} +(2,10.77)\ar@{-->>} +(-5,-8.49)
\endxy 
\vskip5pt
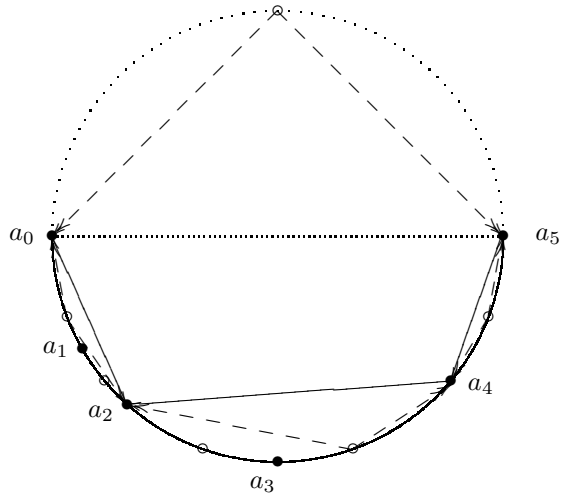
\centerline{\small An 8-gon circumscribing the 4-gon  $(a_0,a_2,a_4,a_5)$ inscribed in the semicircle}
\vskip10pt

The inscribed polygon with vertices at the black points labeled by $a_0,a_2,a_4$ and $a_5$, together with one of the terms arising from the path composition 
formula, gives rise to a dissection of the $R$-deco polygon, i.e., to a term in its coproduct. 

\vskip 20pt
\hskip 90pt
\xy
\POS(0,30) *{\circ} 
\POS(-28,-10.77) *{\circ} 
\POS(-23,-19.26) *{\circ} 
\POS(-10,-28.28) *{\circ}  
\POS(10,-28.28) *{\circ} 
\POS(28,-10.77) *{\circ}
\POS(-17,10) *{a_0}
\POS(17,10) *{a_5}
\POS(-28,-16) *{a_1}
\POS(-18,-26) *{a_2}
\POS(0,-31) *{a_3}
\POS(21,-21) *{a_4}
\POS(-28,-10.77)\ar@{-} +(5,-8.49) 
\POS(10,-28.28)\ar@{-} +(18,17.51)
\POS(-10,-28.28)\ar@{-} +(20,0)
\POS(-23,-19.26)\ar@{-} +(13,-9.02)
\POS(0,30)\ar@{-} +(-28,-40.77) \ar@{-} +(28,-40.77)
\POS(-28,-10.77) \ar@{=>>} +(45,15)

\endxy 
\vskip5pt
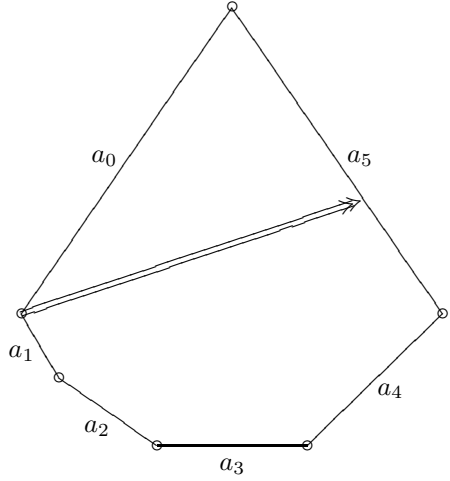
\centerline{\small The $R$-deco polygon associated to the 5-gon $(a_0,a_1,a_2,a_3,a_4,a_5)$}
\vskip10pt
\vskip 3mm

The second picture produces the third.

\vskip20pt
\xy
\hskip 90pt
\POS(0,30) *{\circ} 
\POS(-28,-10.77) *{\circ} 
\POS(-23,-19.26) *{\circ} 
\POS(-10,-28.28) *{\circ}  
\POS(10,-28.28) *{\circ} 
\POS(28,-10.77) *{\circ}
\POS(-24.5,-16.9) *{\bullet}
\POS(-15,-25.) *{\bullet}
\POS(0,-28.28) *{\bullet}
\POS(23,-15.9) *{\bullet}

\POS(-20.6,0) *{\bullet} \POS(20.6,0)*{\bullet},
\POS(-28,-10.77)\ar@{-} +(5,-8.49) 
\POS(10,-28.28)\ar@{-} +(18,17.51)
\POS(-10,-28.28)\ar@{-} +(20,0)
\POS(-23,-19.26)\ar@{-} +(13,-9.02)
\POS(0,30)\ar@{-} +(-28,-40.77) \ar@{-} +(28,-40.77)
\POS(-28,-10.77) \ar@{->>} +(12.9,-13.8)
\POS(-28,-10.77) \ar@{=>>} +(48,10.61)
\POS(10,-28.28)\ar@{->>} +(-24.8,3.92)
\POS(-28,-10.77) *+{}; \POS(-20.6,0)*+{} **\crv{(-20,-9)}
\POS(0,30) *+{}; \POS(-20.6,0)*+{} **\crv{(-10,6)}
\POS(28,-10.77) *+{}; \POS(20.6,0)*+{} **\crv{(20,-9)}
\POS(0,30) *+{}; \POS(20.6,0)*+{} **\crv{(10,6)}
\POS(23,-15.9);\POS(10,-28.28)**\crv{(16,-18)}
\POS(23,-15.9);(28,-10.77)**\crv{(25,-11)}
\POS(-17,10) *{a_0}
\POS(17,10) *{a_5}
\POS(-28,-16) *{a_1}
\POS(-18,-26) *{a_2}
\POS(0,-31) *{a_3}
\POS(21,-21) *{a_4}

\endxy 
\vskip 5pt
\centerline{\small 
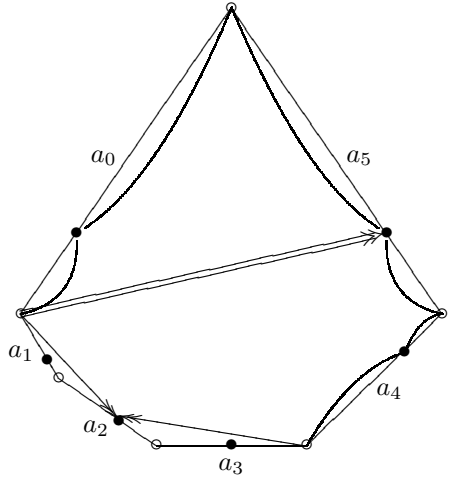
The picture after
pushing arrows inside.}
\vskip20pt

Dropping trivial arrows connecting adjacent vertices (necessarily of different color),
 we produce a picture of an $R$-deco polygon with a 3-fold dissection:
\vskip 20pt
\hskip 90pt
\xy
\POS(-28,-10.77) *{\circ} 
\POS(-23,-19.26) *{\circ} 
\POS(-10,-28.28) *{\circ}  
\POS(10,-28.28) *{\circ} 
\POS(28,-10.77) *{\circ}
\POS(25,-3) *{a_5}
\POS(-28,-16) *{a_1}
\POS(-18,-26) *{a_2}
\POS(0,-31) *{a_3}
\POS(21,-21) *{a_4}
\POS(-28,-10.77)\ar@{-} +(5,-8.49) 
\POS(-28,-10.77) \ar@{=>>} +(48,10.61)
\POS(10,-28.28)\ar@{-} +(18,17.51)
\POS(-10,-28.28)\ar@{-} +(20,0)
\POS(-23,-19.26)\ar@{-} +(13,-9.02)
\POS(20,0)\ar@{-} +(8,-10.77)
\POS(-28,-10.77) \ar@{->>} +(12.9,-13.8)
\POS(10,-28.28)\ar@{->>} +(-24.8,3.92)
\endxy 
\vskip 5pt
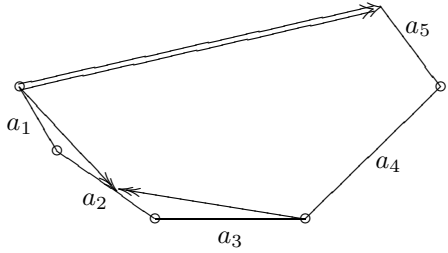
\centerline{\small A dissection associated to the term $\I(a_0;a_2,a_4;a_5)\otimes \I(0;a_1;a_2)
\I(0;a_3;a_2)$.}
\vskip20pt
The root arrow cuts off three trivial arrows and the extra side labeled by $a_0$.
\end{Example}

\subsection{Comparison with the coproduct from \cite{CK1}.}
In \cite{CK1}, Connes and Kreimer considered the 
coproduct on a Hopf algebra on rooted non-planar
trees, which has precisely the same form as our coproduct in Theorem \ref{MT}, 
except that we consider plane 
trees and have signs. Furthermore, they use the convention of 
writing the ``root part'' $R^C$ on the right of the
tensor product, while the remaining (cut-off) trees are denoted by $P_i^C$.
\begin{equation*}
\Delta_{{\rm CK}}(t)= \sum_{{\rm adm\ cut\ }C} \prod_i P_i^C(t) \otimes R^C(t)\,.
\end{equation*}
The commutative product in this case is simply the unordered disjoint union.

This connection
 between the coproduct formula in Theorem \ref{MT} and the one from \cite{CK1} is very similar to the 
relation between the coproduct in the Hopf algebra of iterated integrals in \cite{GonFund} 
and the one from \cite{CK1}.

\section{Associating integrals to the multiple logarithm cycles}
In this section we give another reason 
why the algebraic cycles associated to certain trees are ``avatars'' of 
multiple logarithms.  We first indicate how to associate, when $R=\Fx$ is the multiplicative group of a field $F$, 
an integral to the cycles $Z_{x_1,x_2}$ (cf.~(\ref{Z2})) and  
$Z_{x_1,x_2,x_3}$ (cf.~(\ref{Z3}))
which arise from the polygons $\lceilx x_1,x_2,1\rfloorx$ and  $\lceilx x_1,x_2,x_3,1\rfloorx$.
They deliver the integral
presentations for the double and triple logarithms. 
After these examples, we pass to the general case of 
any generic $R$-deco polygon $\lceilx x_1,x_2,\dots,x_m\rfloorx$ and associate to it an 
element in a certain bar construction, from which the Hodge realization can be read off.

\subsection{Brief review of the Hodge realization construction from \cite{BK}, \S8} \label{review}
Let ${\CN}'$ be a DG-subalgebra  of 
$\CN$.  In our case it is the image of the polygon algebra
$\CP$. Let ${\DP}$ be a larger DGA consisting of
{\it topological}
cycles in $\cub_{\C}^n$ which satisfy a number of conditions; we mention
the ones relevant for us, keeping the notation of loc.cit.:

i) For any cycle of degree $n$ in ${\DP}$, the integral against
 the standard volume form $\omega_n=\displaystyle{\bigwedge_{j=1}^n}(2\pi \sqrt{-1})^{-1}\displaystyle \frac{dz_j}{z_j}$ converges;
denote this integration map by $r: \DP \to \C$.
 
ii) $\CN'$ belongs to ${\DP}$; denote the embedding by $\sigma:\CN'\to \DP$.

iii) For any cycle $\sum_i n_i\lceilx a_1^i|\cdots |a_{k_i}^i\rfloorx$ in $\H^*\B({\cal N}')$,
the element  $\sum_i n_i \lceilx \sigma(a_1^i)|a_2^i|\cdots |a_{k_i}^i\rfloorx$ in 
$\H^*\B({\DP},{\cal N}')$ vanishes; this is a stronger version of the
claim that any cycle in ${\cal N}'$ becomes a boundary in ${\DP}$.

\smallskip
Under their assumptions Bloch and Kriz show that the map
$$\lambda :  \H^0\B({\DP,{\cal N}')\quad\stackrel{\B(r,\id)}{\longrightarrow}}\quad
\H^0\B(\C,{\cal N}')=\C\otimes \H^0\B( {\cal N}')$$
is an embedding. Using $\lambda$, they define a rational structure on 
 $\C\otimes \H^0\B( {\cal N}')$, and the Hodge and weight filtrations are
two natural filtrations corresponding to the Adams grading on this space.

As in \cite{BK} \S9, we construct topological {\em cycles} by the following trick.
Denote by $\delta$ the topological boundary. Let $\Gamma$ be
a small disk around $0$ in $\cub_{\C}$ with the standard orientation.
Consider an (oriented) topological {\it chain} $\eta$ of real dimension 
$2n-i$ in $\cub_{\Bbb C}^{2n-k}$ (in the notation of loc.cit. it is 
either $\eta_*(*)$
for $i=k$ or $q_*(*)$ for $i=k+1$). Associate to it a  {\it cycle} 
$$\theta_n(\eta)=\delta(\eta\times \Gamma)\times (\delta \Gamma)^{\times (i-1)}=
\delta\eta\times \Gamma\times (\delta \Gamma)^{\times (i-1)}
+(-1)^i\eta\times  (\delta \Gamma)^{\times i}
$$
of dimension $2n$ in $\cub_{\Bbb C}^{2n+i-k}$ (in the notation of \cite{BK}, (9.6), it corresponds to
the cycle $\tau_*(*)$ or $v_*(*)$ ). 

For the chains under consideration below we have the following: 
$$\partial(\theta_{n}(\eta))=\theta_{n-1}(\partial \eta)-
\theta_{n-1}(\delta\eta).$$
Indeed, $\delta\Gamma$ does not intersect the facets $0$ and $\infty$.
Therefore intersections of $\theta_n(\eta)$ with facets come 
from the intersections of $\delta(\eta\times \Gamma)$ with the respective facets.
The chain $\Gamma$ intersects the facet $0$.
For our chains $\eta$
the intersection of $\delta\eta$ with a facet coincides with the
topological boundary $\delta$ of the intersection of $\eta$ with the 
facet. Hence
$$\partial(\delta(\eta\times\Gamma))=\partial(\delta\eta\times \Gamma
+(-1)^i\eta\times \delta \Gamma)=
\partial(\delta\eta)\times \Gamma +(-1)^k\delta\eta
+(-1)^i\partial\eta\times \delta \Gamma$$
$$=
\delta\partial\eta\times \Gamma +(-1)^k\delta\eta
+(-1)^i\partial\eta\times \delta \Gamma =
\delta(\partial\eta\times \Gamma)+(-1)^k\delta\eta\,.$$
After multiplication by $(\delta\Gamma)^{\times(i-1)}$ we get
$$\partial(\delta(\eta\times\Gamma))\times(\delta\Gamma)^{\times(i-1)}=
\delta(\partial\eta\times \Gamma)\times(\delta\Gamma)^{\times(i-1)}
+(-1)^i\delta\eta\times\delta\Gamma\times (\delta\Gamma)^{\times(i-2)}$$
$$=
\delta(\partial\eta\times \Gamma)\times(\delta\Gamma)^{\times(i-1)}
+(-1)^{k-i+1}\delta(\delta\eta\times\Gamma)\times (\delta\Gamma)^{\times(i-2)}\,.
$$
In the special case $i=k$, where the cycle $\theta_n(\eta)$ is of real dimension $2n$ in $\cub_\C^{2n}$,
we can integrate it against the volume form $\omega_{2n}$, which gives---up to a factor $(-1)^i$---the 
integral of $\omega_{2n-i}$ against the chain $\eta$, as the restriction of $\omega_{2n-i}$ to $\delta(\eta)$ vanishes for
dimension reasons and the integral of $\frac{d z}{z}$ against $\delta\Gamma$ equals $2\pi \sqrt{-1}$.
 
\smallskip
Therefore the map $\theta_n$ assigns to each chain in $\DP$ a cycle, and 
we can reduce calculations to the case of topological chains
 $\eta$ equipped with two differentials: the algebraic one $\partial$ (intersection
 with facets) and the topological one $\delta$.

\subsection{Algebraic-topological chains}
Following Bloch and \Kriz, we embed the algebraic cycles into a larger 
set-up of ``topological'' chains which have both algebraic
and topological coordinates as well as both types of differentials, and then apply the bar construction. 
Those chains are referred to below as {\em algebraic-topological chains}.
We only consider ``topological'' variables $s_i\in [0,1]\subset \R$ 
subject to the condition $s_i\leq s_j$ if $i<j$,
and taking the topological boundary~$\delta$ for a chain with topological dimension $r$ 
amounts to taking the formal alternating sum over the 
subvarieties where either $s_k=s_{k+1}$ for some $k=1,\dots,r-1$ or $s_1=0$ or $s_r=1$.

\subsection{An example: the double logarithm case}\label{hodge2}
1. In order to bound $Z_{x_1,x_2}$, consider the algebraic-topological chain parametrized by 
$t\in \P_F^1$ and $s_1\in \R$,  $0\leq s_1\leq 1$, as
$$\Big[1-\frac {s_1}t,1-\frac t{x_1},1-\frac t{x_2}\Big]\,,$$
whose topological boundary terms are obtained by putting $s_1=0$ (which produces the empty cycle) or by $s_1=1$
which yields $Z_{x_1,x_2}$. Its algebraic boundary is given by 
$$[1- \frac{s_1}{x_1},1-\frac {s_1}{x_2}] - [1-\frac{s_1}{x_1},1-\frac {x_1}{x_2}] + [1-\frac{ s_1}{x_2},1-\frac{x_2}{x_1}]\,,\quad 0\leq s_1\leq 1\,,$$
where the last two terms are ``negligible'' for the following.

2. Consider the topological chain parametrized by $0\leq s_1\leq s_2\leq 1$, $s_i\in \R$, as
$$[1-\frac{s_1}{x_1}, 1-\frac{s_2}{x_2}]\,,$$
whose boundary terms arise from setting $s_1=0$, $s_1=s_2$ or $s_2=1$, giving the empty cycle, 
$[1-\frac{s_1}{x_1},1-\frac{s_1}{x_2}]$ or $[1-\frac{s_1}{x_1}, 1-\frac 1{x_2}]$, respectively.
What we are after is a cycle $\eta$ in this larger (algebraic-topological) chain complex which bounds the cycle $Z_{x_1,x_2}$,
i.e.,~such that $Z_{x_1,x_2}=(\partial+\delta)\eta$. The ``bounding process'' will give rise to a purely topological
cycle against which we can then integrate the standard volume form $\omega_2$.

In fact, working modulo the ``negligible'' term $[1-\frac{ s_1}{x_1},1-\frac{1}{x_2}]$ above we get 
\begin{equation}\label{boundCab}
Z_{x_1,x_2} = (\partial+\delta)\Big( \big[1- \frac{s_1}t,1-\frac t{x_1},1-\frac t{x_2}\big] + [1-\frac{s_1}{x_1},1-\frac{s_2}{x_2}]\Big)
\end{equation}
(two of the boundary terms cancel), and we associate to $Z_{x_1,x_2}$ the integral
\begin{multline} \label{zx1x2}
\frac 1{(2\pi \sqrt{-1})^2}\int\limits_{[1-\frac{s_1}{x_1},1-\frac{ s_2}{x_2}]\atop 0\leq s_1\leq s_2\leq 1}{} \frac {dz_1}{z_1}\wedge \frac{dz_2}{z_2} = 
\frac 1{(2\pi \sqrt{-1})^2}\int\limits_{0\leq s_1\leq s_2\leq 1} \frac {d(1-\frac{ s_1}{x_1})}{1-\frac{ s_1}{x_1}}
\wedge \frac{d(1-\frac{s_2}{x_2})}{1-\frac{ s_2}{x_2}} \\ = \frac 1{(2\pi \sqrt{-1})^2}\int\limits_{0\leq s_1\leq s_2\leq 1} \frac{ds_1}{s_1-x_1} \wedge \frac{ds_2}{s_2-x_2}\,. 
\end{multline}
Therefore we see that the algebraic cycle $Z_{x_1,x_2}$ corresponds---in a rather precise way---to the iterated
integral $I_{1,1}(x_1,x_2)$.
\sm

Note that the three ``negligible'' terms encountered above can also be covered as part of a boundary if we introduce,
following \cite{BK}, yet another differential $\dbar$ (coming from the bar construction), and in the ensuing
tricomplex all the terms above are taken care of. With the usual bar notation $\big|$ for a certain tensor product, 
the ``correct'' cycle combination is given by
\begin{multline} \label{doubhodge}
C_{1,1} := [1-\frac{s_1}{t},1-\frac t{x_1},1-\frac t{x_2}] + [1-\frac {s_1}{x_1},1-\frac{s_2}{x_2}] \\
+\ \Big(\,[1-\frac {s_1}{x_1}\Bar 1-\frac{1}{x_2}] - [1-\frac {s_1}{x_1}\Bar 1-\frac{x_1}{x_2}] + [1-\frac {s_1}{x_2}\Bar 1-\frac{x_2}{x_1}]\,\Big)\,, 
\end{multline}
and its image under the boundary $\partial + \delta+\dbar$ is precisely the ``bar version'' of $-Z_{x_1,x_2}$, i.e.
$\B(Z_{x_1,x_2})=-Z_{x_1,x_2} + \big(\,[1-\frac{1}{x_1}\Bar 1-\frac{1}{x_2}] - [1-\frac{1}{x_1}\Bar 1-\frac{x_1}{x_2}] + [1-\frac{1}{x_2}\Bar 1-\frac{x_2}{x_1}]\,\big)$. 

A similar treatment  of the triple logarithm can be found in \cite{GGL}.

\subsection{Enhanced polygons as a comodule over the polygon algebra} \label{comodule}
The above examples suggest to encode the algebraic-topological chains using trees
in a similar way as the algebraic cycles, except that we need to distinguish the topological variables $s_i$ above
(which are allowed to run through some real interval $[s_0,s_n]$) from the algebraic variables $t_i$ (which parametrize
$\P^1_F$). For this reason, we introduce ``enhanced trees'' with two types of vertices, the first type encoding
the algebraic variables and the second type encoding the topological ones.

\subsubsection{Enhanced trees} Let $\barbar R=R\cup \{s_0\}$ for some element $s_0$.
An {\bf enhanced $\barbar R$-deco tree} is an $\barbar R$-deco tree with two types of vertices.
All non-root vertices are {\bf of first type} except the ones which lie on the 
path from the first external vertex
to the root vertex. These vertices are {\bf of second type}. 
The root vertex has {\em both} types; this allows to identify
certain enhanced trees after contraction of the root edge with trees in $\CT(R)$.
External vertices of first type have decoration in $R$, while the unique external vertex of second type is
decorated by $s_0$. Typically, we choose $R=\Fx$ for a field $F$, and $s_0=0$.

\begin{Example} \label{enhancedtree}
An enhanced tree with six vertices of first type (marked by a bullet), three vertices of second type (marked
by a square) and the root vertex (marked by both a bullet and a square). The first external vertex
is decorated by $s_0$, while the root vertex is decorated by $s_3$. The path from $s_0$ to $s_3$ is drawn
with dotted lines.
\vskip 5pt
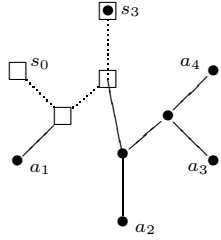
\xy
\hskip 120pt
\POS(-2,9) *{\Box} \ar@{.} +(-6,6)  *{\square} \ar@{-} +(-6,-6)  *{\bullet} \ar@{.} +(6,5)  *{\square}
\POS(12,9) *{\bullet} \ar@{-} +(6,-6)  *{\bullet} \ar@{-} +(6,6)  *{\bullet}\ar@{-} +(-6,-5)  *{\bullet}
\POS(6,4) \ar@{-} +(0,-9)  *{\bullet} \ar@{-} +(-2,10)
\POS(4,14) \ar@{.} +(0,9) *+{\bullet}*{\square}
\POS(9,-6) *{\scriptstyle a_2}
\POS(16,2) *{\scriptstyle a_3}
\POS(15,16) *{\scriptstyle a_4}
\POS(-5,2) *{\scriptstyle a_1}
\POS(-5,16) *{\scriptstyle s_0}
\POS(7,23) *{\scriptstyle s_3}
\endxy
\end{Example}

\begin{Definition} The differential on enhanced $R$-deco trees is defined similalry to the one on $R$-deco trees, 
except that if we contract an edge whose incident vertices
are of different type, then the resulting vertex, if it is not the root, is set to be of second type.
\end{Definition}
\def \gen{{\rm g}} 

In analogy with the DG-module $\DP$ over $\CN'$, we can define a DG-module $\barbar\CT_\bullet(\barbar R)$ over 
$\CT_\bullet(R)$, where $\barbar R= R\cup\{s_0\}$ for some element $s_0$.
Its generic part is denoted $\barbar\CT^\gen_\bullet(\barbar R)$, 
where genericity means that all decorations (including $s_0$) are different.

\subsubsection{The enhanced forest cycling map}
Let $F$ be a subfield of $\C$ and $s_0\in F\cap \R$. 
We pass from $\B\big(\barbar\CT_\smallbullet(F^\times \cup\{s_0\})\big)$ to $\B(\CD P,\CN^\smallbullet)$ as follows.
We modify the forest cycling map $\Phi$ from Definition \ref{fc} by introducing ``topological'' parameters $s_i$ 
(cf. \S\ref{hodge2}).

\begin{Definition} The {\bf enhanced} forest cycling map $\barbar\Phi$ is given as follows. Let $\tau$ be an enhanced $\Fx$-deco tree.
\begin{enumerate}
\item To each of the two external vertices of second type we associate elements $s_0$ and $s_{n+1}$ from 
 $F \cap \R$, while for the external vertices of first type we proceed as in 
Definition \ref{fc};
\item for each internal vertex of first type, we associate a variable in $\P^1_F$;
\item for each internal vertex of second type we associate a ``simplicial'' variable $s_i$ which runs through the interval
$[s_0,s_{n+1}]$; to each edge between two vertices of second type, oriented from $s_0$ to $s_{n+1}$, 
 we associate a ``$\leq$'' constraint;
\item to each oriented edge of $\tau$ from $v$ to $w$, decorated by $y_v$ and $y_w$, respectively, we associate
the expression $\big[1-\frac{y_v}{y_w}\big]$;
\item in the linear ordering of the edges, we concatenate all the coordinates in the previous step, except
for the edges connecting vertices of second type.
\end{enumerate}
The image of the tree with a single edge
 under $\barbar\Phi$ is the unique point in $\cub^0$.
\end{Definition}

\begin{Example} To the above example, we associate the chain 
$$\eta = \Big[1-\frac{s_1}{a_1}, 1-\frac{s_2}{t_1}, 1-\frac{t_1}{a_2}, 1-\frac{t_1}{t_2},1-\frac{t_2}{a_3},1-\frac{t_2}{a_4}\Big]$$
of real dimension 6 in $\cub_\C^6$, with $s_0\leq s_1\leq s_2\leq s_3$, and $t_j\in \P^1_\C$.
For the standard evaluation $s_0=0$ and $s_3=1$, the map $\theta$ assigns to $\eta$ the following cycle in $\cub^8_\C$
$$\theta_4(\eta) = \Big[1-\frac{s_1}{a_1}, 1-\frac{1}{t_1}, 1-\frac{t_1}{a_2}, 1-\frac{t_1}{t_2},1-\frac{t_2}{a_3},1-\frac{t_2}{a_4},r\,e^{i\vartheta_1},e^{i\vartheta_2}\Big] \ + \ \eta\,\times\, [e^{i\vartheta_1},e^{i\vartheta_2}]\,,$$
where $0\leq \vartheta_j <  2\pi$ and  $|r|<\eps$ for some small $\eps$. Here we parametrize the small disk $\Gamma$ 
of radius $\eps$ by polar coordinates.
\end{Example}

\medskip
The map $\barbar\Phi$ induces a morphism of DGAs.  The following Theorem  
is completely similar to Theorem \ref{dga_morphism} and thus its proof is omitted. We put $R=\Fx$ and $s_0=0$.

\begin{Theorem}\label{dga_morphism_hodge} For a subfield $F\hookrightarrow \C$, the enhanced forest cycling map $\barbar\Phi$, followed
by the map $\theta$ above, provides a natural map of differential graded algebras
$$\barbar\CT^\gen_\bullet(F) \to \DP\,.\qed$$
\end{Theorem}

\subsubsection{Enhanced polygons}
By duality, we are led to introduce ``enhanced polygons''.
Consider, e.g., the triangulated polygon dual to the enhanced tree in Example \ref{enhancedtree} above, where
the second type side is drawn by a dotted line:

\vskip10pt
\hskip 120pt
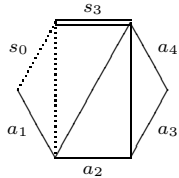
\xy 
\POS(0,0)  \ar@{-} +(10,0)_{a_2}  \ar@{-} +(-5,9)^{a_1}   \ar@{.} +(0,18)
\POS(10,0) \ar@{-} +(5,9)_{a_3}  \ar@{-} +(0,18) 
\POS(15,9) \ar@{-} +(-5,9)_{a_4} 
\POS(10,18) \ar@{=} +(-10,0)_{s_3} 
\POS(-5,9) \ar@{.} +(5,9)^{s_0}   
\POS(0,0) \ar@{-} +(10,18) 
\endxy
\vskip 10pt

Forgetting the triangulation we are led to:

\begin{Definition}
An {\bf enhanced $R$-deco polygon} $\barbar\pi$ is an $R$-deco polygon with one distinguished side $s$ of
the second type adjacent to the root side $r$. A linear order on the sides of $\barbar\pi$ is given by
starting at $s$ and ending at $r$,. It induces an orientation of $\barbar\pi$.
The decoration of $s$ and $r$ is in some ordered set $S$, and $s\leq r$.

The graded vector space $V^{\barbar\pg}$ is freely generated by enhanced polygons. 
The grading is the same as for polygons (i.e., the side of second type does not contribute).
\end{Definition}

In order to construct a differential and a resulting bar complex, we need the analogous notion of dissection 
for enhanced polygons.

A {\bf dissecting arrow} of an enhanced $R$-deco polygon $\barbar\pi$ is defined as for $R$-deco polygons, with two exceptions:
\begin{enumerate}
\item[i)] it is not allowed to start at the common vertex of $s$ and $r$;
\item[ii)] it is not allowed to end in $s$.
\end{enumerate}

We give two examples of non-allowed arrows, violating i) and  ii), respectively.
\vskip10pt
\hskip 60pt
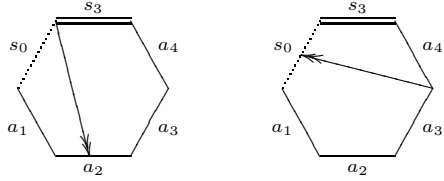
\xy 
\POS(0,0)  \ar@{-} +(10,0)_{a_2}  \ar@{-} +(-5,9)^{a_1}  
\POS(10,0) \ar@{-} +(5,9)_{a_3}  
\POS(15,9) \ar@{-} +(-5,9)_{a_4} 
\POS(10,18) \ar@{=} +(-10,0)_{s_3} 
\POS(-5,9) \ar@{.} +(5,9)^{s_0}  
\POS(0,18) \ar@{->>} +(4.5,-18)  

\hskip 100pt
\POS(0,0)  \ar@{-} +(10,0)_{a_2}  \ar@{-} +(-5,9)^{a_1}  
\POS(10,0) \ar@{-} +(5,9)_{a_3}  
\POS(15,9) \ar@{-} +(-5,9)_{a_4} \ar@{->>} +(-17.5,4.5)
\POS(10,18) \ar@{=} +(-10,0)_{s_3} 
\POS(-5,9) \ar@{.} +(5,9)^{s_0}  
\endxy
\vskip 10pt

An $n$-fold {\bf dissection} of $\barbar\pi$ is given by $n-1$ dissecting arrows (i.e., subject to conditions
i) and ii) above) which do not intersect.

\begin{Example} 
We draw a 4-fold dissection of an enhanced polygon 

\vskip10pt
\hskip 120pt
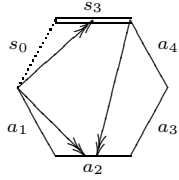
\xy 
\POS(0,0)  \ar@{-} +(10,0)_{a_2}  \ar@{-} +(-5,9)^{a_1}  
\POS(10,0) \ar@{-} +(5,9)_{a_3}  
\POS(15,9) \ar@{-} +(-5,9)_{a_4} 
\POS(10,18) \ar@{=} +(-10,0)_{s_3} 
\POS(-5,9) \ar@{.} +(5,9)^{s_0}   \ar@{->>} +(10,9)  \ar@{->>} +(9,-9)
\POS(10,18) \ar@{->>} +(-4.5,-18)  
\endxy

\end{Example}
\vskip 20pt
\def\bigwedgedot{{\bigwedge{}^{\hskip -3pt\bullet\hskip 3pt}}}

We can view $\barbar\CP :=\barbar\CP(\barbar R) :=V^\barbarpg(\barbar R) \otimes \bigwedgedot V^\pg(R)$ 
(this should not be confused with $\barbar\CP$ in \S\ref{tree_sums}) as a module over the algebra $\CP=\bigwedgedot V^\pg(R)$.
Similarly to the cobracket $\kappa$ of $\CP$, 
there is a comodule map 
$${\barbar\kappa} : V^\barbarpg(\barbar R) \to V^\barbarpg(\barbar R) \otimes \bigwedgedot V^\pg(R)$$
$$\barbar\kappa(\barbar\pi) = 
\sum_{\alpha\text{ arrow}} \sgn(\alpha)\,  \barbar\pi^\rt_\alpha\otimes \barbar\pi^\rst_\alpha\,,$$
where $\sgn(\alpha)$, $\barbar\pi^\rt_\alpha$ and $\barbar\pi^\rst_\alpha$ are as in Definition 
\ref{defposreg}.

The comodule map $\barbar\kappa$ induces a differential on $\barbar\CP$.

\begin{Definition}
There is a differential $\barbar\partial$ on $\barbar\CP$ given by 
$$\barbar\partial (\barbar\pi\otimes \pi') = \barbar\kappa(\barbar\pi)\wedge \pi'\ -\ \barbar\pi\otimes \kappa (\pi')\qquad (\barbar\pi\in V^\barbarpg(\barbar R),
\ \pi'\in \bigwedgedot V^\pg(R))\,.$$
\end{Definition}

\subsubsection{A cocycle in $\H^0\B(\barbar\CP,\CP)$}
The differential on $\barbar\CP$ gives rise, via a right bar resolution, to a right DG-comodule.

\begin{Definition}
The right $\B(\CP)$-comodule $\B(\barbar\CP, \CP)$ is generated as a vector space by the elements 
$$ [\barbar\pi_1\Bar \pi_2\Bar\dots\Bar\pi_n]\qquad (\barbar\pi_1\in \barbar\CP, \pi_i\in \CP).$$
It has a differential  $D_1+D_2$, where 
\begin{eqnarray*}
D_1(\,[\,\barbar\pi_1\Bar \pi_2\Bar\dots\Bar\pi_n\,]\,)&=& - \,[\,\barbar\pi_1\otimes \pi_2\Bar\dots\Bar\pi_n\,] \\
&&  + \sum_{j=2}^{n-1} (-1)^{j}   \,[\,\barbar\pi_1\Bar \pi_2\Bar\dots\Bar\pi_j\wedge \pi_{j+1}\Bar \dots\Bar\pi_n\,]\,, \\
D_2(\,[\,\barbar\pi_1\Bar \pi_2\Bar\dots\Bar\pi_n\,]\,)&=& \,[\,\barbar\kappa( \barbar\pi_1)\Bar \pi_2\Bar\dots\Bar\pi_n\,]+ \\
&& \sum_{j=2}^{n} 
(-1)^{j-1} \,[\,\barbar\pi_1\Bar \pi_2\Bar\dots\Bar \kappa (\pi_j)\Bar \dots\Bar\pi_n\,]\,.
\end{eqnarray*}
\end{Definition}

\medskip
As in  \S\ref{polygonalgebra}, each enhanced polygon provides a 0-cocycle in the bar complex $B(\barbar\CP,\CP)$.
\begin{Definition}
To an enhanced polygon $\barbar\pi$ we associate the following element, where  $\pi_D^i$ the associated subpolygons of a dissection $D$ of $\barbar\pi$ (we include the empty
dissection, corresponding to $\barbar\pi$ itself):
$$\B(\barbar\pi) = \sum_{\text{diss.\ $D$\ of\ }\barbar\pi}\sgn(D) \quad \sum_{\lambda} 
 \ [\barbar\pi_D^{\lambda(1)}\Bar \pi_D^{\lambda(2)}\Bar\dots\Bar\pi_{D}^{\lambda(|D|)}] \ 
\in \B(\barbar\CP,\CP)\,,$$
where the inner sum  runs through all linear orders $\lambda$ on the $\pi_D^i$ compatible with the partial
order on  the
dual tree $\tau(D)$.
\end{Definition}

\begin{Example} For the following enhanced 3-gon $\overline \pi$
\vskip 5pt
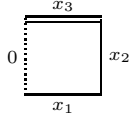
\xy
\hskip 140pt \hskip 20pt 
\POS(10,4) \ar@{=} +(-10,0)_{x_3}
\ar@{-} +(0,-10)^{x_2}
\POS(10,-6) \ar@{-} +(-10,0)^{x_1}
\POS(0,4) \ar@{.} +(0,-10)_0
\endxy
\vskip 5pt \noindent
we get $\B(\barbar\pi)$ as the following combination:

\vskip 5pt
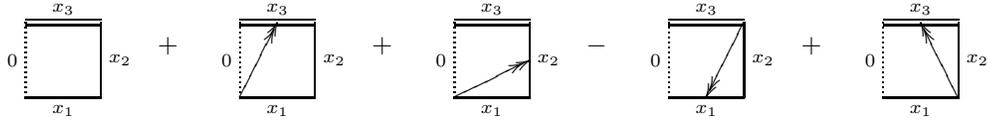
\xy
\POS(10,4) \ar@{=} +(-10,0)_{x_3}
\ar@{-} +(0,-10)^{x_2}
\POS(10,-6) \ar@{-} +(-10,0)^{x_1}
\POS(0,4) \ar@{.} +(0,-10)_0

\hskip 50pt $+$ \hskip 20pt 
\POS(10,4) \ar@{=} +(-10,0)_{x_3}
\ar@{-} +(0,-10)^{x_2}
\POS(10,-6) \ar@{-} +(-10,0)^{x_1}
\POS(0,4) \ar@{.} +(0,-10)_0
\POS(0,-6) \ar@{->>} +(5,10)

\hskip 50pt $+$ \hskip 20pt 
\POS(10,4) \ar@{=} +(-10,0)_{x_3}
\ar@{-} +(0,-10)^{x_2}
\POS(10,-6) \ar@{-} +(-10,0)^{x_1}
\POS(0,4) \ar@{.} +(0,-10)_0
\POS(0,-6) \ar@{->>} +(10,5)

\hskip 50pt $-$ \hskip 20pt 
\POS(10,4) \ar@{=} +(-10,0)_{x_3}
\ar@{-} +(0,-10)^{x_2}
\POS(10,-6) \ar@{-} +(-10,0)^{x_1}
\POS(0,4) \ar@{.} +(0,-10)_0
\POS(10,4) \ar@{->>} +(-5,-10)

\hskip 50pt $+$ \hskip 20pt 
\POS(10,4) \ar@{=} +(-10,0)_{x_3}
\ar@{-} +(0,-10)^{x_2}
\POS(10,-6) \ar@{-} +(-10,0)^{x_1}
\POS(0,4) \ar@{.} +(0,-10)_0
\POS(10,-6) \ar@{->>} +(-5,10)

\endxy
\vskip 5pt
They correspond, in this order, to the terms for $C_{1,1}$ in (\ref{doubhodge}).
\end{Example}

\begin{Proposition}\label{cocyc}
Any enhanced polygon $\barbar\pi$ provides a 0-cocycle, i.e.
$\B(\barbar\pi)$ gives a class in $\H^0 \B(\barbar\CP,\CP)$.
\end{Proposition}
\noindent{\bf Proof.} Analogous to Proposition \ref{cocycle}.\qed

\medskip
\subsubsection{Verifying the condition iii) from \S9.1}
Recall from \cite{BK}, (6.12), that there is a shuffle product on the circular bar construction $\B(\DP,\CN)$ defined by
sending 
$$(a'\otimes a_1\otimes \dots \otimes a_m)\otimes (a''\otimes a_{m+1}\otimes \dots \otimes a_{m+n})$$
to the sum of shuffle terms ($\sigma$ denots a shuffle of $\{1,\dots,m\}$ with $\{1,\dots,n\}$)
$$\pm(a'a''\otimes a_{\sigma(1)}\otimes \dots \otimes a_{\sigma(m+n)})\,.$$

\begin{Remark} \label{insert}
Let $\pi$ be an $R$-deco polygon and $\barbar\pi$ the associated enhanced polygon where the first vertex of $\pi$
has been replaced by a side of second type. 

\vskip 10pt
\hskip 50pt $\pi$ \hskip0pt  \threegon{a_1}{a_2}{a_3} \hskip 50pt $\barbar\pi $ 
 \hskip0pt  \fourgondotted{{s_0}}{a_1}{a_2}{a_3} 
\vskip 10pt
Then there is a bijection of the dissections of $\pi$ and those dissections of $\barbar\pi$ which have one arrow
going from the end of $s$ to the root side of $\barbar\pi$. 
Therefore $\B(\barbar\pi)$ breaks up into two parts,
one of them being $[\barbar 1\Bar \B(\pi)]$. But the bar differential on this expression is $\barbar1\wedge \B(\pi)$,
and since $\B(\barbar\pi)$ is a cocycle, we can conclude that $\B(\pi)$ (which is identified with 
$\barbar1\wedge \B(\pi)$) is the boundary of the other part of $\B(\barbar\pi)$. 

Using the shuffle on the circular bar construction, we can proceed similarly with terms of the form $[\barbar 1\Bar \Shaonly 
\B(\pi_i)]$, i.e., we can bound $\Shaonly \B(\pi_i)$. In this way we see that property (iii) of \S\ref{review} is satisfied for the cycles under consideration
below.
\end{Remark}

\subsubsection{The comodule map}
\medskip
We can now explicitly describe the comodule structure on the 0-cocycle representatives $\B(\barbar\pi)$, again with a Connes-Kreimer like
description. We need a notion of admissible cut for this setting.

\begin{Definition} Let $D$ be a dissection of an enhanced polygon $\barbar\pi$.
A subdissection $D'$ of $D$ is {\bf admissible} if it corresponds to a cut in the dual tree 
of $D$, which 
\begin{enumerate}
\item[i)] is admissible and
\item[ii)] avoids the edges of second type.
\end{enumerate}
\end{Definition}
\begin{Remark}
1. The conditions i) and ii) in the definition rule out the edges in $\tau(D)$ which correspond 
to the arrows of $D$ ending in $s$ or starting at $s\cap r$. These are precisely the ones which 
are ruled out in the conditions for dissecting arrows.

2. The second condition ensures that only the root polygon of a dissection has edges of second type
and therefore is the only polygon which produces an algebraic-to\-po\-lo\-gi\-cal chain in the associated element in the bar 
construction. 

3. There are two distinguished admissible dissections: the empty one and the one cutting off the side of second type.
\end{Remark}
The bar functor induces a comodule map $\B(\barbar\kappa)$ from $\B(\barbar\CP, \CP)$ to 
$\B(\barbar\CP, \CP) \otimes  \B(\CP)$.

\begin{Theorem} The image of $\B(\barbar\pi)$ under the comodule map $\B(\barbar\kappa)$  has the form
\begin{equation} \label{CKlike1}
\sum_{D\text{ \rm adm.\ diss\ of\ }\barbar\pi} \sgn(D)\,
\B(\barbar\pi^D_R) \otimes \Sha{ i\ }{\phantom{i=1}} 
\B(\pi^D_{P_i})\,,
\end{equation}
where $\barbar\pi^D_R$ and $\pi^D_{P_i}$ denote the root polygon and the remaining polygons arising from the dissection $D$.
\end{Theorem}

\noindent{\bf Proof.} Along the same lines as the proof of (\ref{CKlike}). \qed

\medskip
We can reinterpret the statement of the theorem as a $\SS^\bullet V^\pg$-comodule structure on 
$V^\barbarpg$: in view of Proposition \ref{cocycle} we have a map 
$$\SS^\bullet V^\pg \to \H^0\B(\CP)\,,$$
and similarly a map 
$$V^\barbarpg\otimes \SS^\bullet V^\pg \to \H^0\B(\barbar\CP,\CP)\,.$$

\subsection{The Hodge realization for polygons} \label{hodge}
We finally relate $\B(\barbar\CP,\CP)$ to the circular bar construction $\B(\CD\CP,\CN)$ of \cite{BK}.
Mapping $\barbar\CP \to \DP$ and $\CP\to \CN$ induces by functoriality a map
$$\H^0\B(\barbar\CP,\CP) \to \H^0\B(\DP,\CN)$$
The right hand term carries a mixed Hodge structure. It induces a one on the left.

The comodule structure (\ref{CKlike1}) carries over
to algebraic-topological cycles, via trees. The Hodge realization is obtained by integrating 
the first bar component of the left factor of each term in (\ref{CKlike1}) against the 
appropriate volume form on $\cub^\bullet$.

Note that the realization map $r$ (cf. \S\ref{review}) is non-trivial only 
on those trees whose internal vertices
are all of second type. For the associated triangulated polygons,
this corresponds to the unique triangulation for which every triangle is incident with the 
end of the root arrow, see the picture below. 
\vskip 5pt
\hskip 120pt
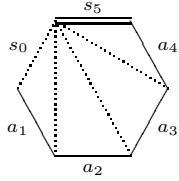
\xy 
\POS(0,0)  \ar@{-} +(10,0)_{a_2}  \ar@{-} +(-5,9)^{a_1}   \ar@{.} +(0,18)
\POS(10,0) \ar@{-} +(5,9)_{a_3}  \ar@{.} +(-10,18) 
\POS(15,9) \ar@{-} +(-5,9)_{a_4} 
\POS(10,18) \ar@{=} +(-10,0)_{s_5} 
\POS(-5,9) \ar@{.} +(5,9)^{s_0}   
\POS(0,18) \ar@{.} +(15,-9) 
\endxy

\vskip 5pt
The integration map $\lambda$ from \S\ref{review} evaluates each enhanced polygon in the coproduct
expression (\ref{CKlike1}) for a given element, and the
Hodge weight can be read off from the number of sides of the polygon. Note that $\lambda$ 
applies only to the first entry $\barbar\pi_1$ of any component $[\barbar\pi_1\Bar\pi_2\Bar\dots\Bar \pi_n]$ of
the bar construction.

\smallskip
The Adams grading defines two filtrations on $V^\pg$ (recall that $\deg(\pi)$ is the number of
non-root sides of $\pi$):
\begin{eqnarray*}
W_{-2(N-1)}(V^\pg) &=& \bigoplus_{\deg(\pi)\geq N} \Q\,\pi \qquad \text{(the weight filtration),}\\
F^{-(N-1)} (V^\pg\otimes \C) &=&  \bigoplus_{\deg(\pi)\leq N} \C\,\pi \qquad\text{(the Hodge filtration)}\,. 
\end{eqnarray*}
They induce filtrations on $\SS^\bullet V^\pg$.

We now define a Hodge-Tate structure on $\SS^\smallbullet V^\pg$ as follows: there is an obvious map embedding
$V^\pg$ into $V^\barbarpg$ by simply adding a side of second type, with a fixed decoration $s_0$, to a polygon (as in Remark \ref{insert} above). Further, there is a map
$$V^\barbarpg \otimes \SS^\smallbullet V^\pg \to \H^0\B(\DP,\CN)$$
where the right hand term is an $\H^0\B(\CN)$-module and the left hand term is an $\SS^\smallbullet V^\pg$-module.
More precisely, we can embed $\SS^\smallbullet V^\pg \hookrightarrow V^\barbarpg \otimes \SS^\smallbullet V^\pg$ by the map
$$ \prod_i \pi_i \mapsto \sum_i \barbar\pi_i \prod_{j\not= i} \pi_j\,.$$
Composing the above with the integration map $r:V^\barbarpg\to \C$ now provides a map
$$\SS^\smallbullet V^\pg \to \C\otimes \SS^\smallbullet V^\pg\,,$$
which gives us the desired rational structure on $\SS^\smallbullet V^\pg$.

\begin{Example} We illustrate the above for the example of the double logarithm. The topological cycle which we need
to consider is $C_{1,1}-[\barbar 1\Bar \B(Z_{x_1,x_2})]$, where $C_{1,1}$ is as in (\ref{doubhodge}) and $Z_{x_1,x_2}$
as in (\ref{Z2}).
Applying $\lambda\circ\theta$ to it annihilates the first term in $C_{1,1}$ (for reasons of type), while the second 
(purely topological) term gives 
$$I  \otimes 1  \ \in \C\otimes \H^0\B(\CP) $$
where $I$ denotes the iterated integral in (\ref{zx1x2}). The three further terms in $C_{1,1}$
give us, up to a factor $(2\pi \sqrt{-1})^{-1}$, the combination
$$ \int_{[1-\frac {s_1}{x_1}]}\frac{dz}{z} \otimes \Big[ 1-\frac{1}{x_2}\Big] - 
\int_{[1-\frac {s_1}{x_1}]}\frac{dz}{z}\otimes \Big[1-\frac{x_2}{x_1}\Big] 
+ \int_{[1-\frac {s_1}{x_2}]}\frac{dz}{z}\otimes \Big[1-\frac{x_1}{x_2}\Big]\,,$$
and finally integrating $[\barbar 1\Bar \B(Z_{x_1,x_2})]$ yields
$1\otimes \B(Z_{x_1,x_2})\,.$
\end{Example}

\section{Algebraic cycles for multiple polylogarithms}

In this section, we sketch how to produce algebraic cycles for multiple {\em poly}\-loga\-rithms,
or rather the associated iterated integrals $I_{n_1,\dots,n_r}(a_1,\dots,a_r)$, in a similar 
fashion using $\Fx$-deco polygons, where we have to allow a third type of side which is undecorated. Roughly, the new type 
corresponds to factors $\frac{dz}z$ 
in the integral representation of polylogarithms. Accordingly, we need to modify the forest cycling map to form the associated 
admissible algebraic cycles. The following condition are required: the first and the last side (i.e., the sides which are
incident with the first vertex) are of the first type.

The iterated integral $I_{n_1,\dots,n_r}(a_1,\dots,a_r)$ corresponds to the polygon with the following sequence of sides (we
denote the undecorated sides by an $\emptyset$):
$$\pi_{n_1,\dots,n_r}(a_1,\dots,a_r):=\lceilx a_1,\underbrace{\emptyset,\dots,\,\emptyset}_{n_1-1\text{ times}},\dots, a_r,\underbrace{\emptyset,\dots,\,\emptyset}_{n_r-1\text{ times}},1\rfloorx\,.$$

The forest cycling map in this setting is modified as follows: let $\tau$ be the dual tree for a triangulation of 
a polygon $\pi$ with both types of sides, then we only need to address the undecorated external edges of $\tau$.
In this case, the associated coordinate of the algebraic cycle is simply given by the parametrizing variable $t$ itself 
(as opposed to, e.g., $1-\frac{t}{y}$ for some edge with external decoration $y$).

If we want to associate an {\em admissible} algebraic cycle, we need to mod out by the ideal
generated by  triangulations containing 
triangles with more than one undecorated side of the original polygon.

\end{document}